\documentclass[pdflatex,sn-mathphys-num]{sn-jnl}

\usepackage{graphicx}%
\usepackage{multirow}%
\usepackage{amsmath,amssymb,amsfonts}%
\usepackage{amsthm}%
\usepackage{mathrsfs}%
\usepackage[title]{appendix}%
\usepackage{xcolor}%
\usepackage{textcomp}%
\usepackage{manyfoot}%
\usepackage{booktabs}%
\usepackage{algorithm}%
\usepackage{algorithmicx}%
\usepackage{algpseudocode}%
\usepackage{listings}%
\usepackage{fullpage}
\usepackage{mathtools}
\usepackage{float}
\usepackage{hyperref}
\usepackage{tikz}
\usepackage{tkz-graph}
\usetikzlibrary{arrows,shapes,patterns,positioning,calc}
\usetikzlibrary{arrows.meta}
\usepackage{subcaption}
\usepackage{siunitx}
\usepackage{pgfplots}

\usepackage{pifont}

\theoremstyle{thmstyleone}%
\newtheorem{theorem}{Theorem}
\theoremstyle{thmstyletwo}%
\newtheorem{example}{Example}%
\newtheorem{remark}{Remark}%

\theoremstyle{thmstylethree}%
\newtheorem{definition}{Definition}%

\begin{document}

\title[Article Title]{Exploring Logistic Functions as Robust Alternatives to Hill Functions in Genetic Network Modeling}


\author*[1]{\fnm{Ismail} \sur{Belgacem}}\email{ismail.belgacem.81@gmail.com}
\affil*[1]{\orgdiv{Independent Researcher, Mezaourou}, \orgname{Ghazaouet}, \city{Tlemcen}, \postcode{13421}, \country{Algeria}}


\abstract{
Hill functions have dominated gene regulatory network (GRN) modeling for decades, yet their fractional exponents create severe analytical obstacles when Hill coefficients assume non-integer values, a ubiquitous occurrence in experimental parameter fitting. We propose systematically replacing Hill functions with logistic functions: the increasing Hill function $h^+(x, \theta, n) = x^n/(x^n + \theta^n)$ with $f^+(x, \theta, \lambda) = 1/(1 + e^{-\lambda(x-\theta)})$, and the decreasing Hill function $h^-(x, \theta, n) = \theta^n/(x^n + \theta^n)$ with $f^-(x, \theta, \lambda) = 1/(1 + e^{\lambda(x-\theta)})$, where the parameter relationship $\lambda = n/\theta$ matches the slopes of the two functions at the half-maximal concentration $x = \theta$, thereby preserving local input-output sensitivity while retaining the biological interpretability of cooperativity. We provide comprehensive mathematical foundations by systematically cataloguing logistic function properties: global $C^\infty$ differentiability, closed-form derivatives with self-referential structure ($\partial f^+/\partial x = \lambda f^+(1-f^+)$), elementary inverse functions via logit transformation, analytically tractable integrals, and multiple approximation regimes. In contrast, we rigorously document Hill function limitations for non-integer $n$: derivative singularities at the origin ($h^{+\prime}(x) \to \infty$ as $x \to 0^+$ for $0 < n < 1$, with higher-order derivatives diverging for $n \in (k,k+1)$), integrals requiring hypergeometric functions, multivalued complex branches in inversion, and logarithmic small-$n$ approximations diverging at low expression. These pathologies propagate into severely ill-conditioned Fisher information matrices when datasets include low-concentration measurements, compromise numerical integration stability near basal expression regimes, and cause observability matrix rank deficiency in control applications. We establish rigorous theoretical foundations for the product-of-logistics GRN model by proving global existence, smoothness, and boundedness of solutions, with an explicit Lipschitz bound $L_F \le M = \max_i(\kappa_i \sum_j L_i^j + \gamma_i)$ guaranteeing well-conditioned Jacobians throughout biologically relevant domains. We demonstrate that Hill and logistic functions are mathematically related by a logarithmic coordinate change---the Hill function is the logistic evaluated at $\ln(x/\theta)$ rather than $x-\theta$---but that this relation does \emph{not} make the two ODE models equivalent: they differ in basal expression, smoothness, Lipschitz constant, and sensitivity structure. The logistic function encodes additive threshold detection (equal absolute increments of $x$ near $\theta$ produce equal response changes), whereas the Hill function encodes multiplicative-increment sensitivity (equal fold-changes in $x/\theta$ produce equal response changes). Here ``fold-change'' refers to the Weber--Fechner sensitivity structure of the function itself; this is distinct from ``fold-change in gene expression'' as used in thermodynamic models of transcriptional regulation~\cite{bintu2005transcriptional,bintu2005transcriptional_app}, where fold-change denotes the ratio of promoter activity in the presence versus absence of a transcription factor. Logistic functions naturally model basal expression ($f^+(0, \theta, \lambda) = 1/(1 + e^{\lambda\theta}) > 0$) without ad hoc modifications, directly capturing the persistent low-level transcription observed experimentally. The uniformly bounded derivative ($|\partial f^\pm/\partial x| \leq \lambda/4$) and global Lipschitz continuity provide substantial computational advantages for gradient-based optimization, parameter estimation with well-conditioned Fisher information matrices, and observer-based state estimation that maintains full-rank observability at low concentrations where Hill derivatives vanish.
}


\maketitle


\section{Introduction}
\label{sec:introduction}

Gene regulatory networks (GRNs) orchestrate cellular decision-making through complex interactions among transcription factors, regulatory elements, and target genes. Mathematical modeling of these networks requires sigmoidal functions that capture cooperative binding, threshold-like regulatory transitions, and saturating responses. For decades, Hill functions have served as the standard formulation, with the increasing form $h^+(x,\theta,n) = x^n/(x^n + \theta^n)$ modeling activation and the decreasing form $h^-(x,\theta,n) = \theta^n/(x^n + \theta^n)$ modeling repression. Their mechanistic foundation in equilibrium binding theory and their ability to encode cooperativity through the Hill coefficient $n$ have established them as the \textit{de facto} choice in systems biology, from bacterial gene circuits to mammalian signaling pathways~\cite{belgacem2018reduction,belgacem2014stability,belgacem2013stability,belgacem2013analysis,belgacem2014mathematical,belgacem2020control,chambon2020qualitative,bernot2012modeling,belgacem2013global,belgacem2012full,kim2011stochastic,polynikis2009comparing}.

However, Hill functions impose subtle yet consequential analytical burdens when the Hill coefficient $n$ assumes non-integer values, a ubiquitous occurrence in experimental parameter fitting~\cite{gottschalk2005five,hernandez2023corrected,reeve2013pharmacodynamic,santillan2008use}. When fitted to empirical dose-response data, Hill coefficients frequently yield non-integer values reflecting incomplete cooperativity, heterogeneous binding-site occupancy, or complex allosteric mechanisms. These fractional exponents create mathematical pathologies that compound in high-dimensional systems: the first derivative diverges at the origin for $0 < n < 1$ (creating a vertical tangent); the second derivative becomes unbounded for $1 < n < 2$; and, more generally, derivatives beyond order $\lfloor n \rfloor$ fail to exist at $x = 0$ whenever $n \in (k, k+1)$. These singularities compromise Jacobian-based stability analysis, prevent rigorous Taylor expansions around low-expression equilibria, and cause numerical integration failures precisely where many biological systems operate---near basal expression levels~\cite{santillan2008use}.

Analytical obstacles extend well beyond differentiation. Hill functions lack closed-form antiderivatives for general $n$, requiring hypergeometric functions or incomplete beta functions~\cite{abramowitz1965handbook}. Inversion requires fractional roots, which introduce multivalued complex branches. Parameter estimation encounters a further difficulty: the gradients $\partial h^+/\partial n$ and $\partial h^+/\partial \theta$ exhibit logarithmic and power-law singularities at low concentrations, propagating into severely ill-conditioned Fisher information matrices when the dataset includes near-zero expression measurements. The resulting ill-conditioning yields unreliable confidence intervals and strong parameter correlations, compromising identifiability. Control-theoretic applications encounter rank-deficient observability matrices when Hill function derivatives vanish at low concentrations, preventing state reconstruction in extended Kalman filters and high-gain observers. These limitations compound significantly in genome-scale models, which require efficient simulation of thousands of regulatory interactions, and in synthetic biology applications that demand real-time model predictive control.

Several alternative formulations have been explored. Step functions provide mathematical simplicity but are undefined at threshold values, precluding rigorous solution theory across regulatory switches~\cite{abou2011theoretical,casey2006piecewise}. Piecewise-linear approximations maintain analytical tractability for qualitative analysis~\cite{de2004qualitative,gouze2002class,belgacem2025chaos,farcot2019chaos}, but they inherit discontinuities that require sophisticated switching frameworks~\cite{belgacem2019probabilistic} and fail to capture the gradual, concentration-dependent responses observed empirically. Ramp function approximations introduce additional threshold parameters without biological justification for the assumed linearity within transition regions~\cite{bottani2017hill,polynikis2009comparing,quee2021ramp}. Critically, none of these alternatives resolves the fundamental tension between biological fidelity and analytical tractability.

We propose a systematic resolution: replace Hill functions with logistic functions that preserve sigmoidal regulatory characteristics while providing superior mathematical structure. Specifically, we substitute the increasing Hill function with $f^+(x,\theta,\lambda) = 1/(1 + e^{-\lambda(x-\theta)})$ and the decreasing Hill function with $f^-(x,\theta,\lambda) = 1/(1 + e^{\lambda(x-\theta)})$. The steepness parameter $\lambda$ is matched to preserve local input-output sensitivity through $\lambda = n/\theta$, ensuring identical slopes at the half-maximal concentration $x = \theta$ while retaining the biological interpretability of cooperativity.

Although logistic functions have demonstrated versatility across statistics, machine learning, and network engineering~\cite{berkson1953statistically,goebbels2020sharpness,ke2014logistic,koizumi2010adaptive,puddu2009artificial}, their application to GRN modeling has been limited. Samuilik et al.~\cite{samuilik2022mathematical} introduced logistic representations but employed a single increasing form with signed weights encoding regulatory direction, a compromise that obscures the distinct dynamics of activation versus repression; subsequent studies have perpetuated this unified framework~\cite{kozlovska2024search,kozlovska2025modeling,sadyrbaev2021modelling,samuilik2022genetic,sadyrbaev2023coexistence,somathilaka2023revealing}. Our approach instead deploys mathematically distinct increasing and decreasing logistic functions wherever biologically appropriate, preserving the interpretability that is essential for understanding regulatory architectures.

This product-of-logistics formulation differs from the weighted-sum framework in four concrete ways. First, regulator-specific thresholds are independently determined from single-regulator dose-response curves, and each parameter maps directly to a measurable biological quantity, whereas the single-sigmoid approach of Samuilik et al.~obscures this through shared thresholds. Second, the explicit decreasing logistic naturally approaches unity at zero repressor concentration, capturing maximal unrepressed expression without artificial offsets. Third, activation and repression are encoded through the functional form rather than the sign of a weight, eliminating a structural pathology in which a negative weight inside an increasing sigmoid places the repressor's critical point at $x_j = \theta_i/w_{ij} < 0$, outside the biologically admissible domain. Fourth, threshold and steepness can be tuned separately without compensatory adjustments---a parameter independence that is essential for systematic control design and absent from weighted-sum formulations, which exhibit ridge-like likelihood surfaces that compromise identifiability.

Our contributions are as follows. \textbf{First}, we systematically catalogue the analytical limitations of Hill functions for non-integer $n$ and trace each pathology to a specific failure mode in numerical integration, stability analysis, parameter estimation, observer design, and control synthesis. \textbf{Second}, we establish the corresponding mathematical properties of the increasing and decreasing logistic functions and show that each Hill failure mode is resolved by a precisely identifiable logistic property: the bounded derivative replaces the divergent power-law slope, the closed-form logit replaces the multivalued fractional-power inverse, the elementary antiderivative replaces hypergeometric integration, and the nonzero basal output $f^+(0,\theta,\lambda) = 1/(1+e^{\lambda\theta})$ replaces the absorbing Hill state $h^+(0)=0$. \textbf{Third}, we prove (Theorem~\ref{thm:smoothness}) that the resulting product-of-logistics GRN model has globally unique, smooth, and uniformly bounded solutions, with an explicit and computable Lipschitz constant $L_F \le M = \max_i(\kappa_i \sum_j L_i^j + \gamma_i)$. \textbf{Fourth}, we clarify the precise mathematical relationship between the two function families: $h^+(x,\theta,n) = \sigma(n\ln(x/\theta))$ is an exact pointwise identity, but the corresponding ODE models are nonequivalent, because the change of variable $s = \ln(x/\theta)$ rescales the production term by a state-dependent factor $e^{-s}$. The two formulations therefore encode genuinely different biological hypotheses---multiplicative-increment versus additive-threshold sensitivity---and the structural advantages of the logistic formulation hold regardless of which interpretation a given system requires.

The paper proceeds as follows. Section~\ref{sec:hill_limitations} establishes the mathematical limitations of Hill functions for non-integer $n$: derivative singularities at the origin, intractable integration, inversion through multivalued complex branches, and logarithmically divergent approximations. Section~\ref{sec:logistic_derivation} derives the increasing and decreasing logistic forms and the parameter-matching relationship $\lambda = n/\theta$. Section~\ref{sec:logistic_foundations} catalogues the mathematical foundations of the logistic functions: $C^\infty$ differentiability, self-referential derivatives, closed-form integrals, analytically invertible logit forms, and global Lipschitz bounds. Section~\ref{sec:equivalence_and_scaling} establishes the coordinate-change identity $h^+(x,\theta,n) = \sigma(n\ln(x/\theta))$, explains why it does not yield ODE-level equivalence, and clarifies the structural differences between the two function families. Section~\ref{sec:practical_consequences} traces the Hill-function limitations of Section~\ref{sec:hill_limitations} into practical failure modes for numerical integration, stability analysis, parameter estimation, observer design, and control synthesis (including LQR, feedback linearization, SMC, MPC, and EKF), identifying in each case the logistic property that resolves the failure. Section~\ref{sec:biological_realism} examines biological realism, including the basal-expression advantage and the sensitivity-structure distinction. Section~\ref{sec:advanced_analysis} presents the general product-of-logistics GRN model, its Boolean-logic extension via De~Morgan's law, the comparison with weighted-sum formulations, and Theorem~\ref{thm:smoothness} on global existence, smoothness, and boundedness. Section~\ref{sec:conclusion} concludes.

\section{Mathematical Limitations of Hill Functions}
\label{sec:hill_limitations}

Hill functions, conventionally formulated as $h^+(x, \theta, n) = \frac{x^n}{\theta^n + x^n}$ for gene activation, emerge from mechanistic equilibrium analysis of ligand-receptor binding kinetics. Related global-dynamics analyses of concave gene-expression models and reduced transcription--translation systems are developed in~\cite{belgacem2014mathematical,belgacem2013stability,belgacem2013analysis,belgacem2014stability,belgacem2018reduction}. Within GRN frameworks, Hill functions encode cooperative binding effects through the Hill coefficient $n > 0$, which simultaneously quantifies the degree of molecular cooperativity and the steepness of the sigmoidal regulatory response. The dissociation constant $\theta > 0$ represents the inducer or ligand concentration at which the response reaches half-maximal, functioning as a biological threshold \cite{santillan2008use}.

The Hill coefficient $n$ merits careful examination. Operationally, experimentalists measure promoter activity over a range of inducer concentrations, fit sigmoidal curves to experimental data, and extract the Hill coefficient as a quantitative measure of response steepness and local input-output sensitivity. This data-driven estimation captures the empirical response without requiring detailed mechanistic assumptions.  Nevertheless, $n$ often admits a mechanistic interpretation: under cooperative binding conditions, where multiple transcription factors or ligands interact or where allosteric effects amplify the response, $n$ approximates the number of interacting binding sites or the degree of cooperativity. The \textit{lac} operon illustrates both the utility and the limits of this picture. Statistical-mechanical (thermodynamic) models of transcriptional regulation~\cite{bintu2005transcriptional,bintu2005transcriptional_app} show that simple repression by LacI at a single operator is described by the expression $F_{\mathrm{reg}} = 1/(1+[R]/K_R)$---a Hill function with $n=1$, not $n > 1$. The dramatic $>$1000-fold repression observed in the wild-type system arises instead from DNA looping between the main operator and auxiliary operators, which is captured by a distinct rational-polynomial expression (Case 9 of Table 1 in~\cite{bintu2005transcriptional}) involving an effective looping concentration $[L]\approx 660$\,nM. The Hill coefficient $n\approx 2$--$3$ commonly cited for the \textit{lac} system refers to the IPTG-induction response---how the inducer releases LacI from DNA---rather than to the repressor--DNA binding thermodynamics. More generally, any Hill coefficient with $n > 1$ obtained by curve-fitting to induction data is a phenomenological summary of the overall network response; it does not directly correspond to a mechanistic step in the equilibrium binding model. Non-integer $n$, which is the norm in experimental fits, has no mechanistic interpretation in terms of integer binding-site multiplicity. Nevertheless, the conceptual framework of cooperativity---even when parameterized phenomenologically---directly informs synthetic biology applications. Engineered toggle switches and oscillators in \textit{E. coli} and yeast rely on precisely tuned cooperative binding (typically $n = 2$--$4$) to achieve bistable memory devices, programmable cellular timers, and biosensors with sharp dose-response thresholds \cite{gardner2000construction,elowitz2000synthetic}. Across naturally occurring and engineered gene regulatory networks, Hill coefficients typically cluster in the range $n = 1$ to $4$, capturing realistic molecular interactions without invoking extreme ultrasensitivity. This moderate range proves instructive. Hemoglobin's cooperative oxygen binding ($n \approx 2.8$) achieves sufficient steepness for efficient oxygen loading in the lungs and unloading in tissues, yet remains gradual enough to prevent catastrophic threshold behavior that would destabilize oxygen homeostasis \cite{santillan2008use}.  Similarly, transcriptional regulatory networks across biology employ Hill coefficients in this range ($n = 1$--$4$), where dimer or tetramer formation of transcription factors generates robust, tunable responses \cite{santillan2008use}. Rare exceptions exist: certain signaling cascades exhibit $n > 4$, producing sharp, switch-like responses. However, these exceptional cases incur costs, including heightened computational instability and reduced biological plausibility in most regulatory contexts. The Hill coefficient need not be an integer. Experimental fitting often yields non-integer values when cooperative mechanisms are incomplete or heterogeneous, in which case $n$ functions as a purely phenomenological parameter extracted from data rather than reflecting a clear molecular mechanism.  Furthermore, Hill coefficients can take values between $0$ and $1$, typically indicating \textit{negative cooperativity}, where the binding of one ligand or transcription factor hinders subsequent bindings, leading to a shallower response curve compared to non-cooperative systems ($n = 1$) \cite{abeliovich2005empirical}. While less common in gene regulatory networks, such sub-unity coefficients are observed in certain receptor systems and enzymatic reactions where conformational changes reduce the binding affinity for subsequent ligands.

The Hill framework has proven invaluable across systems biology, from bacterial gene circuits to mammalian signaling pathways. However, beneath this apparent success lies a fundamental mathematical pathology that has been inadequately addressed in the literature. When Hill coefficients assume non-integer values, a ubiquitous occurrence in experimental data fitting, several mathematical pathologies emerge that compromise both analytical tractability and numerical stability. These limitations are not merely theoretical curiosities but have profound practical consequences for simulation, stability analysis, parameter estimation, and control design.

\subsection{Differentiability Pathologies at the Origin}

Mathematical tractability depends critically on the nature of the Hill coefficient. When $n$ takes positive integer values, Hill functions maintain global infinite differentiability ($C^\infty$) across the entire domain, including at boundary points such as the origin $x = 0$, enabling smooth analytical manipulations, Taylor expansions, and stable numerical implementations. Yet this analytical grace vanishes in practical parameter estimation from real biological data.

When Hill models are fitted to experimental measurements of gene expression versus inducer concentration, the resulting Hill coefficients frequently assume non-integer values \cite{hernandez2023corrected,santillan2008use,reeve2013pharmacodynamic,gottschalk2005five}, a ubiquitous phenomenon reflecting incomplete cooperativity, heterogeneous binding site occupancy, or complex allosteric mechanisms that resist simple integer-valued interpretations. This transition from integer to non-integer Hill coefficients precipitates a catastrophic loss of mathematical structure, as differentiability collapses specifically at the origin, $x = 0$, generating singularities that render the system mathematically ill-posed and computationally unstable.

\subsubsection{Asymptotic behavior near the origin}
To understand the fundamental limitations, consider the standard activation form $h^+(x, \theta, n) = \frac{x^{n}}{\theta^{n} + x^{n}}$. Near $x = 0$, the substitution $t = x^n / \theta^n$ reveals that:
\[
h^+(x, \theta, n) = \frac{t}{1 + t} = t - t^2 + t^3 - \cdots,
\]
yielding leading behavior proportional to $x^n$ for small $x$. The derivatives exhibit increasingly severe singularities as we approach the origin. The first derivative
\[
h^{+\prime}(x) = \frac{n \theta^{n} x^{n-1}}{(\theta^{n} + x^{n})^2}
\]
behaves asymptotically as $h^{+\prime}(x) \sim \frac{n}{\theta^{n}} x^{n-1}$ as $x \to 0^+$. Similarly, differentiating $h^{+\prime}(x)$ with respect to $x$ (using the quotient rule or the asymptotic form) yields the second derivative $h^{+\prime\prime}(x) \sim \frac{n(n-1)}{\theta^{n}} x^{n-2}$ near the origin. More generally, writing $h^{+(m)}(x) := d^m h^+/dx^m$, the $m$-th derivative is proportional to $x^{n-m}$ as $x \to 0^+$:
\[
h^{+(m)}(x) \propto x^{n-m}, \qquad x \to 0^+,
\]
with prefactor given by the falling factorial (Pochhammer symbol):
\[
n(n-1) \cdots (n-m+1) = \frac{\Gamma(n+1)}{\Gamma(n-m+1)}.
\]

\subsubsection{Derivation of the asymptotic form of higher derivatives}
To make this precise, we rewrite $h^+(x) = x^n(\theta^n + x^n)^{-1}$. By Leibniz's rule for the $m$-th derivative of a product $u(x) = x^n$ and $v(x) = (\theta^n + x^n)^{-1}$:
\[
h^{+(m)}(x) = \sum_{k=0}^{m}\binom{m}{k}u^{(k)}(x)v^{(m-k)}(x).
\]
As $x \to 0^+$, the dominant contribution comes from $k = m$ (the term $u^{(m)}(x)v^{(0)}(x) \sim n(n-1)\cdots(n-m+1)\,\theta^{-n}\,x^{n-m}$ has leading order $x^{n-m}$, while every other term with $k < m$ has order $x^{(n-k) + (n-1)(m-k)} = x^{n-m + n(m-k)}$, smaller by a factor of $x^{n(m-k)}$), giving:
\[
h^{+(m)}(x) \sim \frac{\Gamma(n+1)}{\Gamma(n-m+1)}\,\theta^{-n}\,x^{n-m}, \qquad x \to 0^+.
\]

\subsubsection{Case analysis for non-integer Hill coefficients}
These asymptotic forms create severe problems when $n$ is non-integer, a ubiquitous situation in experimental parameter estimation:

\begin{itemize}
    \item \textbf{For $0 < n < 1$:} The first derivative diverges, $h^{+\prime}(x) \sim \frac{n}{\theta^n} x^{n-1} \to \infty$ as $x \to 0^+$, creating a vertical tangent or cusp-like singularity at the origin. The function is not differentiable at $x = 0$. This infinite slope creates insurmountable numerical difficulties for optimization algorithms, ODE solvers, and observer design methods that rely on gradient information. In control applications, the unbounded derivative prevents reliable linearization around low-excitation equilibria, undermining the synthesis of feedback controllers.
    
    \item \textbf{For $1 < n < 2$:} While $h^{+\prime}(0) = 0$ exists and is finite (since $n - 1 > 0$ gives $x^{n-1} \to 0$), the second derivative explodes: $h^{+\prime\prime}(x) \sim \frac{n(n-1)}{\theta^n} x^{n-2} \to \infty$ as $x \to 0^+$ (since $n - 2 < 0$), indicating unbounded curvature. This limits the function to $C^1$ (once continuously differentiable) at the origin, preventing the use of Newton-Raphson methods, Hessian-based optimization, and second-order sensitivity analysis.
    
    \item \textbf{For $n \in (k, k+1)$ where $k \ge 1$:} Derivatives of order $m > \lfloor n \rfloor$ diverge at the origin. The Hill function is only $C^{\lfloor n \rfloor}$ at $x = 0$, not $C^\infty$. This finite-differentiability class limits the applicability of advanced analytical techniques that require arbitrary-order derivatives, such as center manifold theory, normal form analysis, and high-order perturbation methods.
\end{itemize}

Away from the origin on any domain $(\delta, \infty)$ with $\delta > 0$, the function remains analytic as a rational composition of smooth terms. However, the origin-specific pathologies persist and cannot be circumvented by restricting the domain, as biological systems frequently operate in or pass through low-expression regimes. While the function itself is well-defined and continuous at $x = 0$ with $h^+(0) = 0$, the higher-order smoothness failure stems directly from the fractional power-law terms, introducing singularities that cascade through derivatives.

Consider concrete cases: when $n \in (0,1)$, we get an infinite slope at zero; when $n \in (1,2)$, the slope vanishes, but curvature explodes; and for $n \in (k, k+1)$ with $k \ge 1$, all derivatives beyond the $k$-th diverge at the origin. These pathologies become particularly problematic in cooperative regulatory systems, where empirical Hill coefficients commonly fall within the range $n \approx 1.5$ to $n \approx 3.5$, creating second- and third-order derivative singularities that destabilize numerical methods.

\subsubsection{Implications for dynamical analysis}
These undefined higher derivatives create serious complications for dynamical analysis:

\begin{itemize}
    \item They interfere with Jacobian calculations at equilibrium points, essential for determining local stability through eigenvalue analysis.
    \item They obscure the bifurcation structure by preventing the construction of normal forms that require smooth coordinate changes.
    \item They invalidate Taylor series expansions around zero, eliminating a primary tool for local approximation.
    \item They complicate Lyapunov exponent computations, which require derivatives of the flow map.
    \item In observer design, they cause the observability matrix, constructed from successive Lie derivatives of the output function, to lose rank at the origin for $n > 1$, rendering hidden states unobservable precisely in the low-expression regimes where robust state estimation is most crucial.
\end{itemize}

The resulting analytical blind spots can mask important dynamics, such as oscillations, multistability, noise-induced transitions, and time-delay-induced behaviors, that would otherwise be accessible through standard techniques.

\subsection{Integration Challenges}

In stark contrast to the logistic function's elementary antiderivative, the Hill function lacks closed-form antiderivatives for general $n$. Related global-stability analyses for enzymatic chains and Michaelis--Menten reactions appear in~\cite{belgacem2013global,belgacem2012full}. The integral
\[
\int \frac{\xi^n}{\theta^n + \xi^n} \, d\xi
\]
poses significant analytical challenges that escalate rapidly with increasing complexity of the Hill coefficient $n$.

\subsubsection{Case $n=1$: Michaelis-Menten kinetics}
For the simplest non-trivial case, we write $\frac{\xi}{\theta+\xi} = 1 - \frac{\theta}{\theta+\xi}$:
\begin{align}
\int \frac{\xi}{\theta + \xi} \, d\xi &= \int \left(1 - \frac{\theta}{\theta + \xi}\right) d\xi = \xi - \theta \ln|\theta + \xi| + C.
\end{align}

\subsubsection{Case $n=2$: Cooperative binding}
We write $\frac{\xi^2}{\theta^2+\xi^2} = 1 - \frac{\theta^2}{\theta^2+\xi^2}$:
\begin{align}
\int \frac{\xi^2}{\theta^2 + \xi^2} \, d\xi &= \int \left(1 - \frac{\theta^2}{\theta^2 + \xi^2}\right) d\xi = \xi - \theta^2\int\frac{d\xi}{\theta^2 + \xi^2}.
\end{align}
Using the standard formula $\int \frac{d\xi}{a^2 + \xi^2} = \frac{1}{a}\arctan\!\left(\frac{\xi}{a}\right) + C$ with $a = \theta$:
\begin{equation}
\int \frac{\xi^2}{\theta^2 + \xi^2} \, d\xi = \xi - \theta\arctan\!\left(\frac{\xi}{\theta}\right) + C.
\end{equation}
This expression already involves inverse trigonometric functions, significantly complicating symbolic integration in dynamical systems analysis.

\subsubsection{General case: Series representation}
For arbitrary $n$, we employ the substitution $u = \xi/\theta$, yielding $d\xi = \theta\,du$:
\[
\int \frac{\xi^n}{\theta^n + \xi^n} \, d\xi = \theta \int \frac{u^n}{1 + u^n} \, du = \theta\int\!\left(1 - \frac{1}{1+u^n}\right)du.
\]
For $|u^n| < 1$ (corresponding to $|\xi| < \theta$), we apply the geometric series $\frac{1}{1+u^n} = \sum_{k=0}^{\infty}(-1)^k u^{nk}$, which converges uniformly on compact subsets of $\{u : |u| < 1\}$, and integrate term by term:
\begin{align}
\theta\int \frac{1}{1 + u^n} \, du &= \theta\sum_{k=0}^{\infty} (-1)^k \frac{u^{nk+1}}{nk+1}.
\end{align}
Since $\frac{u^n}{1+u^n} = 1 - \frac{1}{1+u^n}$, integrating gives
\[
\theta\int\frac{u^n}{1+u^n}\,du
= \theta u - \theta\sum_{k=0}^{\infty}(-1)^k\frac{u^{nk+1}}{nk+1}.
\]
The $k=0$ term of the series equals $\theta u^1/(1) = \theta u$, which cancels the standing $\theta u$:
\[
= \theta u - \theta u + \theta\sum_{k=1}^{\infty}(-1)^{k+1}\frac{u^{nk+1}}{nk+1}
= \theta\sum_{k=1}^{\infty}(-1)^{k+1}\frac{u^{nk+1}}{nk+1}.
\]
Substituting back $u = \xi/\theta$:
\begin{equation}
\int \frac{\xi^n}{\theta^n + \xi^n} \, d\xi = \sum_{k=1}^{\infty} (-1)^{k+1} \frac{\xi^{kn+1}}{\theta^{kn}(kn+1)} + C, \qquad |\xi| < \theta.
\end{equation}

\subsubsection{Connection to special functions}
For non-integer $n$, this series cannot be summed in terms of elementary functions. Instead, the antiderivative requires hypergeometric functions $_2F_1$ or incomplete beta functions $B_z(a,b)$. Through the substitution $v = u^n$ (so $u = v^{1/n}$, $du = \frac{1}{n}v^{1/n-1}\,dv$):
\[
\int \frac{u^n}{1+u^n} \, du = \frac{1}{n}\int \frac{v}{1+v}\,v^{1/n-1}\,dv = \frac{1}{n}\int\frac{v^{1/n}}{1+v}\,dv,
\]
which evaluates, via the integral representation of the hypergeometric function, to:
\[
\int \frac{u^n}{1+u^n} \, du = \frac{u^{n+1}}{n+1}\,{}_2F_1\!\left(1,\,1+\tfrac{1}{n};\,2+\tfrac{1}{n};\,-u^n\right) + C
\]
after further manipulation \cite{abramowitz1965handbook}.

\subsubsection{Computational implications}
The complexity of these integrals has several profound consequences:
\begin{itemize}
    \item \textbf{Computationally expensive:} Numerical evaluation requires adaptive quadrature methods (e.g., Gauss-Kronrod rules) with careful error control. For non-integer $n$, each evaluation may require iterative computation of special functions, multiplying cost by factors of 10--100 compared to elementary operations.
    
    \item \textbf{Not amenable to symbolic manipulation:} The lack of closed forms prevents the derivation of analytical formulas for equilibrium points, basin boundaries in multistable systems, or bifurcation curves where parameter-dependent integrals determine stability transitions.
    
    \item \textbf{Limited qualitative insights:} Without explicit closed forms, the dependence on parameters $(\theta, n)$ remains opaque, preventing intuitive understanding and asymptotic analysis.
    
    \item \textbf{Numerical instability near $\xi = 0$:} The derivative singularities propagate into integration routines. When $0 < n < 1$, the integrand's derivative diverges as $\xi \to 0^+$, causing adaptive integrators to refine step sizes excessively.
    
    \item \textbf{Series convergence issues:} The power series converges only for $|\xi| < \theta$. For $\xi > \theta$ (supra-threshold concentrations), alternative representations are mandatory. Convergence deteriorates as $\xi \to \theta^-$, requiring $O(1/\epsilon)$ terms for accuracy $\epsilon$.
\end{itemize}

\subsection{Inversion Problems}

The increasing Hill function $h^+(x) = \frac{x^n}{\theta^n + x^n}$ lacks a general closed-form inverse for arbitrary $n$. To find the formal inverse, we solve $y = x^n/(\theta^n + x^n)$ for $x$:
\begin{align}
y(\theta^n + x^n) &= x^n \notag\\
y\theta^n &= x^n(1 - y) \notag\\
x^n &= \frac{y}{1-y}\,\theta^n \notag\\
x &= \theta\left(\frac{y}{1-y}\right)^{1/n}.
\end{align}
For integer values, the inverse can be written explicitly. For example, when $n=2$:
\[
(h^+)^{-1}(y) = \theta\sqrt{\frac{y}{1-y}},
\]
which still involves an irrational (square root) function. For $n = 3$, $(h^+)^{-1}(y) = \theta(y/(1-y))^{1/3}$, involving cube roots.

For non-integer $n$ (e.g., $n \approx 1.39$ from empirical fits \cite{santillan2008use}), the formal inverse
\[
(h^+)^{-1}(y) = \theta \left( \frac{y}{1-y} \right)^{1/n}
\]
requires computing fractional powers $(\cdot)^{1/n}$ for non-integer $1/n$. This is a multi-valued complex function for non-integer exponents: for a real argument $r > 0$, we have $r^{1/n} = e^{(\ln r)/n}$, which is well-defined and real; however, for $r < 0$ (which can occur along complex paths), $r^{1/n}$ introduces branch cuts and numerical ambiguities. While $y/(1-y) > 0$ for $y \in (0,1)$, this necessitates numerical inversion, introducing approximation errors that hinder exact linearization and limit the precision of real-time control.

Iterative methods such as Newton-Raphson can compute the inverse numerically; however, this approach introduces computational overhead and potential convergence failures, particularly near $y = 0$ and $y = 1$. As $y \to 1^-$ (equivalently, $x \to \infty$), $h^{+\prime}(x) \sim n\theta^n/x^{n+1} \to 0$ for every $n > 0$, making each Newton step diverge in magnitude regardless of the cooperativity. As $y \to 0^+$ (equivalently, $x \to 0^+$), $h^{+\prime}(x) \to 0$ for $n > 1$ (Newton steps blow up, again) while $h^{+\prime}(x) \to \infty$ for $0 < n < 1$ (Newton steps collapse to zero, stalling progress). Both regimes therefore obstruct reliable inversion, with the obstruction at $y \to 0^+$ depending qualitatively on the cooperativity range.

\subsection{Approximation Difficulties}

We explicitly derive the Hill function's small-$n$ expansion. Write $x^n = e^{n\ln x}$ and $\theta^n = e^{n\ln\theta}$, so:
\[
h^+(x, \theta, n) = \frac{e^{n\ln x}}{e^{n\ln x} + e^{n\ln\theta}} = \frac{1}{1 + e^{n(\ln\theta - \ln x)}} = \frac{1}{1 + e^{-n\ln(x/\theta)}}.
\]
This is the logistic function in $z = \ln(x/\theta)$ with steepness $n$, evaluated at the single point. Taylor-expanding in $n$ around $n = 0$: at $n = 0$,
\[
h^+(x,\theta,0) = \frac{1}{1 + e^0} = \frac{1}{2}.
\]
The derivative with respect to $n$ at $n = 0$ is:
\[
\left.\frac{\partial h^+}{\partial n}\right|_{n=0} = \left.\frac{-e^{-n\ln(x/\theta)}\cdot(-\ln(x/\theta))}{(1 + e^{-n\ln(x/\theta)})^2}\right|_{n=0} = \frac{\ln(x/\theta)}{4}.
\]
Therefore:
\begin{equation}\label{eq:hill_small_n}
h^+(x, \theta, n) \approx \frac{1}{2} + \frac{n}{4}\ln\!\left(\frac{x}{\theta}\right) + O(n^2),
\end{equation}
which introduces logarithmic dependence on $x$. Since $\ln(x/\theta) \to -\infty$ as $x \to 0^+$, this approximation is ill-defined near the origin, rendering it numerically unstable precisely where many biological systems operate under basal conditions.

\subsubsection{Behavior near the origin}
Around $x = 0$, the Hill function exhibits power-law behavior:
\[
h^+(x, \theta, n) = \frac{x^n}{\theta^n + x^n} = \frac{x^n}{\theta^n}\cdot\frac{1}{1 + (x/\theta)^n} = \frac{x^n}{\theta^n}\left[1 - \left(\frac{x}{\theta}\right)^n + O\!\left(\left(\frac{x}{\theta}\right)^{2n}\right)\right] \approx \frac{x^n}{\theta^n} + O\!\left(\left(\frac{x}{\theta}\right)^{2n}\right),
\]
indicating that the function vanishes at the origin with increasing flatness for larger $n$. The behavior of the first derivative at $x = 0$ depends critically on $n$:
\begin{itemize}
\item For $n > 1$: $h^{+\prime}(x) \sim \frac{n}{\theta^n}x^{n-1} \to 0$ as $x \to 0^+$, so $h^{+\prime}(0) = 0$.
\item For $n = 1$: $h^{+\prime}(x) \to 1/\theta$ as $x \to 0^+$, giving a non-zero finite slope.
\item For $0 < n < 1$: $h^{+\prime}(x) \sim \frac{n}{\theta^n}x^{n-1} \to \infty$ as $x \to 0^+$, leading to an infinite slope.
\end{itemize}

\subsubsection{Behavior around the threshold}
Around the midpoint $x = \theta$, we expand $h^+(x,\theta,n)$ by writing $x = \theta + \epsilon$ for small $\epsilon$:
\begin{align}
h^+(\theta+\epsilon) &= \frac{(\theta+\epsilon)^n}{(\theta+\epsilon)^n + \theta^n} = \frac{1}{1 + (\theta/(\theta+\epsilon))^n} = \frac{1}{1 + (1 + \epsilon/\theta)^{-n}} \notag\\
&\approx \frac{1}{1 + 1 - \frac{n\epsilon}{\theta} + O(\epsilon^2/\theta^2)} \approx \frac{1}{2} + \frac{n}{4\theta}\epsilon + O(\epsilon^2),
\end{align}
giving:
\begin{equation}
h^+(x, \theta, n) \approx \frac{1}{2} + \frac{n}{4\theta}(x - \theta) + O\!\left((x-\theta)^2\right).
\end{equation}
Comparing with the logistic expansion $f^+(x) \approx \frac{1}{2} + \frac{\lambda}{4}(x-\theta)$, equating slopes yields the parameter matching condition $\lambda = n/\theta$. The Hill function's slope $n/(4\theta)$ depends on both $n$ and $\theta$, making it sensitive to the threshold scale in a way that the logistic function's slope $\lambda/4$ is not (since $\lambda$ is an independent parameter). These approximations become functionally equivalent when parameters are matched through $\lambda = n/\theta$, ensuring identical slopes at the midpoint and enabling similar local stability analysis in dynamical models.

\section{Derivation of the Equivalent Logistic Forms for Gene Regulatory Network Modeling}
\label{sec:logistic_derivation}

This section derives explicit increasing and decreasing logistic representations for activation and repression that preserve biological interpretability while ensuring analytical tractability. The parameter-matching conditions linking logistic steepness to Hill cooperativity are established along the way.

\medskip
\noindent\textbf{Construction of the increasing form.} For activation dynamics, we begin with the increasing standard logistic function, as detailed in \cite{chen2013properties,goebbels2020sharpness,kyurkchiev2015sigmoid}: 
\[
f(s) = \frac{1}{1 + e^{-s}}.
\]
 This canonical form is centered at $s = 0$, where the function attains the value $1/2$, and exhibits unit steepness in the exponential argument. To adapt this function for gene regulation, we must introduce two modifications:
\begin{enumerate}
    \item \textbf{Adjustable steepness:} We incorporate a parameter $\lambda > 0$ to control the transition sharpness, yielding
    \[
    f(s, \lambda) = \frac{1}{1 + e^{-\lambda s}}.
    \]
    Larger values of $\lambda$ produce steeper sigmoids, analogous to higher Hill coefficients $n$ in cooperative binding.
    
    \item \textbf{Threshold alignment:} Since biological regulation responds to absolute concentration levels rather than deviations from zero, we shift the inflection point to the regulatory threshold $\theta$ by substituting $s = x - \theta$. This translation ensures the function is centered at the biologically meaningful concentration $x = \theta$.
\end{enumerate}

Therefore, for activation dynamics, we define $f^+(x, \theta, \lambda) = f(x - \theta, \lambda)$, where the argument $s = x - \theta$ ensures the desired monotonicity: as the regulator concentration $x$ increases, so does $s$, and consequently $f^+(x, \theta, \lambda)$ increases strictly monotonically from values near zero (when $x \ll \theta$) to values approaching unity (when $x \gg \theta$).

\medskip
\noindent\textbf{Construction of the decreasing form.} For repression, we require a function that decreases with increasing regulator concentration. We achieve this by reversing the sign in the exponential. Specifically, we define the \textbf{standard decreasing logistic function} by:
\[
f(s) = \frac{1}{1 + e^{s}},
\]
As before, we introduce steepness through $\lambda$ and shift the inflection point to $\theta$, yielding:

\[
f^-(x, \theta, \lambda) = \frac{1}{1 + e^{\lambda (x - \theta)}}.
\]
The positive coefficient in the exponential ensures strict monotonic decrease: as $x$ rises above $\theta$, the denominator grows exponentially, driving the function value toward zero. Equivalently, we may write this as
\[
f^-(x, \theta, \lambda) = \frac{1}{1 + e^{-\lambda (\theta - x)}},
\]
making explicit that repression responds to the signed deviation $\theta - x$ rather than $x - \theta$.

\medskip
\noindent In summary: The \textbf{activation} is modeled by the increasing Hill function $h^+(x, \theta, n) = \frac{x^n}{x^n + \theta^n}$, now replaced by the increasing logistic function $f^+(x, \theta, \lambda) = \frac{1}{1 + e^{-\lambda (x - \theta)}}$. The \textbf{repression} is modeled by the decreasing Hill function $h^-(x, \theta, n) = \frac{\theta^n}{x^n + \theta^n}$, now replaced by the decreasing logistic function $f^-(x, \theta, \lambda) = \frac{1}{1 + e^{\lambda (x - \theta)}} = \frac{1}{1 + e^{-\lambda (\theta - x)}}$.  In both cases, the parameters $n$ and $\lambda$ govern the steepness of the sigmoidal response, quantifying how sharply the regulatory function transitions between its lower and upper asymptotes. The threshold parameter $\theta$ marks the concentration at which the regulatory effect is half-maximal, serving as the inflection point of the sigmoid curve.

\section{Increasing and Decreasing Logistic Functions: Mathematical Foundations}
\label{sec:logistic_foundations}

Logistic and Hill functions both generate smooth sigmoidal response curves with tunable steepness and threshold behavior. We present the mathematical properties of the increasing and decreasing logistic functions here to establish the toolkit deployed in Section~\ref{sec:practical_consequences} (where the resolution of each Hill-function failure mode is traced to a specific logistic property) and in Section~\ref{sec:advanced_analysis} (where these properties are lifted to the network level). Because the decreasing logistic $f^-$ satisfies $f^-(x,\theta,\lambda) = 1 - f^+(x,\theta,\lambda)$ (the symmetry property derived in Section~\ref{sec:logistic_symmetry} below), its mathematical properties mirror those of $f^+$; detailed derivations are given for the increasing case, with the decreasing case following by analogy.

\begin{figure}[htbp]
  \centering
  \includegraphics[width=\linewidth]{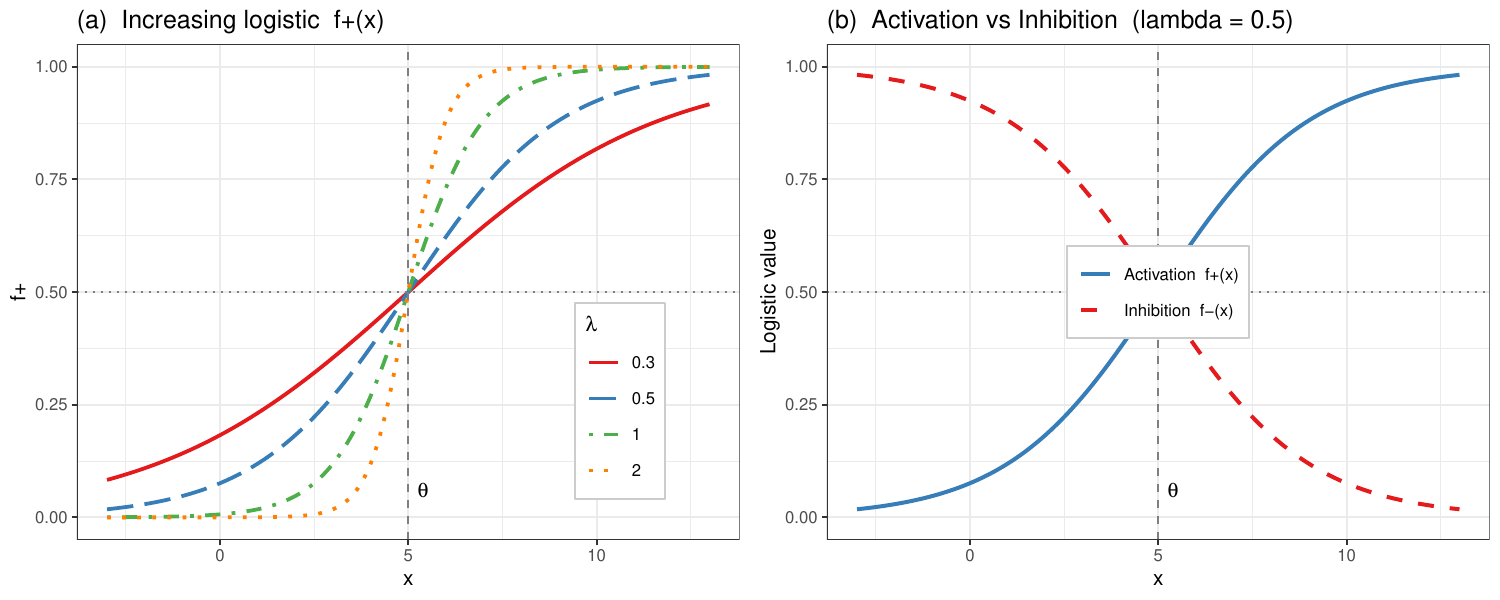}
  \caption{Logistic functions used in the GRN modeling framework.
    \textbf{(a)} Increasing logistic $f^+(x,\theta,\lambda)$ for
    $\theta=5$ and four values of steepness
    $\lambda\in\{0.3,\,0.5,\,1.0,\,2.0\}$; higher $\lambda$ yields a
    sharper threshold response.
    \textbf{(b)} Comparison of the increasing (activation) logistic
    $f^+(x,\theta,\lambda)$ and the decreasing (repression) logistic
    $f^-(x,\theta,\lambda)=1-f^+(x,\theta,\lambda)$ for $\theta=5$,
    $\lambda=0.5$.  The vertical dashed line marks the inflection point
    $x=\theta$ where both functions attain their half-maximum value
    $\tfrac{1}{2}$ and their maximum absolute slope $\lambda/4$
    (equation~\eqref{eq:maxslope}).}
  \label{fig:logistic}
\end{figure}

We recall the two logistic primitives that underpin the entire framework and illustrate their behavior in Figure~\ref{fig:logistic} before cataloguing the analytical properties in the subsections below.

\begin{definition}[Increasing logistic --- activation]
\label{def:fplus}
\begin{equation}
  f^+(x,\theta,\lambda)
  \;=\; \frac{1}{1+e^{-\lambda(x-\theta)}}.
  \label{eq:fplus_def}
\end{equation}
This function increases monotonically from $f^+\approx 0$ for $x\ll\theta$ to $f^+\approx 1$ for $x\gg\theta$,
with inflection point $f^+(\theta,\theta,\lambda)=\tfrac{1}{2}$.  It models the activation of a gene (or
process) by a regulator $x$: the gene is essentially off below the threshold $\theta$ and fully on above it,
with sharpness controlled by $\lambda>0$.
\end{definition}

\begin{definition}[Decreasing logistic --- repression]
\label{def:fminus}
\begin{equation}
  f^-(x,\theta,\lambda)
  \;=\; \frac{1}{1+e^{\lambda(x-\theta)}}
  \;=\; 1-f^+(x,\theta,\lambda).
  \label{eq:fminus_def}
\end{equation}
This function decreases monotonically from $f^-\approx 1$ for $x\ll\theta$ to $f^-\approx 0$ for
$x\gg\theta$, implementing Boolean NOT: high regulator concentration suppresses expression.
\end{definition}

Both functions are $C^\infty(\mathbb{R})$ for all $\lambda>0$, strictly monotonic, and bounded in $(0,1)$,
as illustrated in Figure~\ref{fig:logistic}.  Their derivatives satisfy the self-referential identities
\begin{equation}
  \frac{df^+}{dx} = \lambda\,f^+(1-f^+) > 0,
  \qquad
  \frac{df^-}{dx} = -\lambda\,f^-(1-f^-) < 0,
  \label{eq:selfref}
\end{equation}
with maximum absolute slope $\lambda/4$ attained at the threshold $x=\theta$:
\begin{equation}
  \left.\left|\frac{df^\pm}{dx}\right|\right|_{x=\theta}
  \;=\; \frac{\lambda}{4}.
  \label{eq:maxslope}
\end{equation}
Identity~\eqref{eq:selfref} makes Jacobian computation algebraic---no chain rule beyond the logistic values
themselves---and is the key property enabling the closed-form linearizations in
Section~\ref{sec:approx}.  The global $C^\infty$ smoothness stands in sharp contrast to the $C^k$
restriction of non-integer Hill functions discussed in Section~\ref{sec:hill_limitations}, and restores
access to center manifold reduction, normal form analysis, and Hessian-based parameter estimation
regardless of the cooperativity value.  Full proofs, antiderivative formulae, and Lipschitz constants are
derived in the subsections below.

\subsection{Closed-Form Derivatives with Self-Referential Structure}
\label{sec:logistic_derivatives}

The derivatives of the logistic functions admit remarkably elegant closed forms that follow from straightforward differentiation. We derive these explicitly to establish the self-referential structure that underpins all subsequent analysis.

\subsubsection{First derivative of the increasing logistic function}
Let $u = -\lambda(x - \theta)$, so that $f^+(x) = (1 + e^u)^{-1}$. By the chain rule:
\[
\frac{d f^+}{d x} = -(1 + e^u)^{-2} \cdot e^u \cdot \frac{du}{dx} = -(1 + e^u)^{-2} \cdot e^u \cdot (-\lambda) = \frac{\lambda e^{-\lambda(x-\theta)}}{\bigl(1 + e^{-\lambda(x-\theta)}\bigr)^2}.
\]
To express this in self-referential form, we factor the expression as a product and note that $f^+(x) = (1 + e^{-\lambda(x-\theta)})^{-1}$:
\[
\frac{d f^+}{d x} = \lambda \cdot \frac{1}{1 + e^{-\lambda(x-\theta)}} \cdot \frac{e^{-\lambda(x-\theta)}}{1 + e^{-\lambda(x-\theta)}} = \lambda f^+(x) \cdot \bigl(1 - f^+(x)\bigr),
\]
where the last step uses $e^{-\lambda(x-\theta)}/(1 + e^{-\lambda(x-\theta)}) = 1 - 1/(1 + e^{-\lambda(x-\theta)}) = 1 - f^+(x)$. Therefore:
\begin{equation}\label{eq:deriv_fplus}
\frac{\partial f^+}{\partial x} = \lambda f^+(x, \theta, \lambda)\bigl(1 - f^+(x, \theta, \lambda)\bigr).
\end{equation}

\subsubsection{First derivative of the decreasing logistic function}
By identical reasoning with $u = \lambda(x-\theta)$:
\[
\frac{d f^-}{d x} = -(1 + e^u)^{-2} \cdot e^u \cdot \lambda = -\frac{\lambda e^{\lambda(x-\theta)}}{\bigl(1 + e^{\lambda(x-\theta)}\bigr)^2}.
\]
Since $1 - f^-(x) = e^{\lambda(x-\theta)}/(1 + e^{\lambda(x-\theta)})$, this simplifies to:
\begin{equation}\label{eq:deriv_fminus}
\frac{\partial f^-}{\partial x} = -\lambda f^-(x, \theta, \lambda)\bigl(1 - f^-(x, \theta, \lambda)\bigr).
\end{equation}

\subsubsection{Second derivative and higher-order structure}
The self-referential form in \eqref{eq:deriv_fplus} enables systematic computation of all higher derivatives. Differentiating once more with respect to $x$, using the product rule and \eqref{eq:deriv_fplus}:
\begin{align}
\frac{d^2 f^+}{d x^2} &= \lambda \frac{d}{dx}\bigl[f^+(1 - f^+)\bigr]
= \lambda \left[ \frac{df^+}{dx}(1 - f^+) + f^+\left(-\frac{df^+}{dx}\right) \right] \notag \\
&= \lambda \frac{df^+}{dx}\bigl[(1 - f^+) - f^+\bigr]
= \lambda \cdot \lambda f^+(1 - f^+) \cdot (1 - 2f^+),
\end{align}
giving:
\begin{equation}\label{eq:second_deriv}
\frac{d^2 f^+}{d x^2} = \lambda^2 f^+(1 - f^+)(1 - 2f^+).
\end{equation}
At the inflection point $x = \theta$ where $f^+(\theta) = 1/2$, the factor $(1 - 2f^+(\theta)) = 0$, confirming that the second derivative vanishes at $x = \theta$. For $x < \theta$, $f^+ < 1/2$ so $(1-2f^+) > 0$, giving $f^{+\prime\prime} > 0$ (convexity). For $x > \theta$, $f^+ > 1/2$ so $(1-2f^+) < 0$, giving $f^{+\prime\prime} < 0$ (concavity).

\subsubsection{Third derivative}
Differentiating \eqref{eq:second_deriv} and applying the product rule, with $p = f^+$, $p' = \lambda p(1-p)$:
\begin{align}
\frac{d^3 f^+}{d x^3} &= \lambda^2 \frac{d}{dx}\bigl[p(1-p)(1-2p)\bigr] \notag \\
&= \lambda^2 \bigl[p'(1-p)(1-2p) + p(-p')(1-2p) + p(1-p)(-2p')\bigr] \notag \\
&= \lambda^2 p'\bigl[(1-p)(1-2p) - p(1-2p) - 2p(1-p)\bigr] \notag \\
&= \lambda^2 \cdot \lambda p(1-p) \cdot \bigl[(1-2p)(1-2p) - 2p(1-p)\bigr] \notag \\
&= \lambda^3 p(1-p)\bigl[(1-2p)^2 - 2p(1-p)\bigr].
\end{align}
Since $p(1-p) = (1 - (1-2p)^2)/4$, letting $q = 1 - 2f^+$ gives $p(1-p) = (1-q^2)/4$, and the third derivative takes the compact form:
\begin{equation}
\frac{d^3 f^+}{d x^3} = \lambda^3 f^+(1-f^+)\bigl(1 - 6f^+(1-f^+)\bigr).
\end{equation}
More generally, all derivatives of $f^+$ can be expressed as polynomials in $f^+$ (equivalently, in $\tanh(\lambda(x-\theta)/2)$), confirming global $C^\infty$ regularity.

\subsubsection{Uniform derivative bound}
Since $f^+(x) \in (0,1)$ for all $x \in \mathbb{R}$, the product $f^+(1-f^+)$ is maximized when $f^+ = 1/2$, giving $f^+(1-f^+) \leq 1/4$. Therefore, from \eqref{eq:deriv_fplus}:
\begin{equation}\label{eq:uniform_bound}
\left|\frac{\partial f^+}{\partial x}\right| = \lambda f^+(1-f^+) \leq \frac{\lambda}{4},
\end{equation}
with equality at $x = \theta$. The identical bound holds for $f^-$:
\[
\left|\frac{\partial f^-}{\partial x}\right| \leq \frac{\lambda}{4}.
\]
This uniform bound holds across the entire domain $\mathbb{R}$, guaranteeing that the system's Jacobian remains well-conditioned everywhere.

\subsection{Symmetry Property Around the Inflection Point}
\label{sec:logistic_symmetry}

The parametric logistic functions exhibit a fundamental symmetry. We have:
\[
f^+(x, \theta, \lambda) + f^-(x, \theta, \lambda) = \frac{1}{1 + e^{-\lambda(x-\theta)}} + \frac{1}{1 + e^{\lambda(x-\theta)}}.
\]
Multiplying the second term numerator and denominator by $e^{-\lambda(x-\theta)}$:
\[
f^+(x, \theta, \lambda) + f^-(x, \theta, \lambda) = \frac{1}{1 + e^{-\lambda(x-\theta)}} + \frac{e^{-\lambda(x-\theta)}}{e^{-\lambda(x-\theta)} + 1} = \frac{1 + e^{-\lambda(x-\theta)}}{1 + e^{-\lambda(x-\theta)}} = 1.
\]
Therefore:
\begin{equation}
f^+(x, \theta, \lambda) = 1 - f^-(x, \theta, \lambda), \qquad f^-(x, \theta, \lambda) = 1 - f^+(x, \theta, \lambda),
\end{equation}
centered at the inflection point $x = \theta$ where both functions equal $1/2$. This symmetry provides a natural probabilistic interpretation, as the functions balance positive and negative deviations from the threshold symmetrically, facilitating the transformation of mixed activation-repression networks into unified mathematical formulations.

\subsection{Analytically Invertible Closed Form}

The inverse functions possess elementary closed forms, derivable by direct algebraic manipulation. For the increasing function, we solve $y = f^+(x, \theta, \lambda)$:
\begin{align}
y &= \frac{1}{1 + e^{-\lambda(x-\theta)}} \notag\\
e^{-\lambda(x-\theta)} &= \frac{1-y}{y} \notag\\
-\lambda(x - \theta) &= \ln\!\left(\frac{1-y}{y}\right) \notag\\
x &= \theta + \frac{1}{\lambda}\ln\!\left(\frac{y}{1-y}\right) = \theta - \frac{1}{\lambda}\ln\!\left(\frac{1}{y} - 1\right), \quad y \in (0,1).
\end{align}
For the decreasing function, solving $y = f^-(x, \theta, \lambda)$ proceeds similarly:
\begin{align}
y &= \frac{1}{1 + e^{\lambda(x-\theta)}} \notag\\
e^{\lambda(x-\theta)} &= \frac{1-y}{y} \notag\\
\lambda(x - \theta) &= \ln\!\left(\frac{1-y}{y}\right) \notag\\
x &= \theta + \frac{1}{\lambda}\ln\!\left(\frac{1-y}{y}\right) = \theta + \frac{1}{\lambda}\ln\!\left(\frac{1}{y} - 1\right).
\end{align}
The logit transformation $\operatorname{logit}(y) = \ln\!\left(\frac{y}{1-y}\right)$ forms the foundation of logistic regression, a powerful tool for inferring gene regulatory networks from binary expression data.

\subsection{Global Infinite Differentiability ($C^\infty$)}

The parametric logistic functions are infinitely differentiable ($C^\infty$) across their entire domain, with all derivatives continuous everywhere, including boundary points. This follows from the fact that $f^+(x)= (1 + e^{-\lambda(x-\theta)})^{-1}$ is a composition of the globally analytic functions $x \mapsto -\lambda(x-\theta)$, $u \mapsto e^u$, and $v \mapsto (1+v)^{-1}$ (which is analytic since $1 + e^u > 0$ for all $u \in \mathbb{R}$). Compositions of analytic functions are analytic, hence $C^\infty$, and the same argument applies to $f^-$. The monotonicity, convexity, and inflection-point properties of $f^+$ follow directly from the first and second derivative expressions derived in Section~\ref{sec:logistic_derivatives}.

\subsection{Approximation and Expansion Theory}
\label{sec:approx}

\subsubsection{Hyperbolic identities and well-behaved linear approximations}

We first establish the exact identity relating logistic functions to the hyperbolic tangent, then derive the linear approximations. Starting from the definition:
\[
\tanh(z) = \frac{e^z - e^{-z}}{e^z + e^{-z}} = \frac{e^{2z} - 1}{e^{2z} + 1} = 1 - \frac{2}{e^{2z} + 1}.
\]
Setting $z = \lambda(x-\theta)/2$:
\[
\tanh\!\left(\frac{\lambda(x-\theta)}{2}\right) = 1 - \frac{2}{e^{\lambda(x-\theta)} + 1} = 1 - 2f^-(x, \theta, \lambda),
\]
and rearranging:
\begin{equation}\label{eq:tanh_fminus}
f^-(x,\theta,\lambda) = \frac{1}{2} - \frac{1}{2}\tanh\!\left(\frac{\lambda(x-\theta)}{2}\right).
\end{equation}
Using the symmetry $f^+ = 1 - f^-$:
\begin{equation}\label{eq:tanh_fplus}
f^+(x,\theta,\lambda) = \frac{1}{2} + \frac{1}{2}\tanh\!\left(\frac{\lambda(x-\theta)}{2}\right).
\end{equation}

Applying the Taylor expansion of $\tanh$, which follows from the power series $\tanh(z) = z - z^3/3 + 2z^5/15 - \cdots$ (with radius of convergence $\pi/2$, set by the nearest singularities of $\tanh$ at $z = \pm i\pi/2$), to the argument $z = \lambda(x-\theta)/2$ gives the linear approximations:
\[
f^+(x, \theta, \lambda) \approx \frac{1}{2} + \frac{\lambda(x-\theta)}{4}, \qquad
f^-(x, \theta, \lambda) \approx \frac{1}{2} - \frac{\lambda(x-\theta)}{4}.
\]
The full Taylor series converges for $|x-\theta| < \pi/\lambda$, but the linear approximation itself is accurate only when the cubic term is negligible relative to the linear term, i.e.\ when $\lambda|x-\theta| \ll 1$. The convergence radius $\pi/\lambda$ should therefore be understood as an upper bound on the validity of the truncated series, not a guarantee that the linear truncation is accurate throughout that range.

Multiple useful linear approximations facilitate different aspects of analysis:

\begin{itemize}
    \item \textbf{Near the origin ($x \approx 0$):} A first-order Taylor expansion of $f^+(x) = (1 + e^{-\lambda(x-\theta)})^{-1}$ around $x = 0$ gives, using $f^{+\prime}(x) = \lambda f^+(x)(1-f^+(x))$:
    \[
    f^+(x) \approx f^+(0) + f^{+\prime}(0)\cdot x = \frac{1}{1 + e^{\lambda\theta}} + \frac{\lambda e^{\lambda\theta}}{(1 + e^{\lambda\theta})^2}\, x.
    \]
    This captures basal expression levels (the constant term) and sensitivity to small concentration increases (the linear coefficient), both essential for analyzing low-expression regimes.
    
    \item \textbf{Around the threshold ($x \approx \theta$):}
    \[
    f^+(x) \approx \frac{1}{2} + \frac{\lambda}{4}(x - \theta).
    \]
    This midpoint linearization characterizes switch-like behavior with slope $\lambda/4$, providing the foundation for local stability analysis.
    
    \item \textbf{Near the threshold for small perturbations ($\lambda|x-\theta| \ll 1$):}
    Expanding in powers of $\lambda(x-\theta)$, since $e^{-\lambda(x-\theta)} = 1 - \lambda(x-\theta) + O(\lambda^2(x-\theta)^2)$:
    \[
    f^+(x) = \frac{1}{1 + 1 - \lambda(x-\theta) + O(\lambda^2(x-\theta)^2)}
    = \frac{1}{2 - \lambda(x-\theta) + O(\lambda^2(x-\theta)^2)}
    = \frac{1}{2}\cdot\frac{1}{1 - \tfrac{\lambda(x-\theta)}{2} + O(\lambda^2(x-\theta)^2)},
    \]
    and applying $(1-u)^{-1} \approx 1+u$ for small $u = \lambda(x-\theta)/2$:
    \[
    f^+(x) \approx \frac{1}{2}\!\left(1 + \frac{\lambda(x-\theta)}{2}\right) = \frac{1}{2} + \frac{\lambda(x-\theta)}{4} + O(\lambda^2(x-\theta)^2),
    \]
    valid whenever $\lambda|x-\theta| \ll 1$ (i.e., for deviations small relative to $1/\lambda$, regardless of whether $\lambda$ itself is small or large). In the special case $\lambda \ll 1$ this holds for $|x-\theta|$ of order unity.
    
    \item \textbf{Exponential approximation:} For $x \ll \theta$ with $\lambda\theta$ large, $e^{-\lambda(x-\theta)} = e^{\lambda(\theta - x)} \gg 1$, so:
    \[
    f^+(x) = \frac{1}{1 + e^{-\lambda(x-\theta)}} \approx e^{\lambda(x-\theta)},
    \]
    capturing the exponential decay of gene expression far below the threshold.
\end{itemize}

\subsubsection{Taylor series expansion}

Around the inflection point $x = \theta$, we substitute $z = \lambda(x-\theta)/2$ into the known Taylor series for $\tanh$:
\[
\tanh(z) = z - \frac{z^3}{3} + \frac{2z^5}{15} - \frac{17z^7}{315} + \cdots
\]
and apply identity \eqref{eq:tanh_fplus}. With $z = \lambda(x-\theta)/2$:
\begin{align}
f^+(x) &= \frac{1}{2} + \frac{1}{2}\left[\frac{\lambda(x-\theta)}{2} - \frac{1}{3}\!\left(\frac{\lambda(x-\theta)}{2}\right)^3 + \frac{2}{15}\!\left(\frac{\lambda(x-\theta)}{2}\right)^5 - \cdots \right] \notag\\
&= \frac{1}{2} + \frac{\lambda}{4}(x-\theta) - \frac{\lambda^3}{48}(x-\theta)^3 + \frac{\lambda^5}{480}(x-\theta)^5 + \mathcal{O}\bigl((x-\theta)^7\bigr),
\end{align}
where the coefficients are obtained as follows:
\begin{itemize}
    \item Coefficient of $(x-\theta)$: $\frac{1}{2} \cdot \frac{\lambda}{2} = \frac{\lambda}{4}$.
    \item Coefficient of $(x-\theta)^3$: $\frac{1}{2} \cdot \left(-\frac{1}{3}\right) \cdot \frac{\lambda^3}{8} = -\frac{\lambda^3}{48}$.
    \item Coefficient of $(x-\theta)^5$: $\frac{1}{2} \cdot \frac{2}{15} \cdot \frac{\lambda^5}{32} = \frac{\lambda^5}{480}$.
\end{itemize}
Note that only odd powers of $(x-\theta)$ appear, consistent with the odd-function symmetry of $f^+(x) - 1/2$ about $x = \theta$. The even-order derivatives of $f^+$ vanish at $x = \theta$. This rapidly converging series (with coefficients that decay faster than exponentially) facilitates Jacobian-based stability analysis, higher-order perturbation expansions in bifurcation studies, and reduced-order modeling.

\subsection{Closed-Form Integral}

We explicitly derive the antiderivative of $f^+(x)$. Let $u = x - \theta$, so $du = dx$:
\begin{align}
\int f^+(x)\,dx &= \int \frac{1}{1 + e^{-\lambda u}}\,du.
\end{align}
Multiply numerator and denominator by $e^{\lambda u}$:
\begin{align}
\int \frac{1}{1 + e^{-\lambda u}}\,du &= \int \frac{e^{\lambda u}}{e^{\lambda u} + 1}\,du.
\end{align}
Recognizing that the numerator $e^{\lambda u}$ is $(1/\lambda)$ times the derivative of the denominator $e^{\lambda u} + 1$:
\begin{align}
\int \frac{e^{\lambda u}}{e^{\lambda u} + 1}\,du &= \frac{1}{\lambda}\int \frac{\lambda e^{\lambda u}}{e^{\lambda u} + 1}\,du = \frac{1}{\lambda}\ln\!\left(e^{\lambda u} + 1\right) + C.
\end{align}
Substituting back $u = x - \theta$:
\begin{equation}\label{eq:integral_form1}
\int f^+(x)\,dx = \frac{1}{\lambda}\ln\!\bigl(1 + e^{\lambda(x-\theta)}\bigr) + C.
\end{equation}

An equivalent form is obtained by writing:
\begin{align}
\frac{1}{\lambda}\ln\!\bigl(1 + e^{\lambda(x-\theta)}\bigr) 
&= \frac{1}{\lambda}\ln\!\bigl(e^{\lambda(x-\theta)}\bigl(e^{-\lambda(x-\theta)} + 1\bigr)\bigr) \notag\\
&= \frac{1}{\lambda}\Bigl[\lambda(x-\theta) + \ln\!\bigl(1 + e^{-\lambda(x-\theta)}\bigr)\Bigr] \notag\\
&= (x - \theta) + \frac{1}{\lambda}\ln\!\bigl(1 + e^{-\lambda(x-\theta)}\bigr).
\end{align}
Therefore:
\begin{equation}\label{eq:integral_both}
\int f^+(x)\,dx = (x-\theta) + \frac{1}{\lambda}\ln\!\bigl(1 + e^{-\lambda(x-\theta)}\bigr) + C = \frac{1}{\lambda}\ln\!\bigl(1 + e^{\lambda(x-\theta)}\bigr) + C,
\end{equation}
confirming the equivalence stated in the text. This closed-form antiderivative enables exact solutions for transient dynamics in simple GRN motifs and proves invaluable when computing integrals arising in moment calculations for stochastic models.

\subsection{Rigorous Bounds and Stability}

\subsubsection{Uniform derivative bounds}
As established in Section~\ref{sec:logistic_derivatives} (see \eqref{eq:uniform_bound}), the derivatives satisfy:
\[
\left|\frac{\partial f^+}{\partial x}\right| \leq \frac{\lambda}{4}, \qquad \left|\frac{\partial f^-}{\partial x}\right| \leq \frac{\lambda}{4}.
\]
These bounds ensure well-conditioned Jacobian matrices across the entire domain, thereby guaranteeing numerical stability in eigenvalue analysis, perturbation expansions, and bifurcation studies.

\subsubsection{Lipschitz continuity}
The bounded derivatives immediately imply Lipschitz continuity with explicit constants. By the mean value theorem, for any $x_1, x_2 \in \mathbb{R}$ there exists $\xi$ between $x_1$ and $x_2$ such that:
\[
|f^+(x_1) - f^+(x_2)| = |f^{+\prime}(\xi)| \cdot |x_1 - x_2| \leq \frac{\lambda}{4}\,|x_1 - x_2|.
\]
Hence $f^+$ is Lipschitz continuous with constant $L = \lambda/4$. This Lipschitz property guarantees, by the Picard--Lindel\"{o}f theorem, the existence and uniqueness of solutions to ordinary differential equations involving logistic regulatory functions, with explicit control over solution sensitivity to initial conditions and parameters.

\subsubsection{Boundedness guarantees}
Since $f^+(x) \in (0, 1)$ and $f^-(x) \in (0, 1)$ for all $x \in \mathbb{R}$, any dynamical system constructed from logistic regulatory functions inherits natural boundedness properties. For a gene regulatory network with production terms of the form $\kappa f(x)$ and linear degradation $\gamma x$, the state space is bounded: $0 \leq x \leq \kappa/\gamma$. Analogous invariance and reduction arguments for transcription--translation systems and concave gene-expression models, in which the regulatory function is replaced by a Hill or piecewise-linear surrogate, are developed in~\cite{belgacem2018reduction,belgacem2014mathematical}; the logistic formulation here recovers the same invariance with the additional benefit of global $C^\infty$ smoothness across the entire state space.

\section{Mathematical Relationship Between Hill and Logistic Functions}
\label{sec:equivalence_and_scaling}

\subsection{The Coordinate-Change Identity and Its Meaning}
\label{sec:equivalence}

The Hill activation function and the logistic function are members of the same abstract family, related by a nonlinear change of input variable. Starting from the Hill function $h^+(x,\theta,n) = x^n/(x^n+\theta^n)$, substitute $u = x/\theta$ to write
\[
h^+(x,\theta,n) = \frac{u^n}{u^n + 1} = \frac{1}{1 + u^{-n}} = \frac{1}{1 + e^{-n\ln u}} = \frac{1}{1 + e^{-n\ln(x/\theta)}}.
\]
This is precisely the standard logistic function with steepness $n$ evaluated at the log-ratio $s = \ln(x/\theta)$:
\begin{equation}
\boxed{h^+(x,\theta,n) = \sigma\!\bigl(n\ln(x/\theta)\bigr)}, \qquad \sigma(z) = \frac{1}{1+e^{-z}}.
\label{eq:hill_logistic_identity}
\end{equation}
Equivalently, the logistic GRN function $f^+(x,\theta,\lambda) = \sigma(\lambda(x-\theta))$ applies the same sigmoid $\sigma$ to the \emph{affine} argument $\lambda(x-\theta)$ rather than the \emph{logarithmic} argument $n\ln(x/\theta)$. Identity~\eqref{eq:hill_logistic_identity} is exact and holds for all $x > 0$, $\theta > 0$, $n > 0$. Its meaning is straightforward: Hill and logistic are the \emph{same shape of function} ($\sigma$), but evaluated on inputs that differ by the coordinate change
\begin{equation}
\varphi : (0,\infty) \to (-\infty,\infty), \qquad \varphi(x) = n\ln(x/\theta),
\label{eq:coordinate_change}
\end{equation}
which maps the positive half-line to the whole real line via the natural logarithm.

\subsection{Why Mathematical Relatedness Does Not Imply Model Equivalence}
\label{sec:why_not_equivalent}

Identity~\eqref{eq:hill_logistic_identity} might suggest that choosing between Hill and logistic functions is merely a matter of re-parameterization, with no substantive biological consequence. This conclusion is incorrect. The two formulations give rise to \emph{different dynamical systems} when inserted into a GRN model, for four decisive reasons.

\subsubsection{1. Different behavior at zero concentration}
\[
h^+(0,\theta,n) = 0 \quad\text{for all } n > 0, \qquad f^+(0,\theta,\lambda) = \frac{1}{1+e^{\lambda\theta}} > 0 \quad\text{for all } \lambda,\theta > 0.
\]
The Hill-based ODE $\dot{x}_i = \kappa_i h^+(x_j) - \gamma_i x_i$ has an absorbing fixed point at $x_i = 0$ when all activators are absent: $h^+(0)=0$ gives zero production, so the state cannot escape zero by intrinsic production alone. The logistic ODE always maintains $f^+(0) > 0$, so the system cannot be permanently silenced. This difference drives qualitatively different long-run behavior (Section~\ref{sec:hill_pathologies}).

\subsubsection{2. Different smoothness and Lipschitz regularity}
For non-integer $n$, $h^+$ is only $C^{\lfloor n \rfloor}$ at $x=0$ (derivatives of order $> \lfloor n \rfloor$ diverge), whereas $f^+$ is $C^\infty$ on all of $\mathbb{R}$. The logistic ODE vector field is globally Lipschitz with finite constant $L_F \leq M$ (Theorem~\ref{thm:smoothness}), guaranteeing unique solutions and stable numerical integration everywhere. The Hill vector field, by contrast, fails to be Lipschitz with uniformly bounded derivatives of all orders on any neighborhood of the positive-orthant boundary: for $0 < n < 1$ the first derivative $\partial h^+/\partial x$ already diverges as $x \to 0^+$, so $h^+$ is not even locally Lipschitz; for $n > 1$, $h^+$ is locally Lipschitz but only $C^{\lfloor n \rfloor}$, since the $(\lfloor n \rfloor + 1)$-th derivative diverges at the origin. In both regimes the standard convergence theorems for Runge--Kutta methods of order $p > \lfloor n \rfloor$ cease to apply near $x=0$ (Section~\ref{sec:hill_pathologies}).

\subsubsection{3. Different sensitivity structures}
The logistic derivative $\lambda f^+(1-f^+) \leq \lambda/4$ is uniformly bounded and achieves its maximum at $x = \theta$, declining symmetrically on both sides. The Hill derivative $n\theta^n x^{n-1}/(x^n+\theta^n)^2$ scales as $x^{n-1}/\theta^n$ for $x \ll \theta$. The Hill function is sensitive to \emph{multiplicative} changes in $x/\theta$---a fixed fold-ratio in $x/\theta$ produces the same response increment regardless of the starting level---whereas the logistic function is sensitive to \emph{additive} changes $x - \theta$. (Here ``fold-change'' is in the Weber-Fechner / signal-processing sense; see Section~\ref{sec:biological_realism} for the distinction from ``fold-change in gene expression'' in thermodynamic models~\cite{bintu2005transcriptional,bintu2005transcriptional_app}.) These are distinct biological hypotheses about how a cell processes a regulatory signal, with direct consequences for parameter identifiability (see Section~\ref{sec:biological_implications}).

\subsubsection{4. The coordinate change does not commute with the ODE dynamics}
The nonlinear change $\varphi(x) = n\ln(x/\theta)$ does not transform the linear ODE structure into another linear ODE structure of the same form. Setting $s_i = \ln(x_i/\theta_i)$ (so $x_i = \theta_i e^{s_i}$ and $\dot{x}_i = \theta_i e^{s_i}\dot{s}_i$), the Hill ODE $\dot{x}_i = \kappa_i \sigma(n s_j) - \gamma_i x_i$ becomes, after dividing by $\theta_i e^{s_i}$,
\[
\dot{s}_i = \frac{\kappa_i}{\theta_i}\,e^{-s_i}\,\sigma(n s_j) - \gamma_i.
\]
This is \emph{not} of logistic ODE form: the production term acquires the multiplicative factor $e^{-s_i}$, which is a nonlinear function of $s_i$ itself and cannot be written as a logistic function of $s_j$ alone. The linear degradation $-\gamma_i x_i = -\gamma_i\theta_i e^{s_i}$ in the original $x$-equation reduces to the \emph{constant} $-\gamma_i$ in the transformed $\dot{s}_i$ equation, so it is the production side---not the degradation side---that is responsible for the non-equivalence. The two systems have different equilibria, different stability properties, and different trajectories.

The structural mathematical advantages of the logistic formulation---$C^\infty$ regularity, globally finite Lipschitz constant, nonzero basal output, and decoupled parameters $\lambda$ and $\theta$---hold independently of any biological interpretation and constitute the primary grounds for preferring the logistic framework in applications where numerical stability, bifurcation analysis, and observer design are required.

\subsection{Low-Steepness Approximation and Parameter Matching}

Expanding both functions around their shared inflection point $x = \theta$ (where both equal $1/2$) provides a useful local comparison. Writing $x = \theta + \epsilon$:
\[
h^+(x,\theta,n) \approx \tfrac{1}{2} + \tfrac{n}{4\theta}\epsilon + O(\epsilon^2/\theta^2), \qquad
f^+(x,\theta,\lambda) \approx \tfrac{1}{2} + \tfrac{\lambda}{4}\epsilon + O(\epsilon^2).
\]
The two are locally identical at $x = \theta$ when $\lambda = n/\theta$, the parameter matching condition used throughout this paper. However, the local equivalence extends only within a small neighborhood of $\theta$; away from the threshold, the responses diverge, with the Hill function extending along the logarithmic scale and the logistic along the linear scale.

The low-steepness limit also illuminates a structural difference. As $n \to 0$ or $\lambda \to 0$, both functions reduce to $\frac{1}{2}$ plus a small perturbation: $\frac{n}{4}\ln(x/\theta)$ for the Hill function, which diverges logarithmically as $x\to 0^+$, and $\frac{\lambda}{4}(x-\theta)$ for the logistic, which is linear and bounded below. The Hill function cannot be regularized at the origin by taking a small-steepness limit; the logistic can.

\subsection{Biological Implications of the Hill--Logistic Relationship}
\label{sec:biological_implications}

The coordinate-change identity reveals that choosing between the two models is equivalent to choosing whether gene regulatory responses are governed by \emph{multiplicative-increment} (fold-change) sensitivity in concentration or by \emph{additive} deviations from an absolute threshold. We use ``fold-change'' here in the Weber-Fechner / signal-processing sense (equal multiplicative increments produce equal responses); see Section~\ref{sec:biological_realism} for a precise distinction from ``fold-change in gene expression'' as used in the thermodynamic-model literature~\cite{bintu2005transcriptional,bintu2005transcriptional_app}.

\subsubsection{Evidence for additive threshold detection}
Several lines of experimental evidence favor additive threshold models in transcriptional regulation. First, the half-occupancy point of a transcription factor binding site---directly measured as a dissociation constant $K_d$ or $\text{EC}_{50}$---is an absolute concentration; the regulatory switch responds when $x$ crosses $\theta$, not when $x/\theta$ crosses some ratio~\cite{naqvi2023precise}. Second, studies on transcription factor dosage sensitivity, including haploinsufficiency and bistability in the \textit{E.~coli} galactose system, are more consistent with additive threshold crossing than with multiplicative-increment sensitivity~\cite{son2021nf}. Third, in risk assessment and pharmacology, the logistic function's linear concentration axis yields more interpretable points of departure and avoids logarithmic distortions that inflate variance at low doses and create singularities near zero expression~\cite{kappenberg2023guidance,waddell2004dose}. Fourth, single-cell genomics measurements (RNA-seq TPM, absolute transcript counts by smFISH) are naturally expressed in absolute units; logistic models integrate with these data without log-transformation artifacts that can distort inference at low-expression levels~\cite{escher2018advantages}. Multiplicative-increment sensitivity remains well established in bacterial chemotaxis, photoreceptor adaptation, and certain signaling cascades~\cite{frank2013input}; neither interpretation is universally correct, and both frameworks can be fitted to the same dose-response data with appropriate parameterization.

\subsubsection{Practical consequence of the bounded, decoupled derivative}
The maximum slope of the Hill function at threshold is $n/(4\theta)$, which couples steepness and threshold: a high-threshold gene ($\theta$ large) has a shallow response unless $n$ is inflated to compensate. This compensatory coupling is biologically undesirable and statistically harmful ---it creates ridge-like likelihood surfaces and high parameter correlations, as documented in Section~\ref{sec:practical_consequences}. The logistic formulation decouples these entirely: $\lambda$ controls steepness independently of $\theta$, so a steep response at a high threshold is prescribed by setting $\lambda$ large without any constraint from $\theta$. The uniform bound $|\partial f^\pm/\partial x| \leq \lambda/4$ established in Section~\ref{sec:logistic_foundations} then guarantees a well-conditioned Fisher information matrix across the full parameter range.

\subsubsection{Basal expression and low-expression dynamics}
The non-zero basal output $f^+(0,\theta,\lambda) = 1/(1+e^{\lambda\theta})$ reflects the persistent leaky transcription ($0.1$--$1$ mRNAs per cell in bacteria; $1$--$10$ TPM in mammalian cells) documented across gene regulatory systems~\cite{ozbudak2004multistability,fantom2014promoter}. This basal activity prevents irreversible trapping in off-states, primes the transcriptional machinery for rapid induction, and maintains responsiveness in bistable circuits under stochastic fluctuations~\cite{joanito2020basal,becskei2000engineering}. The logistic formulation captures this biology with a single interpretable parameter, product $\lambda\theta$, rather than an ad hoc additive offset. The Hill function's structural commitment to $h^+(0)=0$ and $h^-(0)=1$ requires artificial modifications for both activation and repression, adding parameters without biophysical grounding and complicating Jacobian-based stability analysis as documented in Section~\ref{sec:hill_limitations}.

In summary, while neither the fold-change nor the additive threshold interpretation is universally correct, the logistic formulation is better aligned with experimental evidence from transcriptional switches operating through molecular titration, provides a self-consistent treatment of basal expression across activation and repression, decouples steepness from threshold for better parameter identifiability, and preserves the structural mathematical advantages---$C^\infty$ regularity, global Lipschitz bound, decoupled parameters---that enable rigorous analysis regardless of the biological interpretation adopted.

\subsection{Summary of Model Differences}

Table~\ref{tab:hill_logistic_comparison} collects the principal mathematical distinctions between the two frameworks.

\begin{table}[ht]
\centering
\caption{Principal mathematical distinctions between Hill and logistic GRN models. Both functions take the same absolute concentration $x$ as input; the differences are structural, not a matter of ``absolute'' versus ``relative'' scaling.}
\label{tab:hill_logistic_comparison}
\begin{tabular}{lll}
\toprule
Property & Hill $h^+(x,\theta,n)$ & Logistic $f^+(x,\theta,\lambda)$ \\
\midrule
Value at $x=0$ & $0$ (absorbing state) & $1/(1+e^{\lambda\theta}) > 0$ (basal expression) \\
Smoothness & $C^{\lfloor n \rfloor}$ at $x=0$ (non-integer $n$) & $C^\infty$ everywhere \\
Lipschitz constant & $+\infty$ near $x=0$ & $\lambda/4 < \infty$ globally \\
Sensitivity structure & Fold-change (multiplicative in $x/\theta$) & Additive deviation from $\theta$ \\
Max slope & $n/(4\theta)$ (couples $n$ and $\theta$) & $\lambda/4$ (independent of $\theta$) \\
Closed-form inverse & No (fractional power for non-integer $n$) & Yes: $\theta + \lambda^{-1}\ln(y/(1-y))$ \\
Closed-form integral & No (hypergeometric for non-integer $n$) & Yes: $\lambda^{-1}\ln(1+e^{\lambda(x-\theta)})$ \\
Basal expression & Requires ad hoc offset $\varepsilon$ & Built in; controlled by $\lambda\theta$ \\
Parameter coupling & Steepness $n/\theta$ tied to threshold & $\lambda$ and $\theta$ independent \\
\bottomrule
\end{tabular}
\end{table}

The coordinate change $x \mapsto \ln(x/\theta)$ converts Hill into logistic at the level of a single function value, but does not convert one ODE model into the other and does not remove the structural differences listed above. The two frameworks are \emph{related but not equivalent} as models of gene regulation; the logistic formulation is preferred on the grounds of the structural mathematical advantages in the table.


\section{Practical Consequences of Hill Function Limitations}
\label{sec:practical_consequences}

The mathematical pathologies catalogued in Section~\ref{sec:hill_limitations}
are not isolated analytical curiosities; each propagates into a specific failure
mode in simulation, stability analysis, parameter estimation, and control design.
This section traces these failures explicitly, identifying in each case the exact
property of the logistic function---derived in Section~\ref{sec:logistic_foundations}---that
eliminates the failure.  The arguments are entirely analytical; systematic
quantitative comparisons are deferred to future computational work.

\subsection{Numerical Integration}
\label{sec:numerical_integration}

The fundamental obstacle for ODE solvers is the absence of a global Lipschitz
constant for Hill-function vector fields with non-integer~$n$.  For the
scalar activation case, the Jacobian entry is
$\partial h^+/\partial x \sim (n/\theta^n)\,x^{n-1}$ as $x \to 0^+$
(Section~\ref{sec:hill_limitations}). The behavior splits into two regimes.
For $0 < n < 1$, $\partial h^+/\partial x$ diverges as $x \to 0^+$, so the
supremum of $|\partial h^+/\partial x|$ over any neighborhood of the
positive-orthant boundary is infinite, and the Hill vector field is
\emph{not} locally Lipschitz there: the classical existence--uniqueness
theorem does not apply directly, and the standard error bound
$\|\mathbf{e}(t)\| \leq \|\mathbf{e}(0)\|\,e^{Lt}$ provides no useful
guarantee since $L = \infty$.
For $n > 1$, the first derivative remains bounded near the origin (it
vanishes as $x \to 0^+$ for $n > 1$), so $h^+$ is locally Lipschitz, but
the higher-order derivatives diverge: for $n \in (k,k+1)$ the
$(k{+}1)$-th derivative $\partial^{k+1} h^+/\partial x^{k+1}$ blows up at
$x=0$, leaving $h^+$ only $C^{\lfloor n \rfloor}$ at the boundary.
Standard convergence and stability theorems for Runge--Kutta and multistep
methods of order $p > \lfloor n \rfloor$ require derivatives of order
$p+1$ to be bounded along the trajectory; these fail near $x=0$, and
adaptive solvers that internally estimate higher-order derivatives produce
spuriously large local error estimates and are forced to reduce step size
dramatically. In both regimes, then, the usual quantitative convergence
guarantees break down on any neighborhood of the positive-orthant boundary.

Furthermore, any trajectory that passes through a neighborhood of zero---as
low-expression dynamics necessarily do---encounters a vector field for which
$x^n$ is real-valued only for $x \geq 0$ (a related class of singular/Zeno difficulties in piecewise-linear and hybrid GRN models is addressed by probabilistic-convolution regularization in~\cite{belgacem2019probabilistic}); rounding errors of order $10^{-15}$
that produce a momentarily negative component introduce complex arithmetic
into every subsequent integration step, corrupting the trajectory silently and
irrecoverably.

The logistic vector field eliminates both obstacles by construction.  The
uniform bound $|\partial f^+/\partial x| \leq \lambda/4$ (equation~\eqref{eq:uniform_bound})
implies that the Jacobian entries of the logistic GRN vector field satisfy
$|\partial F_i/\partial x_j| \leq \kappa_i L_i^j + \gamma_i\delta_{ij}$
uniformly over all $\mathbf{x} \in \mathbb{R}^n$, yielding the finite global
Lipschitz constant $L_F \leq M$ established in Theorem~\ref{thm:smoothness},
Part~(ii).  Standard convergence theorems therefore apply to all of
$\mathbb{R}^n$, including at and below the positive-orthant boundary.  The
stiffness ratio (ratio of the largest to the smallest absolute eigenvalue of the
Jacobian) is bounded by $M/\min_i\gamma_i$, a finite, parameter-determined
constant, in contrast to the Hill case, where no such bound exists near zero.

\subsection{Stability and Bifurcation Analysis}

Computing equilibrium stability by eigenvalue analysis of the Jacobian, and
performing bifurcation analysis by normal form theory both require
derivatives of the vector field to a prescribed order.  For the Hill function
with $n \in (k, k+1)$, derivatives of order greater than $k$ diverge at
$x = 0$ (Section~\ref{sec:hill_limitations}), so the vector field is
only $C^k$ at the boundary of the positive orthant.  This finite
differentiability class has three direct analytical consequences.

First, eigenvalue computations at equilibria with components near zero are
ill-conditioned: the Jacobian entries $\partial F_i/\partial x_j \propto x_j^{n-1}$
either diverge (for $n < 1$) or vanish (for $n > 1$) as $x_j \to 0^+$,
inflating numerical errors in eigenvalue decompositions and potentially
reversing the sign of eigenvalues, thereby misclassifying stability.

Second, normal form theory for local bifurcation analysis requires Taylor
expansions of the vector field to order equal to the normal form degree.
For a pitchfork or Hopf bifurcation, at least third-order terms are
essential; these terms involve derivatives that are undefined for
$n \in (1, 2)$, so the normal form computation cannot be completed near
the boundary.  This masks degenerate bifurcation structures---including
degenerate Hopf bifurcations and homoclinic orbits---that require terms of
order three or higher.

Third, parameter continuation of equilibria along bifurcation branches
requires the Jacobian to be a continuously differentiable function of
parameters.  For a Hill system, $D\mathbf{F}$ is discontinuous as a function
of parameters whenever an equilibrium component approaches zero for
non-integer $n$; continuation algorithms will generically encounter
apparent termination points or spurious bifurcations at such loci.

The logistic vector field is globally $C^\infty$ (Section~\ref{sec:logistic_foundations}
and Theorem~\ref{thm:smoothness}, Part~(i)), and all derivatives can be
expressed as polynomials in the function values $f^\pm \in (0,1)$ through
the self-referential structure \eqref{eq:deriv_fplus}--\eqref{eq:second_deriv}.
The Taylor expansion around any equilibrium converges (Section~\ref{sec:approx} gives
the explicit coefficients in terms of~$\lambda$), normal forms of all
orders are computable, and $D\mathbf{F}$ is analytic in all parameters, so
continuation is unobstructed throughout the biologically relevant domain. Stability and bifurcation results for related smooth and piecewise-linear gene-network models---including reduction-based stability analysis of transcription--translation systems~\cite{belgacem2013stability,belgacem2014stability,belgacem2018reduction} and bifurcation/chaos analysis in Glass-network ring circuits~\cite{belgacem2025chaos,farcot2019chaos}---transfer directly to the logistic framework once the Hill regulator is replaced, since each of these analyses relies on the smoothness and Lipschitz regularity that the logistic formulation guarantees globally.

\subsection{Parameter Estimation}

Fitting Hill models to experimental data requires the parameter gradients
$\partial h^+/\partial\theta$ and $\partial h^+/\partial n$.  From
Section~\ref{sec:hill_limitations}, as $x \to 0^+$:
\[
  \frac{\partial h^+}{\partial n}   \sim \left(\frac{x}{\theta}\right)^n \ln\!\left(\frac{x}{\theta}\right), \qquad
  \frac{\partial h^+}{\partial\theta} \sim -\frac{n}{\theta}\left(\frac{x}{\theta}\right)^n,
\]
so their ratio diverges:
\[
  \frac{|\partial h^+/\partial n|}{|\partial h^+/\partial\theta|} \approx \frac{\theta\,|\ln(x/\theta)|}{n} \xrightarrow{x\to 0^+} \infty.
\]
This logarithmically divergent ratio is the direct cause of Fisher information
ill-conditioning.  The Fisher information matrix
\[
  \mathcal{I}_{ij} \approx \frac{1}{\sigma^2}\sum_{k=1}^N
  \frac{\partial h^+(x_k,\mathbf{p})}{\partial p_i}\,
  \frac{\partial h^+(x_k,\mathbf{p})}{\partial p_j}
\]
accumulates this divergent ratio from any design points~$x_k$ near zero,
causing one eigenvalue of~$\mathcal{I}$ to dominate while the other
approaches zero: $\kappa(\mathcal{I}) = \lambda_{\max}/\lambda_{\min} \to \infty$.
By the Cram\'{e}r--Rao bound~\cite{rao1945information,cramer1999mathematical},
$\operatorname{Cov}(\hat{\mathbf{p}}) \geq \mathcal{I}^{-1}$, so the
near-singular~$\mathcal{I}$ yields unbounded parameter variances and strong
cross-correlation between $\hat{n}$ and $\hat{\theta}$, producing ridge-like
likelihood surfaces on which gradient-based optimizers stall.  For
$0 < n < 1$, the concentration sensitivity $\partial h^+/\partial x \propto x^{n-1} \to \infty$
simultaneously destabilizes gradient updates from near-zero observations,
compounding the identifiability failure.

The logistic formulation resolves this by decoupling steepness from threshold.
The parameter gradients of $f^+(x,\theta,\lambda)$ are
\[
  \frac{\partial f^+}{\partial \lambda} = (x-\theta)\,f^+(1-f^+), \qquad
  \frac{\partial f^+}{\partial \theta} = -\lambda f^+(1-f^+),
\]
which, using the uniform bound~\eqref{eq:uniform_bound}, satisfy
\[
  \left|\frac{\partial f^+}{\partial\lambda}\right| \leq \frac{|x-\theta|}{4},
  \qquad
  \left|\frac{\partial f^+}{\partial\theta}\right| \leq \frac{\lambda}{4},
\]
on the entire real line.  Crucially, \emph{neither gradient diverges} at any
concentration---including near zero, where Hill gradients become pathological.
Moreover, the two gradients vanish at \emph{different} points:
$\partial f^+/\partial\lambda$ vanishes only at $x = \theta$, while
$\partial f^+/\partial\theta = -\lambda/4$ attains its maximum magnitude there
and is bounded away from zero on every compact set.  As a consequence, the
gradient direction
$\bigl(\partial f^+/\partial\lambda,\ \partial f^+/\partial\theta\bigr) =
f^+(1-f^+)\,(x-\theta,\,-\lambda)$ varies across the design points $\{x_k\}$,
and the Fisher information matrix
$\mathcal{I} \approx \sigma^{-2}\sum_k \nabla f^+(x_k)\,\nabla f^+(x_k)^\top$
is full rank with bounded condition number whenever the design contains at
least two distinct concentrations.  This stands in sharp contrast to the Hill
case, where the gradient direction collapses to a single asymptotic ray as
$x \to 0^+$---driving one eigenvalue of $\mathcal{I}$ to zero and the other
to infinity.  Moreover, the Jacobian of~$\mathbf{F}$ with respect to parameters
is itself globally Lipschitz (Remark~\ref{rem:lipschitz_jacobian}), so
Gauss--Newton and Levenberg--Marquardt iterations converge reliably from any
initialization.

\subsection{Observer Design and Observability}

State estimation from partial measurements requires that the observability
matrix
\[
  \mathcal{O} = \bigl[\nabla h,\;\nabla(\mathcal{L}_F h),\;\nabla(\mathcal{L}_F^2 h),\;\ldots\bigr]^\top,
\]
formed from successive Lie derivatives of the output function~$h$ along
the vector field~$F$, be of full rank.  Each row inherits the derivative
structure of the regulatory functions.  For the repression function
$h^-(x_j,\theta,n)$ with cooperative exponent $n > 1$, the gradient that
enters the second row of~$\mathcal{O}$ satisfies
\[
  h^{-\prime}(x_j) = -\frac{n\theta^n x_j^{n-1}}{(\theta^n+x_j^n)^2}
  \xrightarrow{x_j \to 0^+} 0,
\]
so~$\mathcal{O}$ loses rank in the low-expression regime ($x_j \to 0^+$),
which is precisely where regulatory transitions are most active and where measurement noise is most severe.  For $0 < n < 1$, the same derivative diverges, maintaining formal rank but introducing cusp-like singularities
that create severe ill-conditioning in the gain matrices of any
observer that relies on inverting or weighting the observability structure.

The logistic repression derivative~\eqref{eq:deriv_fminus},
\[
  f^{-\prime}(x_j) = -\lambda f^-(x_j)\bigl(1 - f^-(x_j)\bigr),
\]
is strictly nonzero for all $x_j \in \mathbb{R}$, since $f^-(x_j) \in (0,1)$
implies $f^-(x_j)(1-f^-(x_j)) > 0$ everywhere.  The corresponding row of
$\mathcal{O}$ therefore never vanishes, and $\operatorname{rank}(\mathcal{O}) = n$
is maintained throughout the biologically relevant domain, including the
limit $x_j \to 0$.  Furthermore, the global Lipschitz property
$L_F \leq M$ of Theorem~\ref{thm:smoothness} provides the condition
required for high-gain observer synthesis: a gain $L > L_F$ guarantees
error convergence at rate $e^{-(L - L_F)t}$, with the explicit and
computable bound $L_F \leq M$ anchoring the gain selection to biological
parameters $(\kappa_i, \gamma_i, \lambda_{i,m}, \theta_{i,m})$.

\subsection{Control Design}

The differentiability and invertibility properties documented in
Sections~\ref{sec:logistic_foundations}--\ref{sec:hill_limitations} have
direct consequences for several control-theoretic paradigms used in
synthetic and systems biology: linear quadratic regulation (LQR), exact
feedback linearization, sliding mode control (SMC), model predictive
control (MPC), extended Kalman filtering (EKF), and linear-approximation
based gain design. We treat each in turn.

\subsubsection{LQR and Riccati-based synthesis}
LQR requires solving an algebraic Riccati equation whose coefficients
depend on the system Jacobian $J = D\mathbf{F}(\mathbf{x}^*)$ at the
operating point. For a logistic-based GRN, each off-diagonal Jacobian entry takes the form $\pm\kappa_i\lambda_{i,m}\,g_{i,m}(1-g_{i,m})\!\prod_{q\neq m} g_{i,q}$, expressed entirely in terms of the logistic factor values themselves and uniformly bounded by $\kappa_i\lambda_{i,m}/4$ via~\eqref{eq:bound_partial_fi} (Section~\ref{sec:logistic_foundations}). The Riccati equation is therefore posed with a well-conditioned coefficient matrix throughout the biologically relevant domain. For a Hill-based formulation with non-integer $n$, the Jacobian entry scales as $x^{n-1}$ near the origin, so the Riccati equation inherits an ill-conditioned structure precisely in the low-expression regime where regulatory transitions most often occur~\cite{santillan2008use}.

\subsubsection{Feedback linearization and the logit coordinate change}
Exact feedback linearization requires a smooth, invertible coordinate
change constructed from the inverse of the regulatory function. For the
Hill function, the formal inverse $(h^+)^{-1}(y) = \theta(y/(1-y))^{1/n}$
involves a fractional power that is multi-valued for non-positive
arguments. Direct differentiation gives
\[
\frac{d(h^+)^{-1}}{dy} = \frac{\theta}{n}\cdot\frac{1}{y(1-y)}\,\left(\frac{y}{1-y}\right)^{1/n},
\]
which by the inverse function theorem equals $1/h^{+\prime}(x)$ where $x = (h^+)^{-1}(y)$.
Since $h^{+\prime}(x) \to 0$ as $x \to \infty$ (i.e., as $y \to 1^-$) for every $n > 0$,
$d(h^+)^{-1}/dy$ diverges at $y \to 1^-$ in all cases. At $y \to 0^+$ (i.e., $x \to 0^+$),
the behavior is dichotomous: for $n > 1$, $h^{+\prime}(0) = 0$, so $d(h^+)^{-1}/dy \to \infty$;
for $0 < n < 1$, $h^{+\prime}(x) \to \infty$, so $d(h^+)^{-1}/dy \to 0$, an equally problematic
\emph{loss of sensitivity} that prevents the inverse map from resolving small concentration changes.
In every non-trivial case ($n \neq 1$), the inverse map is pathological at one or both endpoints
of $(0,1)$, making the linearizing control law singular precisely in the low- and high-expression
regimes relevant to biological switching.

For the logistic function, the inverse is the elementary logit
transformation (Section~\ref{sec:logistic_foundations}):
\begin{equation}
(f^+)^{-1}(y) = \theta + \frac{1}{\lambda}\ln\!\left(\frac{y}{1-y}\right),
\qquad
\frac{d(f^+)^{-1}}{dy} = \frac{1}{\lambda\,y(1-y)},
\qquad y\in(0,1),
\label{eq:logit_inverse_control}
\end{equation}
which is single-valued, smooth, and bounded on every compact subset of
$(0,1)$. The linearizing control law can therefore be written in closed
form without numerical iteration; the coordinate change exactly
linearizes the GRN dynamics, enabling pole placement, LQR optimal
control, and $\mathcal{H}_\infty$ synthesis through standard linear
systems theory~\cite{belgacem2020control}.

\subsubsection{Sliding mode control and basal activity}
Sliding mode control (SMC) stabilizes a system at a target state by
imposing a sliding surface $s_i = 0$ and choosing a discontinuous control
law satisfying the reaching condition $s_i\dot{s}_i < 0$. Computing the
equivalent control---the smooth component that keeps the system on the
surface---requires inverting the input-to-state map at the desired
setpoint. For logistic-based models, this inversion reduces to evaluating
the closed-form logit~\eqref{eq:logit_inverse_control}, yielding
equivalent-control expressions bounded and well-defined across the entire
operating range. For the corresponding Hill-based system, the equivalent
control involves fractional powers of the form
$\theta\bigl(\gamma x^*/(\kappa - \gamma x^*)\bigr)^{1/n}$, which is
undefined when $x^* = 0$ and ill-conditioned when $x^*$ is small, making
SMC design analytically fragile at low-expression operating
points~\cite{belgacem2020control,chambon2020qualitative}. The logistic
function's strictly positive basal output
$f^+(0,\theta,\lambda) = 1/(1+e^{\lambda\theta}) > 0$ ensures that the
system never enters an absorbing zero state, so equivalent-control
expressions remain well-defined at all biologically relevant setpoints.
A boundary-layer saturation function can replace the discontinuous sign
term to reduce chattering while preserving the reaching condition, and
the logistic function's smooth, bounded structure ensures that small
perturbations in $\kappa_i$, $\gamma_i$, or $\theta_{i,m}$ do not produce
unbounded or discontinuous responses in the control law.

\subsubsection{Model predictive control and gradient-based optimization}
Gradient-based MPC solvers require the Hessian of the objective with
respect to the control inputs, which through the chain rule involves
second derivatives of the production terms. For the Hill function near
the origin, $|h^{+\prime\prime}(x)| \propto x^{n-2}$ diverges for $n < 2$
(Section~\ref{sec:hill_limitations}), rendering the Hessian
ill-conditioned whenever any predicted state passes through a
low-expression phase during the optimization horizon. For the logistic
function, the second derivative
$|f^{+\prime\prime}(x)| = \lambda^2|f^+(1-f^+)(1-2f^+)|$ is uniformly
bounded by $\lambda^2\rho$ with $\rho = \sqrt{3}/18 \approx 0.096$
(Remark~\ref{rem:lipschitz_jacobian}), so the Hessian is
well-conditioned at every operating point. The closed-form
inverse~\eqref{eq:logit_inverse_control} further enables feedback
linearization within the MPC prediction model, reducing each online
optimization to a well-posed problem amenable to real-time
solution~\cite{lugagne2024deep}. By contrast, Hill-based MPC solvers
require extensive line searches, increased iterations, or fail to
converge within real-time constraints when any gene operates near basal
expression during the prediction horizon.

\subsubsection{Extended Kalman filtering}
For extended Kalman filters, the same Jacobian enters the covariance
propagation equation. The logistic bound
$|\partial f^\pm/\partial x| \leq \lambda/4$ guarantees that the predicted
covariance matrix $P_{k+1|k}$ remains positive definite throughout, in
contrast to Hill-based EKFs where Jacobian entries proportional to
$\hat{x}_j^{n-1}$ diverge as $\hat{x}_j \to 0^+$ for $0 < n < 1$, causing
$P_{k+1|k}$ to lose positive semi-definiteness and the filter to fail
irrecoverably~\cite{santillan2008use}.

\subsubsection{Linearization-based gain design}
The Taylor expansion of the logistic function around its inflection
point, $f^+(x,\theta,\lambda) \approx \tfrac{1}{2} + \tfrac{\lambda}{4}(x-\theta)$,
is accurate when $\lambda|x-\theta|\ll 1$ (Section~\ref{sec:approx}; the underlying series converges on $|x-\theta|<\pi/\lambda$). Linearizing system~\eqref{eq:logistic_system_thm} about an operating point $\mathbf{x}^*$ (an equilibrium, or any state of interest) yields an affine system
$\dot{\boldsymbol{\delta}} = A\boldsymbol{\delta}$ in deviation
coordinates $\boldsymbol{\delta} = \mathbf{x}-\mathbf{x}^*$, with matrix entries
\[
A_{ij} = \kappa_i\,\frac{\partial f_i}{\partial x_j}\bigg|_{\mathbf{x}^*} - \gamma_i\,\delta_{ij}
       = \sigma_{i,m}\,\kappa_i\lambda_{i,m}\,g_{i,m}^*(1-g_{i,m}^*)\!\!\prod_{q\neq m} g_{i,q}^* - \gamma_i\,\delta_{ij},
\]
where $g_{i,m}^* := g_{i,m}(x^*_{j(i,m)})$ and the sum collapses to a single $m$ when each regulator drives a unique target index $j$. The bound $|A_{ij}-(-\gamma_i\delta_{ij})|\le \kappa_i\lambda_{i,m}/4$ from~\eqref{eq:bound_partial_fi} holds at every operating point. In the canonical single-regulator case ($M_i = 1$, no self-loops) and at the threshold $\mathbf{x}^*$ where $g_{i,1}^* = 1/2$, the off-diagonal entry simplifies to $\sigma_{i,1}\kappa_i\lambda_{i,1}/4$ and the diagonal entry is exactly $-\gamma_i$; for genes with $M_i$ regulators evaluated simultaneously at threshold, the off-diagonal entry is $\sigma_{i,m}\kappa_i\lambda_{i,m}/2^{M_i+1}$, and self-regulation contributes an additional term $\kappa_i\partial f_i/\partial x_i|_{\mathbf{x}^*}$ to the diagonal.
This linear structure enables symbolic derivation of stability
conditions, pole placement, and LQR gains directly from biological
parameters, and supports controllability analysis via the rank of the
controllability matrix $\mathcal{C} = [B,\,AB,\,\ldots,\,A^{n-1}B]$. A
state-feedback law $\mathbf{u} = K(\mathbf{x}_d - \mathbf{x})$ yields a
closed-loop matrix $A - BK$ whose eigenvalues can be placed at any
desired locations through standard pole-placement algorithms, enabling
systematic tuning of convergence speed and robustness margins.

By contrast, the small-$n$ approximation
$h^+(x,\theta,n) \approx \tfrac{1}{2} + \tfrac{n}{4}\ln(x/\theta)$
(equation~\eqref{eq:hill_small_n}) diverges logarithmically as
$x \to 0^+$, preventing a globally valid affine reduction and confining
linear approximations to a neighborhood bounded strictly away from
zero---excluding precisely the basal-expression regimes relevant to
bistability analysis and induction dynamics. Even the local
linearization around $\theta$,
$h^+(x) \approx \tfrac{1}{2} + \tfrac{n}{4\theta}(x-\theta)$, is valid
only in a neighborhood of $\theta$ bounded strictly away from zero.

\medskip
Each of the control- and estimation-theoretic failure modes identified above is therefore resolved by a specific structural property of the logistic function established in Section~\ref{sec:logistic_foundations}. Theorem~\ref{thm:smoothness} (Section~\ref{sec:advanced_analysis}) lifts these pointwise properties---closed-form derivatives, integrals, inverses, and Lipschitz bounds---to the level of the full GRN vector field, securing global existence, smoothness, and uniform boundedness of solutions.

\section{Biological Realism of Logistic Functions}
\label{sec:biological_realism}

The preceding sections established the mathematical foundations of the logistic framework (Section~\ref{sec:logistic_foundations}), the precise mathematical relationship to the Hill function (Section~\ref{sec:equivalence_and_scaling}), and the practical computational and control-theoretic advantages it confers over Hill-based models (Section~\ref{sec:practical_consequences}). This section addresses the complementary question of biological fidelity: do logistic functions capture the regulatory phenomena that matter biologically? We show that they do, and in some respects do so more naturally than Hill functions.

\subsection{Sensitivity Structure: Fold-Change Detection Versus Additive Threshold Crossing}

It is important to clarify a potential confusion about what distinguishes Hill and logistic functions as \emph{models}. Both functions take the absolute concentration $x$ (in nM or molecules per cell) as their argument, and both produce a dimensionless output in $(0,1)$. The distinction between them is therefore \emph{not} a matter of ``absolute'' versus ``relative'' concentrations---both use absolute concentrations. The distinction lies instead in the \emph{structure of their sensitivity}: how the regulatory response changes as $x$ varies. (Throughout this section, ``fold-change'' in the phrase ``multiplicative-increment sensitivity'' or ``fold-change sensitivity'' refers to equal ratios of $x/\theta$ producing equal response changes---the Weber-Fechner notion. This is entirely distinct from ``fold-change in gene expression'' as used in thermodynamic models of transcriptional regulation~\cite{bintu2005transcriptional,bintu2005transcriptional_app}, where fold-change denotes the ratio of promoter activity in the presence versus absence of a transcription factor.)

To see this precisely, note the exact algebraic identity (derived in Section~\ref{sec:equivalence}):
\begin{equation}
h^+(x,\theta,n) = \frac{x^n}{x^n+\theta^n} = \frac{1}{1 + e^{-n\ln(x/\theta)}} = f^+\!\bigl(\ln(x/\theta),\, 0,\, n\bigr).
\label{eq:hill_as_logistic_of_log}
\end{equation}
The Hill function is therefore the standard logistic function evaluated not at the absolute concentration $x$ but at the log-ratio $\ln(x/\theta)$. This single difference in the input variable is the source of the distinct sensitivity structures of each model:

\begin{itemize}
\item \textbf{Hill / multiplicative-increment sensitivity.} The slope of $h^+$ with respect to $x$ is
\[
{h^+}'(x) = \frac{n}{\theta}\cdot\frac{(x/\theta)^{n-1}}{(1+(x/\theta)^n)^2},
\]
which scales as $n\, x^{n-1}/\theta^n$ for $x \ll \theta$ and as $n\,\theta^n/x^{n+1}$ for $x \gg \theta$. At the threshold $x=\theta$, the slope attains its maximum value $n/(4\theta)$. In the log-input variable $s = \ln(x/\theta)$, the chain rule gives $\partial h^+/\partial s = x\,h^{+\prime}(x) = n\sigma'(ns)$, which is bounded by $n/4$ uniformly in $x$. The Hill function therefore has constant slope per equal multiplicative increment of $x$: equal log-ratio increments produce equal response changes, regardless of the starting level. This is its \emph{multiplicative-increment sensitivity}: a doubling from $1$\,nM to $2$\,nM produces the same response change as a doubling from $10$\,nM to $20$\,nM. (As noted in the abstract, this is distinct from ``fold-change in gene expression'' as used in thermodynamic models~\cite{bintu2005transcriptional,bintu2005transcriptional_app}.)

\item \textbf{Logistic / additive threshold sensitivity.} The slope of $f^+$ with respect to $x$ is $\lambda f^+(1-f^+)$, which is maximized at $x = \theta$ with value $\lambda/4$ and decays symmetrically on either side. Equal \emph{additive} increments of $x$ produce the same shift in $f^+$ when they occur at the same distance from $\theta$. A gain of $1$\,nM starting at $x = \theta - 0.5$\,nM produces the same response change as a gain of $1$\,nM starting at $x = \theta + 0.5$\,nM, but a markedly different response from a gain of $1$\,nM starting at $x = \theta + 10$\,nM.
\end{itemize}

\subsubsection{Why the mathematical equivalence \eqref{eq:hill_as_logistic_of_log} does \emph{not} make the models equivalent}
Equation~\eqref{eq:hill_as_logistic_of_log} identifies the Hill function as a special case of the logistic family, but this does \emph{not} mean the two ODE models are equivalent. The gene regulatory system
\[
\dot{x}_i = \kappa_i\, h^+(x_j, \theta, n) - \gamma_i x_i
\]
is a \emph{different dynamical system} from
\[
\dot{x}_i = \kappa_i\, f^+(x_j, \theta, \lambda) - \gamma_i x_i,
\]
even when $\lambda = n/\theta$. In particular:
\begin{enumerate}
\item $h^+(0) = 0$ while $f^+(0) = 1/(1+e^{\lambda\theta}) > 0$: the models make different predictions about basal expression, which drives qualitatively different long-term dynamics (Section~\ref{sec:hill_pathologies}).
\item $h^+$ is only $C^{\lfloor n \rfloor}$ at $x = 0$ for non-integer $n$, while $f^+$ is $C^\infty$ everywhere (Theorem~\ref{thm:smoothness}): the models have different regularity, with direct consequences for numerical integration and bifurcation analysis.
\item The Lipschitz constant of the Hill vector field is infinite near $x = 0$, while the logistic vector field is globally Lipschitz: the models have fundamentally different stability guarantees (Section~\ref{sec:hill_pathologies}).
\end{enumerate}
The coordinate change $s = \ln(x/\theta)$ converts one function value into the other at a single point, but it does not convert one ODE into the other: the production term acquires a state-dependent factor that breaks the logistic ODE structure (see Section~\ref{sec:why_not_equivalent}, item~4).

\subsubsection{Which sensitivity structure is biologically appropriate?}
Neither is universally correct; the choice reflects an assumption about how the cell ``computes'' its regulatory response. Multiplicative-increment sensitivity (logarithmic input scale) is well established in adaptation systems---bacterial chemotaxis, photoreceptor adaptation, and some signaling pathways---where the relevant quantity is the ratio of current to previous ligand levels~\cite{frank2013input}. Additive threshold detection is more natural for gene regulatory switches in which a fixed number of DNA-bound molecules is titrated as activator concentration rises. The half-occupancy point, directly measured as a dissociation constant $K_d$ or $\text{EC}_{50}$, is an absolute concentration; the biological switch fires when $x$ crosses $\theta$, not when $x/\theta$ crosses some ratio. Experimental evidence on transcription factor dosage effects---haploinsufficiency, gene copy number sensitivity, and threshold-dependent bistability in the \textit{E.~coli} gal system---is more consistent with additive threshold detection~\cite{naqvi2023precise,son2021nf}. However, the issue remains open across many systems, and both frameworks can fit the same experimental data with appropriate parameterization.

What can be said with certainty is that the \emph{mathematical} advantages of the logistic formulation---$C^\infty$ regularity, finite Lipschitz constant, nonzero basal expression, decoupled parameters---hold regardless of which sensitivity interpretation is adopted, and these advantages directly affect the reliability of numerical simulation, the tractability of bifurcation analysis, and the feasibility of observer and control design.

\subsection{Basal Expression and Non-Zero Output}
\label{sec:basal_expression}

Biological systems rarely exhibit true zero expression. Even under full repression, promoters display leakiness, stochastic transcription initiation events producing 0.1--1 mRNA molecules per cell in bacteria \cite{ozbudak2004multistability}, or 1--10 transcripts per million (TPM) in mammalian cells \cite{fantom2014promoter}. This persistent basal activity serves critical functions:
\begin{itemize}
\item It reduces phenotypic noise by shifting expression distributions from multimodal to unimodal.
\item It prevents irreversible trapping in off-states during stochastic fluctuations.
\item It enables rapid induction responses by maintaining the transcriptional machinery in a primed state.
\end{itemize}

The logistic function captures this fundamental biology naturally through its non-zero output at $x = 0$:
\[
f^+(0, \theta, \lambda) = \frac{1}{1 + e^{\lambda \theta}} > 0,
\]
which can be tuned via the product $\lambda \theta$ without introducing additional parameters. 

In contrast, Hill functions yield exactly zero at zero input: $h^+(0, \theta, n) = 0$. This forces modelers to add artificial constants:
\[
h_{\text{modified}}^+(x, \theta, n) = \varepsilon + (1 - \varepsilon) \frac{x^n}{x^n + \theta^n},
\]
where $\varepsilon > 0$ represents basal expression. However, this modification introduces several problems:

\begin{itemize}
    \item \textbf{Compresses dynamic range:} The range of $h_{\text{modified}}^+$ is $[\varepsilon, 1]$ rather than $[0, 1]$: the minimum value is lifted to $\varepsilon > 0$ and all intermediate values are scaled by $(1-\varepsilon)$, compressing the effective dynamic range and complicating quantitative interpretation of expression levels.
    
    \item \textbf{Adds an extra parameter:} The system now has $(\theta, n, \varepsilon)$ instead of $(\theta, n)$, increasing the parameter space dimension and reducing identifiability from limited data.
    
    \item \textbf{Obscures mechanistic interpretation:} $\varepsilon$ has no clear biophysical meaning (it is not a dissociation constant or cooperativity coefficient), making it difficult to relate to experimental observables.
    
    \item \textbf{Complicates derivative calculations:} The modified derivative is $h'^+_{\text{modified}}(x) = (1-\varepsilon) h'^+(x)$, introducing the scaling factor throughout stability analysis and control design.
\end{itemize}

\begin{figure}[htb!]
    \centering
    \includegraphics[height=0.5\linewidth]{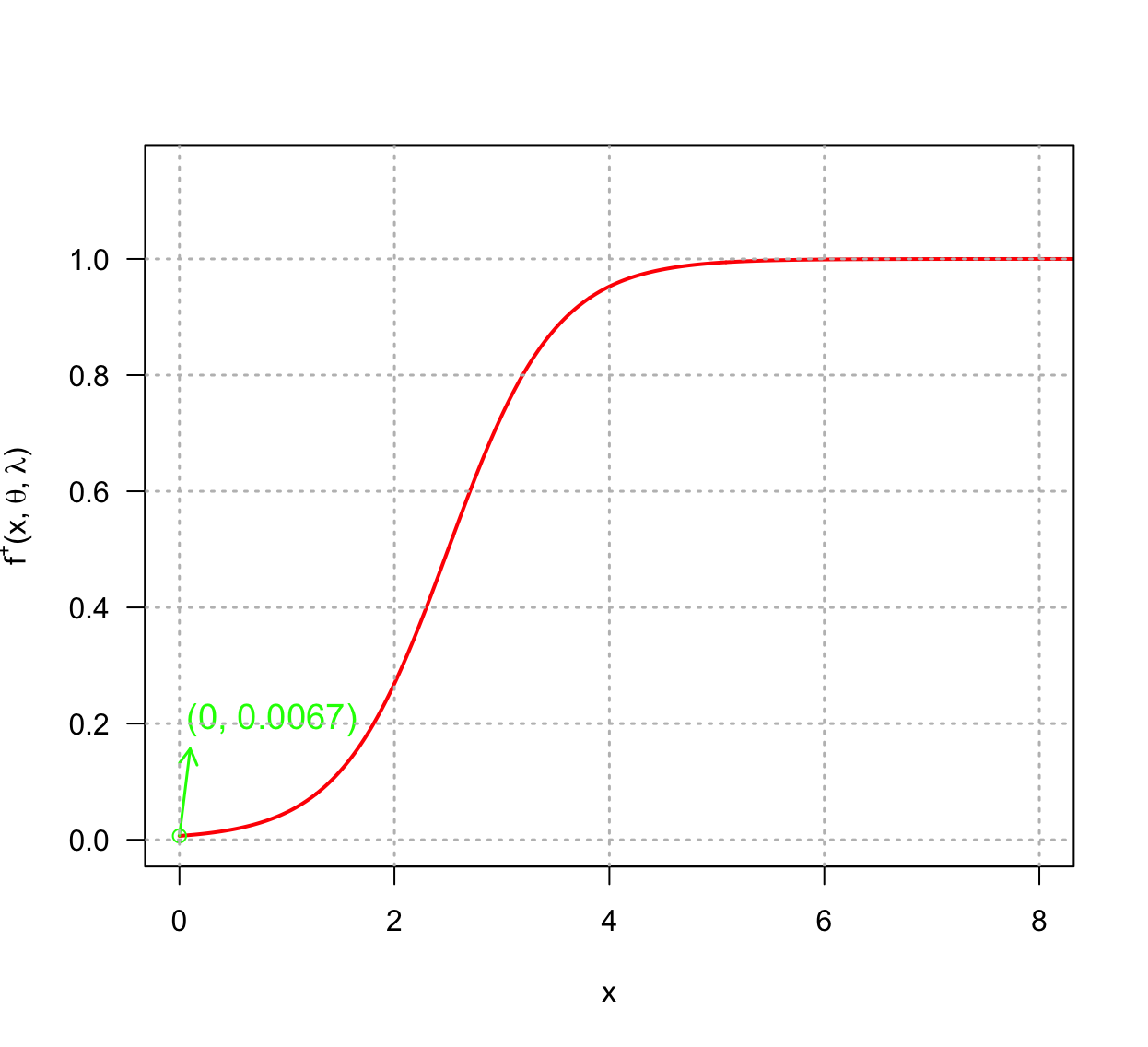}
    \caption{The logistic function $f^+(x, \theta, \lambda) =
    \frac{1}{1+e^{-\lambda(x-\theta)}}$ with $\lambda = 2$ and
    $\theta = 2.5$. Unlike the Hill function $h^+(x,\theta,n) =
    \frac{x^n}{x^n+\theta^n}$, which vanishes at $x=0$, the logistic
    ensures $f^+(0)>0$, naturally representing basal transcription.}
    \label{fig:logistic_hill}
\end{figure}

The logistic function, by contrast, incorporates basal activity as an
inherent feature. For the activation form
$f^+(x,\theta,\lambda) = \frac{1}{1+e^{-\lambda(x-\theta)}}$, the output
at $x=0$ is:
\begin{equation}
    f^+(0, \theta, \lambda) = \frac{1}{1 + e^{\lambda\theta}} > 0.
    \label{eq:basal_logistic}
\end{equation}
This value is nonzero (and small when the product $\lambda\theta$ is large), reflecting the minimum transcription level
even in the absence of inducers. By tuning $\lambda$ and $\theta$, modelers
can adjust this baseline: larger $\lambda\theta$ products drive $f^+(0)$
toward zero (approaching Hill-like behavior), while smaller products yield
higher basal rates. Crucially, this flexibility requires no additional
parameters; the basal level emerges directly from the function's shape,
determined solely by the steepness and threshold that already define the
regulatory response.

\begin{remark}[Basal expression range: leaky and constitutive regimes]
\label{rem:basal_range}
A natural question is whether high basal expression produced by the logistic
function is biologically realistic, or whether it signals a model
deficiency. The answer depends entirely on the system being modeled, and
this is precisely where the logistic framework demonstrates its versatility.

Consider the two-gene network of Vinoth~et~al.~\cite{vinoth2025extreme},
whose basal production rates are $g_A = g_B = 50\,\text{nM/min}$ — values
that are deliberately large, reflecting a strongly constitutively active
promoter architecture. In this regime, high basal expression is not a
modeling artifact but a biological fact: the genes in question are never
near silence. The logistic function handles this naturally, since the product
$\lambda\theta$ can be set small enough to place $f^+(0)$ close to
$1/(1+e^1) \approx 0.27$ or higher, correctly representing a promoter that
retains substantial activity even without its activator. Conversely, for the
lac or gal operons in \textit{E.~coli}, where leakiness is on the order of
$0.1$--$1\%$ of maximal expression~\cite{ozbudak2004multistability}, a
large $\lambda\theta$ product drives $f^+(0)$ toward zero. The same
functional form, with the same parameters $\lambda$ and $\theta$, thus
covers the full biological spectrum from tightly repressed promoters to
constitutively active genes, without any structural modification. This
stands in sharp contrast to Hill functions, which are structurally committed
to zero basal activity and require ad hoc additive offsets to enter either
regime.
\end{remark}

\textbf{Why basal expression matters.}
The biological importance of basal transcription is well-documented. In
bacteria such as \textit{E.~coli}, the \textit{lac} operon exhibits
measurable leakiness even under full repression, producing approximately
$0.1$--$1$ mRNA molecules per cell~\cite{ozbudak2004multistability}. This
low-level expression enables stochastic induction and bistable switching,
allowing cells to respond rapidly when lactose becomes available. Mammalian
housekeeping genes maintain moderate constitutive activity, typically
$10$--$100$ transcripts per million (TPM), ensuring continuous support for
essential cellular functions~\cite{fantom2014promoter}. Synthetic promoters
in engineered systems display controllable basal expression ranging from
$1$--$10\%$ of maximal induction, depending on promoter architecture and
regulatory design~\cite{lutz1997independent,kelly2009measuring}. In
synthetic oscillators, basal leakage at $5$--$20\%$ of maximum expression
can stabilize or destabilize dynamics depending on the circuit
architecture~\cite{joanito2020basal}.

Research by Flouriot~et~al.~\cite{flouriot2020basal} demonstrates that
basal expression correlates with chromatin accessibility and is essential
for controlling cellular differentiation, particularly in pluripotent stem
cells, where high basal activity facilitates reprogramming. Lu~et~al.~\cite{lu2018characterization}
showed that basal expression patterns can confound responses to
environmental cues, emphasizing the need to account for temporal and
context-dependent fluctuations. In cancer cell lines, Rees~et~al.~\cite{rees2016correlating}
found that basal expression profiles correlate with drug sensitivity,
revealing mechanistic insights that would be obscured if basal levels were
ignored.

The work of Joanito~et~al.~\cite{joanito2020basal} is particularly
illustrative. They examined synthetic repressilator-like circuits and
demonstrated that small basal leakage, modeled as imperfect repression or
cellular noise, can either destabilize or stabilize oscillatory dynamics,
depending on the regulatory architecture. In pure transcriptional models,
even modest basal leakage ($5\%$ of the maximum) disrupts oscillations by
pushing the system toward stable equilibria, mimicking the trapping of
low-expression states due to high degradation. However, adding
post-translational controls allows circuits to tolerate higher leakage while
sustaining oscillations, underscoring the need for minimal non-zero
production to maintain responsiveness in noisy environments. A complementary
quantitative analysis is given in the companion paper~\cite{belgacem2026logistic2},
where biophysically grounded \textit{E.~coli} parameters for the galactose
positive-autoregulation circuit yield a noise-driven off-state escape time
of approximately 44~minutes for the logistic formulation, whereas the
parameter-matched Hill formulation remains permanently confined near zero
with no intrinsic recovery mechanism.

\begin{example}[Two-node network: mutual activation and self-repression
as a concrete illustration of basal expression coverage]
\label{ex:two_node}

To make the above remarks concrete, we examine the two-gene regulatory
network of Vinoth~et~al.~\cite{vinoth2025extreme}, which exhibits
periodic oscillations, deterministic chaos, and extreme events —
abrupt, large-amplitude spikes in protein concentration reminiscent of
``rogue waves'' in nonlinear dynamics. The network comprises two proteins,
$A$ and $B$, governed by:
\begin{itemize}
    \item \textbf{Mutual activation:} Protein $A$ activates gene $B$,
          and protein $B$ activates gene $A$.
    \item \textbf{Cooperative self-repression:} Each gene represses its
          own expression via multimeric binding (e.g.\ dimers or
          tetramers), traditionally captured by a Hill coefficient $n>1$.
    \item \textbf{Time delays:} The system incorporates delays
          $\tau_1,\tau_2,\tau_{12},\tau_{21}$ representing
          transcription/translation lags or signalling delays.
\end{itemize}

These delays interact with the cooperative nonlinearities to produce a
rich bifurcation structure, leading to bistability, oscillatory bursting,
and even chaotic dynamics, in line with related observations of complex periodic and chaotic behavior in low-dimensional regulatory ring circuits~\cite{belgacem2025chaos,farcot2019chaos}. The regulatory architecture is depicted in
Figure~\ref{fig:two_gene_regulation}.

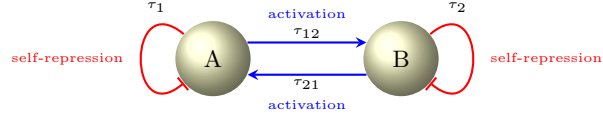
\begin{figure}[h]
\centering
\begin{tikzpicture}[font=\small, >=stealth, node distance=2.5cm]
  \node[circle, shading=ball, ball color=yellow!20,
        minimum size=1cm] (A) {$A$};
  \node[circle, shading=ball, ball color=yellow!20,
        minimum size=1cm, right of=A] (B) {$B$};
  \draw[-|, thick, red]
      (A) .. controls +(140:1.5) and +(220:1.5) .. (A)
      node[midway, left, xshift=-0.1cm] {\tiny self-repression};
  \draw[-|, thick, red]
      (B) .. controls +(40:1.5) and +(320:1.5) .. (B)
      node[midway, right, xshift=0.1cm] {\tiny self-repression};
  \draw[->, thick, blue]
      (A.25) -- (B.155) node[midway, above=0.2cm] {\tiny activation};
  \draw[->, thick, blue]
      (B.205) -- (A.335) node[midway, below=0.2cm] {\tiny activation};
  \node[above left=0.2cm of A]         {\tiny $\tau_1$};
  \node[above right=0.2cm of B]        {\tiny $\tau_2$};
  \node[above=0.15cm of $(A)!0.5!(B)$] {\tiny $\tau_{12}$};
  \node[below=0.15cm of $(A)!0.5!(B)$] {\tiny $\tau_{21}$};
\end{tikzpicture}
\caption{Regulatory topology of the two-gene network. Red bars indicate
cooperative self-repression; blue arrows represent mutual activation.
Delays $\tau_1,\tau_2$ govern self-repression lags; $\tau_{12},\tau_{21}$
govern cross-activation delays.}
\label{fig:two_gene_regulation}
\end{figure}

\subsubsection*{Original formulation}
In~\cite{vinoth2025extreme} the model combines additive linear activation
with Hill-type repression:
\begin{align}
    \frac{dA}{dt}
    &= \Bigl(g_A + g_{AB}\,B(t-\tau_{12})\Bigr)
       \frac{1}{1+\bigl(A(t-\tau_1)/A_0\bigr)^n}
       - \gamma_A A(t),
    \label{eq:vinoth_A}\\[4pt]
    \frac{dB}{dt}
    &= \Bigl(g_B + g_{BA}\,A(t-\tau_{21})\Bigr)
       \frac{1}{1+\bigl(B(t-\tau_2)/B_0\bigr)^n}
       - \gamma_B B(t),
    \label{eq:vinoth_B}
\end{align}
with $g_A=g_B=50\,\text{nM/min}$, $g_{AB}=g_{BA}=3\,\text{nM/min}$,
$\gamma_A=0.20\,\text{min}^{-1}$, $\gamma_B=0.24\,\text{min}^{-1}$,
$A_0=B_0=100\,\text{nM}$.

The basal rates $g_A=g_B=50\,\text{nM/min}$ are large by design:
these promoters are constitutively active at half-maximal rate even in
the complete absence of the cross-activator. This is precisely the
regime identified in Remark~\ref{rem:basal_range} where a small
$\lambda\theta$ product is biologically appropriate, and the logistic
model captures it without any structural modification.

\subsubsection*{Logistic reformulation}
We replace the Hill repression term $1/(1+(A/A_0)^n)$ by the decreasing
logistic $f^-(A,A_0,\lambda_3)$, and the unbounded linear activation
$(g_A + g_{AB}B)$ by a bounded increasing logistic scaled by $\kappa_1$:
\begin{align}
    \dot{A}(t)
    &= \kappa_1\,
       f^+\!\bigl(B(t-\tau_{12}),\,\theta_B,\,\lambda_1\bigr)\,
       f^-\!\bigl(A(t-\tau_1),\,A_0,\,\lambda_3\bigr)
       - \gamma_A A(t),
    \label{eq:A_logistic}\\[4pt]
    \dot{B}(t)
    &= \kappa_2\,
       f^+\!\bigl(A(t-\tau_{21}),\,\theta_A,\,\lambda_2\bigr)\,
       f^-\!\bigl(B(t-\tau_2),\,B_0,\,\lambda_4\bigr)
       - \gamma_B B(t),
    \label{eq:B_logistic}
\end{align}
or, written out explicitly:
\begin{align}
    \dot{A}(t)
    &= \frac{\kappa_1}
            {1+e^{-\lambda_1(B(t-\tau_{12})-\theta_B)}}
       \cdot
       \frac{1}{1+e^{\,\lambda_3(A(t-\tau_1)-A_0)}}
       - \gamma_A A(t),
    \label{eq:A_explicit}\\[4pt]
    \dot{B}(t)
    &= \frac{\kappa_2}
            {1+e^{-\lambda_2(A(t-\tau_{21})-\theta_A)}}
       \cdot
       \frac{1}{1+e^{\,\lambda_4(B(t-\tau_2)-B_0)}}
       - \gamma_B B(t).
    \label{eq:B_explicit}
\end{align}

\subsubsection*{Basal expression in this system}
The basal production rate of gene $A$ under the logistic formulation,
evaluated at zero activator concentration ($B=0$), is
\begin{equation}
    \dot{A}\big|_{B=0}^{\text{basal}}
    = \frac{\kappa_1}{1+e^{\lambda_1\theta_B}}
      \cdot f^-(A,A_0,\lambda_3).
    \label{eq:Abasal}
\end{equation}
This quantity is strictly positive for all finite $\lambda_1$ and
$\theta_B$, in agreement with the large basal rate $g_A=50\,\text{nM/min}$
of the original model. Far from being a deficiency of the logistic
formulation, a high basal value here is the correct biological answer:
these genes are constitutively active promoters, not leaky ones. The
logistic model naturally accommodates both regimes — tightly repressed
(large $\lambda\theta$) and constitutively active (small $\lambda\theta$)
— within a single unified formalism.

\subsubsection*{Parameter identification for this system}
The logistic parameters $\kappa_1,\lambda_1,\theta_B$ (and their
counterparts for gene $B$) can be estimated simultaneously from time-series
data by minimizing
\begin{equation}
    \min_{\kappa_i,\,\gamma_i,\,\lambda_{ij},\,\theta_{ij}}
    \sum_{t}
    \bigl\|\mathbf{x}(t)
    - \hat{\mathbf{x}}(t;\kappa_i,\gamma_i,\lambda_{ij},\theta_{ij})
    \bigr\|_2^2.
    \label{eq:ls_example}
\end{equation}
A key observation is that the observed basal expression level — the
value of $A(t)$ or $B(t)$ measured under null-activator conditions or at
$t=0$ before induction — provides a natural anchor for this system. Under
null-activator conditions and the approximation $f^-(A_{\mathrm{basal}},A_0,\lambda_3)\approx 1$
(valid when $\lambda_3(A_0-A_{\mathrm{basal}})\gg 1$, which holds in particular when
$\lambda_3 A_0 \gg 1$ and $A_{\mathrm{basal}}\ll A_0$, as is typical at basal levels),
\eqref{eq:Abasal} reduces to a single scalar equation relating $\kappa_1$,
$\lambda_1$, and $\theta_B$:
\begin{equation}
    r_{A,\mathrm{basal}}
    := \frac{\kappa_1}{1+e^{\lambda_1\theta_B}},
    \label{eq:basal_constraint}
\end{equation}
where $r_{A,\mathrm{basal}}$ (units: nM/min) is the basal production rate of gene $A$
inferred from the data (distinguishing it from the steady-state concentration $A^*$,
which satisfies $\kappa_1 f^+(A^*_B)\,f^-(A^*,A_0,\lambda_3)=\gamma_A A^*$ at equilibrium).
Equation~\eqref{eq:basal_constraint} enters the least-squares objective as an additional data
point rather than as an analytical formula from which one parameter is solved
in terms of the others. This is an important distinction from the approach
in the companion paper~\cite{belgacem2026logistic}, where the threshold $\theta$ is itself
derived from the basal expression: specifically, $\theta$ is set at the
concentration at which gene expression first exceeds the basal level,
so that the basal measurement determines $\theta$ directly, and $\lambda$
is then obtained using $\theta$.
In practice, all three quantities $\kappa_1$, $\lambda_1$, $\theta_B$
contribute distinct, independently observable signatures to the trajectory shape:
$\kappa_1$ sets the ceiling of the production rate, $\theta_B$ positions
the sigmoid's inflection point along the concentration axis, and $\lambda_1$
governs the steepness of the transition. The basal measurement at $B=0$
acts as an \emph{anchor constraint} that breaks the degeneracy among these
three parameters, making the identification problem better conditioned than
it would be from dynamic data alone. This simultaneous identifiability is
a structural advantage of the logistic parameterization that Hill-linear
hybrid models — where the basal rate $g_A$ and the Hill parameters
$(A_0, n)$ enter different functional forms and cannot be jointly anchored
by a single basal observation, do not share.
\end{example}

\textbf{Comparison with Hill functions.}
Standard Hill functions lack this flexibility entirely. The activation
Hill function $h^+(x,\theta,n)=x^n/(x^n+\theta^n)$ yields exactly zero
at $x=0$. One might attempt to remedy this by adding a constant offset:
$h^+(x,\theta,n)+\varepsilon$. But this modification is awkward and
problematic: $\varepsilon$ is arbitrary, lacking biological justification;
it shifts the entire response curve upward, violating the normalized
$[0,1]$ range for large $x$; and it adds a parameter with no clear
biophysical meaning, reducing identifiability from limited data.
To maintain normalization, Hill modelers must instead use
\[
h_{\text{modified}}^+(x,\theta,n)
= \varepsilon + (1-\varepsilon)\frac{x^n}{x^n+\theta^n},
\]
which further complicates derivative calculations
($h'^+_{\text{modified}}=(1-\varepsilon)h'^+$) throughout stability
analysis and control design. The logistic function requires none of these
modifications: its basal output is governed by the same parameters
$\lambda$ and $\theta$ that define the regulatory response, maintaining
the $[0,1]$ range without structural alteration and covering the full
spectrum from leaky to constitutively active promoters within a single
unified form.

\begin{figure}[htb!]
    \centering
    \begin{subfigure}{0.46\linewidth}
        \centering
        \includegraphics[width=\linewidth]{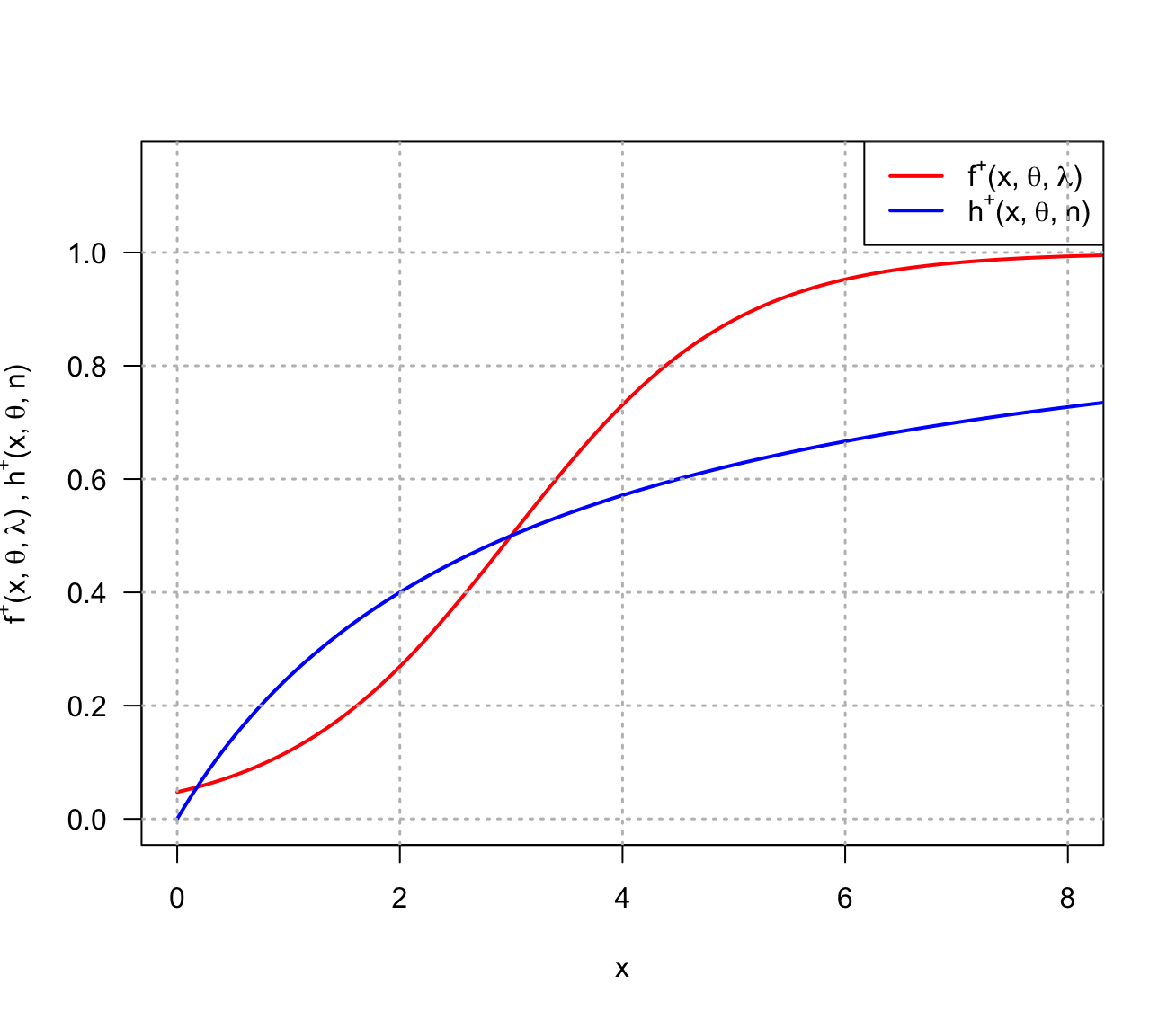}
        \caption{$\lambda = n = 1$}
    \end{subfigure}
    \hfill
    \begin{subfigure}{0.46\linewidth}
        \centering
        \includegraphics[width=\linewidth]{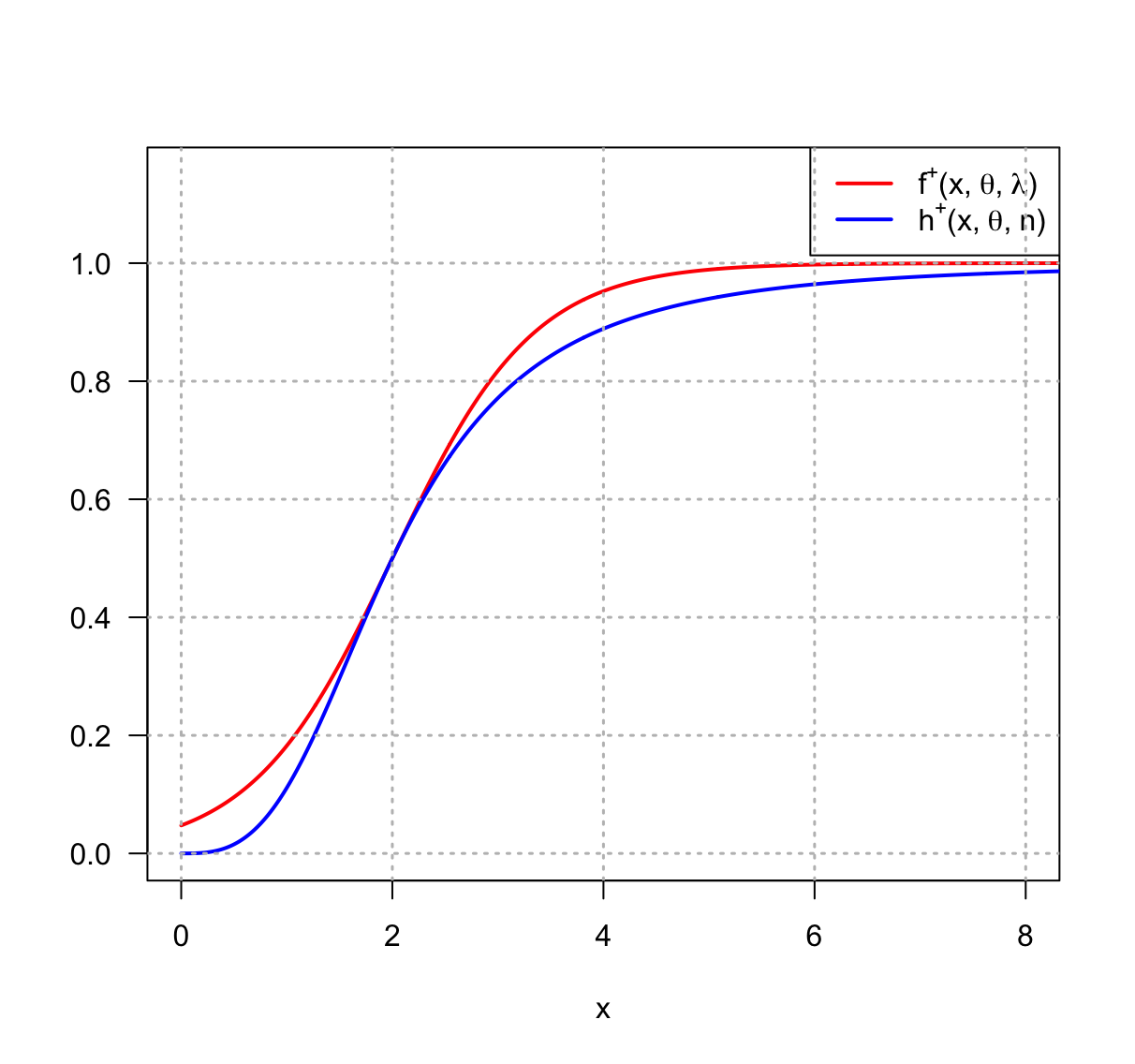}
        \caption{$n=3,\;\theta=2,\;\lambda=\tfrac{3}{2}$}
    \end{subfigure}
    \hfill
    \begin{subfigure}{0.46\linewidth}
        \centering
        \includegraphics[width=\linewidth]{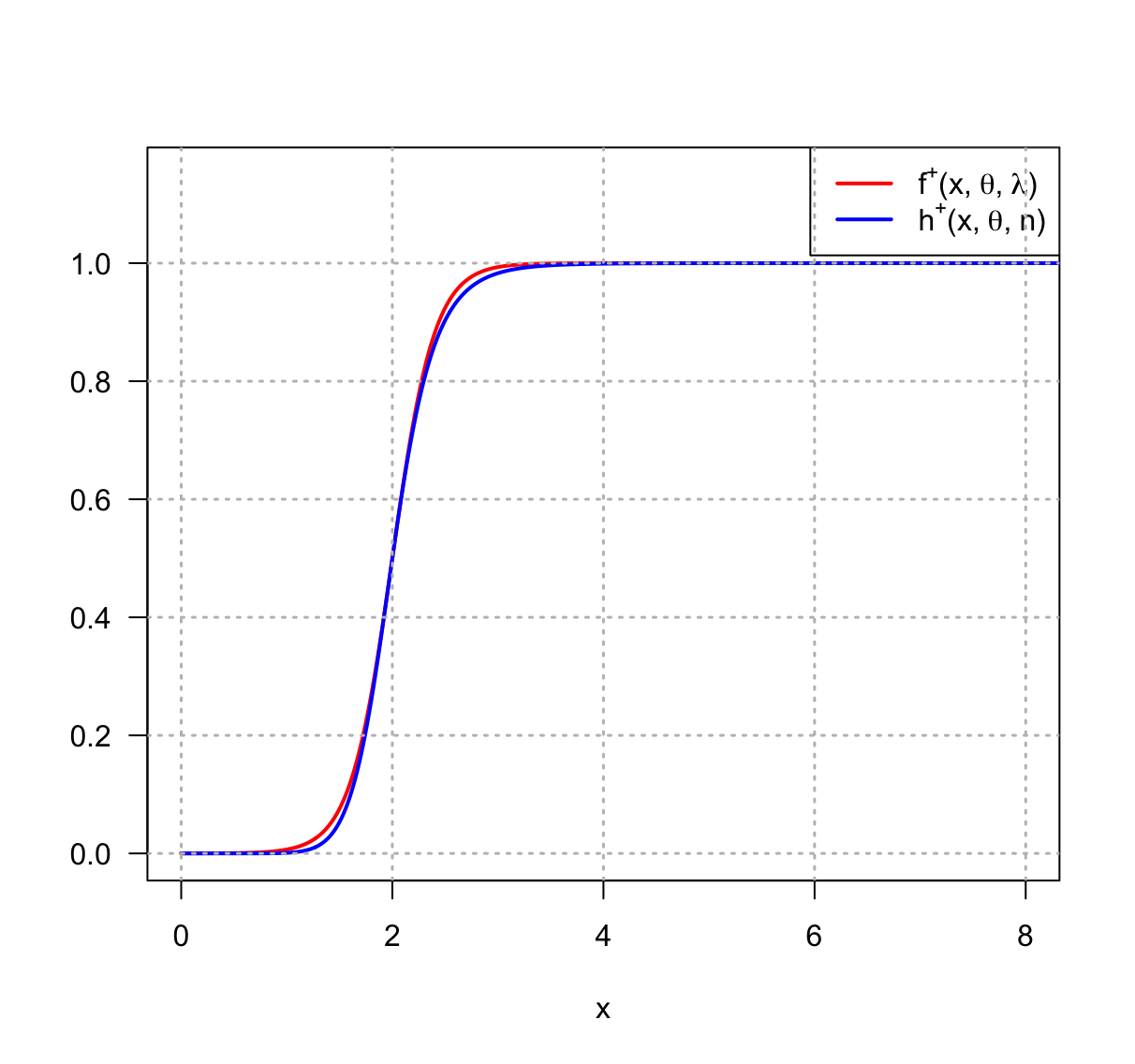}
        \caption{$n=10,\;\theta=2,\;\lambda=\tfrac{n}{\theta}=5$}
    \end{subfigure}
    \caption{Comparison of logistic and Hill activation functions with
    increasing steepness. As $\lambda$ increases, $f^+(x,\theta,\lambda)$
    approximates the steepness of $h^+(x,\theta,n)$ while maintaining
    $f^+(0)>0$, capturing basal transcription without artificial offsets.}
    \label{fig:logistic_hill_comparison}
\end{figure}

For repression, the decreasing logistic $f^-(0) = 1/(1+e^{-\lambda\theta})$ provides a
tunable, unrepressed expression that approaches but does not rigidly equal unity,
reflecting polymerase saturation and competition for transcriptional resources.
The decreasing Hill function, by contrast, enforces $h^-(0)=1$ unconditionally,
overstating achievable expression levels in resource-limited cellular environments
and requiring artificial offsets when modeling realistic imperfections.
The parameter regime consequences of this asymmetry are treated quantitatively
in Section~\ref{sec:logistic_regime_considerations}, and a full analysis of repression
dynamics in stochastic environments and the scaled logistic variant appears in
Section~\ref{sec:repression_noisy}.

\subsection{Parameter Regime Considerations: Basal Expression and Repression Bounds}
\label{sec:logistic_regime_considerations}

A complete assessment of the logistic framework requires acknowledging two parameter-regime-dependent behaviors that modelers should consider explicitly when selecting $\lambda$ and $\theta$. We show here that both behaviors admit direct biological interpretations and are fully controllable through the single parameter product $\lambda\theta$, which is itself anchored to experimentally measurable quantities.

\subsubsection{Basal expression and the symmetry of the logistic function}

The logistic function is symmetric about its inflection point $x = \theta$, so the
basal output is governed entirely by the product $\lambda\theta$
(Section~\ref{sec:basal_expression}, equation~\eqref{eq:basal_logistic}).
Concretely: $\lambda\theta = 6$ gives $f^+(0) = 1/(1+e^6) \approx 0.247\%$,
appropriate for tightly repressed promoters such as the \textit{lac} or \textit{gal}
operons; $\lambda\theta = 1$ gives $f^+(0) \approx 26.9\%$, appropriate for
constitutively active promoters. Modelers should select $\lambda$ and $\theta$ so
that equation~\eqref{eq:basal_logistic} matches the experimentally measured basal
transcription level of the system under study~\cite{joanito2020basal,ozbudak2004multistability}.

\subsubsection{Upper bound of the repression function at zero repressor}

The decreasing logistic satisfies
\begin{equation}\label{eq:repression_at_zero}
f^-(0, \theta, \lambda) = \frac{1}{1 + e^{-\lambda\theta}} < 1
\end{equation}
for all finite positive $\lambda$ and $\theta$, in contrast to the decreasing Hill function, which attains exactly $h^-(0) = 1$. For applications requiring full expression saturation in the total absence of repressor, this sub-unity value may appear restrictive. However, it reflects a biologically grounded reality: even under fully unrepressed conditions, polymerase saturation, competition for transcriptional machinery and ribosomes, stochastic promoter switching, and chromatin accessibility fluctuations prevent deterministic maximal expression~\cite{frank2013input,reeve2013pharmacodynamic,kelly2009measuring}. Experimental studies consistently report that strong promoters achieve $70$--$99\%$ of their theoretical maximum under typical cellular conditions~\cite{madar2011negative,oehler1994quality}, values naturally representable by~\eqref{eq:repression_at_zero} through appropriate parameter selection. Indeed, the Hill function's rigid assumption $h^-(0) = 1$ is itself a modeling idealization that overstates achievable expression levels in resource-limited or noisy cellular environments~\cite{hernandez2023corrected}.

When exact normalization to unity is nonetheless required for theoretical consistency or comparative analysis with Hill-based models, the scaled repression variant
\begin{equation}\label{eq:scaled_fminus_section5}
f^-_{\mathrm{scaled}}(x, \theta, \lambda) = \bigl(1 + e^{-\lambda\theta}\bigr) \cdot f^-(x, \theta, \lambda)
\end{equation}
restores $f^-_{\mathrm{scaled}}(0) = 1$ exactly while preserving global $C^\infty$ smoothness, the sigmoidal shape, and all analytical tractability of the logistic framework. Under typical parameter regimes ($\lambda\theta \geq 4$), the scaling factor $(1 + e^{-\lambda\theta})$ deviates from unity by less than $2\%$~\cite{gottschalk2005five}, making the scaled and unscaled forms practically indistinguishable in model dynamics. The choice between forms is therefore guided entirely by biological context: the unscaled form is preferable when modeling resource-limited, stochastic, or noisy cellular environments where $f^-(0) < 1$ accurately captures imperfect unrepressed expression, while the scaled form is appropriate for comparative benchmarking against Hill-based formulations or for theoretical frameworks requiring exact boundary conditions.

\subsubsection{Summary}

Neither behavior constitutes a fundamental flaw of the logistic framework. Both are consequences of a single interpretable parameter product $\lambda\theta$ that maps directly to measurable biological quantities---the basal transcription level for activation~\eqref{eq:basal_logistic} and the maximal unrepressed expression level for repression~\eqref{eq:repression_at_zero}---and both can be addressed by principled parameter selection grounded in experimental measurements, or, in the case of repression, by the scaled variant~\eqref{eq:scaled_fminus_section5}. This explicit and tunable control over boundary behavior---entirely absent from Hill functions, which are structurally committed to $h^+(0) = 0$ and $h^-(0) = 1$ regardless of biological context---represents a modeling advantage. We encourage modelers to report the chosen $\lambda\theta$ product alongside fitted parameter values, thereby making the implied basal expression level transparent and directly verifiable against experimental measurements~\cite{joanito2020basal,ozbudak2004multistability}.

\subsection{Preserving Cooperativity and Biological Fidelity}
Having demonstrated the mathematical, computational, and biological advantages of logistic functions, we address a critical question: Can logistic formulations preserve the cooperativity that Hill functions capture mechanistically?  The answer is affirmative, through careful parameter matching that preserves phenomenological behavior while improving analytical tractability.

\subsubsection{Parameter Matching Strategies}
The Hill coefficient $n$ quantifies the degree of cooperativity and is typically obtained by fitting Hill functions directly to gene expression data \cite{madar2011negative, oehler1994quality}. This cooperativity can be preserved in logistic formulations through the steepness parameter $\lambda$, which controls the sharpness of the sigmoidal transition analogously to $n$.

\subsubsection{Derivative matching at the half-maximal point}
We derive the matching condition $\lambda = n/\theta$ by equating the slopes of the Hill and logistic functions at their respective half-maximal points.

For the Hill activation function $h^+(x) = x^n/(x^n + \theta^n)$, the derivative at $x = \theta$ is:
\begin{align}
\left.\frac{d h^+}{dx}\right|_{x=\theta}
&= \left.\frac{n x^{n-1} \theta^n}{(x^n + \theta^n)^2}\right|_{x=\theta}
= \frac{n \theta^{n-1} \cdot \theta^n}{(\theta^n + \theta^n)^2}
= \frac{n \theta^{2n-1}}{4\theta^{2n}}
= \frac{n}{4\theta}.
\end{align}

For the increasing logistic function $f^+(x) = (1 + e^{-\lambda(x-\theta)})^{-1}$, the derivative at $x = \theta$ is (using \eqref{eq:deriv_fplus} and $f^+(\theta) = 1/2$):
\[
\left.\frac{d f^+}{dx}\right|_{x=\theta}
= \lambda f^+(\theta)\bigl(1 - f^+(\theta)\bigr)
= \lambda \cdot \frac{1}{2} \cdot \frac{1}{2}
= \frac{\lambda}{4}.
\]

Equating these slopes:
\begin{equation}\label{eq:param_match}
\frac{\lambda}{4} = \frac{n}{4\theta} \implies \boxed{\lambda = \frac{n}{\theta}},
\end{equation}
ensuring identical local responsiveness to concentration changes near the activation threshold. The same result holds for the decreasing case.

\subsubsection{Phenomenological Versus Mechanistic Interpretations}

A clarification is warranted concerning the relationship between logistic and Hill functions and the concept of mechanistic modeling. Both functions are almost always deployed as \emph{phenomenological} approximations in systems biology: the Hill function is formally derived from equilibrium binding theory for a single step with $n$ identical independent sites, but in practice it is routinely fitted to experimental dose-response data as a convenient sigmoidal curve with tunable steepness, without requiring the underlying molecular mechanism to conform to multi-site equilibrium binding~\cite{santillan2008use,bintu2005transcriptional}. Non-integer Hill coefficients, which are the norm in experimental fits, have no simple mechanistic interpretation in terms of integer binding sites; they serve as phenomenological descriptors of ultrasensitivity.

Logistic functions are similarly phenomenological in their standard deployment in GRN models. However, it should be noted that logistic-type responses can also be derived from first principles in appropriate settings. The thermodynamic framework of Bintu et al.~\cite{bintu2005transcriptional,bintu2005transcriptional_app}, which uses statistical mechanics to compute the probability that RNA polymerase occupies a promoter as a function of regulator concentrations, produces regulation factors of the general form
\[
F_{\mathrm{reg}} = \frac{\text{(sum of Boltzmann weights for states with RNAP bound)}}{\text{(sum of Boltzmann weights for states without RNAP bound)}},
\]
which for a simple repressor reduces to $F_{\mathrm{reg}} = 1/(1+[R]/K_R)$ (a Hill function with $n=1$), and for more complex architectures (cooperative repression, DNA looping, dual activators) gives rational polynomials of higher degree in the TF concentrations. Notably, $F_{\mathrm{reg}}$ for simple activation takes the form $(1 + [A] f/K_A)/(1+[A]/K_A)$, where $f$ is an enhancement factor; this can be written $f^+(\ln([A]/K_A), 0, 1)$ in the log-concentration variable, consistent with the logistic-of-log identity~\eqref{eq:hill_as_logistic_of_log}. More generally, combinations of activators and repressors in the thermodynamic model yield rational polynomials in $[TF]$ (see Table~1 of~\cite{bintu2005transcriptional} and Equations~1--2 of~\cite{bintu2005transcriptional_app}), providing a biophysical foundation for approximating each interaction term by a logistic factor and assembling them into the product-of-logistics formulation employed here.

The key distinction between the two frameworks, therefore, is not mechanistic versus phenomenological but rather mathematical structure. The Hill function's power-law form introduces derivative singularities for non-integer exponents, couples steepness and threshold in a single expression ($n/(4\theta)$), and lacks a closed-form inverse for general $n$. The logistic function avoids all these analytical pathologies while capturing the same sigmoidal switching behavior. The steepness parameter $\lambda$ plays the same operational role as the Hill coefficient $n$ (quantifying effective cooperativity from dose-response data) but does so without creating smoothness restrictions that compromise downstream analysis. Both parameters are ultimately fitted to experimental data and should be understood as summaries of ultrasensitivity mechanisms that may include cooperative binding, allosteric effects, phosphorylation cascades, competition for limiting resources, or other cellular nonlinearities.

\subsubsection{When to prefer which framework}
Hill functions retain their natural connection to thermodynamic parameters in settings where the mechanistic detail of equilibrium binding constants is the primary scientific question. In settings where analytical tractability, scalable inference, or real-time control are the primary goals, the logistic framework provides the same sigmoidal switching behavior without the smoothness restrictions imposed by non-integer exponents.

\subsection{Repression Dynamics in Noisy Environments and Maximal Expression Scaling}
\label{sec:repression_noisy}

Repression functions suppress gene expression as the concentration of repressors $x$ increases.
The biochemically grounded Hill repression function $h^-(x, \theta, n) = \frac{\theta^n}{x^n +
\theta^n}$ originates from equilibrium binding models and incorporates cooperativity through the
exponent $n$, which governs the steepness of the repression curve. By contrast, the logistic
repression function $f^-(x, \theta, \lambda) = \frac{1}{1 + e^{\lambda(x - \theta)}}$ delivers
a smoother, exponentially decaying profile with steepness controlled by the parameter $\lambda$.

The Hill function achieves absolute unity at $x = 0$, that is, $h^-(0, \theta, n) = 1$, encoding full gene expression in the absence of repressor. The logistic function, however, starts \textit{near} but not exactly at 1 when $x = 0$. For instance, with $\lambda = 4$ and $\theta = 3$,
\[
f^-(0, 3, 4) = \frac{1}{1 + e^{-\lambda\theta}} = \frac{1}{1 + e^{-12}} \approx 0.999994,
\]
very close to 1 but not quite reaching it. This slight deficit from unity is not a defect but a \textbf{biological feature}, capturing the reality that even without repressors, maximal expression rarely reaches 100\% due to polymerase saturation, resource competition (ribosomes, NTPs, ATP/GTP), transcriptional noise from stochastic bursting, and chromatin/promoter-accessibility fluctuations. The logistic's initial value depends on the product $\lambda\theta$ and is precisely tunable: increasing $\lambda\theta$ pushes $f^-(0)$ closer to 1, capturing prolonged yet imperfect expression in repressor-absent states, as observed experimentally in systems with promoter leakiness, polymerase limitations, or stochastic bursting~\cite{el2005stochastic}.

When repressors are absent, the logistic's adjustable asymptote near 1 enables fine-tuned representation of low-level constitutive activity, crucial in noisy or fluctuating environments. This capability enhances robustness in \textbf{stochastic simulations}, particularly those employing the Gillespie algorithm: the logistic's smooth, exponential form yields more stable propensity functions and minimizes variance in rate calculations at low copy-number regimes, where stochastic effects dominate and mean-field approximations break down. Hill functions, derived from mean-field assumptions that presuppose large molecular populations, frequently require stochastic corrections to account for fluctuations~\cite{hernandez2023corrected}. The logistic's exponential structure inherently mitigates these issues, integrating seamlessly with Gillespie algorithms. In \textbf{synthetic biology control design}, the logistic's smoothness facilitates stable feedback linearization and observer-based estimation, improving circuit predictability under uncertainty~\cite{del2016control,barajas2022synthetic,teo2015synthetic,sarpeshkar2014analog}; the continuous derivatives of all orders make it amenable to differential geometric control methods.

\textbf{Specific synthetic-biology architectures.} In \textbf{auto-regulatory circuits} (negative auto-regulation in bacterial quorum-sensing systems and metabolic pathway controllers), the logistic's smooth form minimizes propensity-rate variance during low-copy fluctuations, yielding more accurate steady-state distribution predictions. In \textbf{toggle switches} (mutual repression between two genes), Hill models may overestimate the sharpness of bistability thresholds and the separation between stable states under stochastic fluctuations; the logistic provides a more robust approximation that naturally accounts for resource availability and intrinsic noise. In \textbf{repressilators} (cyclic repression among three or more genes, e.g., the canonical LacI--TetR--$\lambda$CI circuit), Hill models can overestimate oscillation damping at low molecular counts ($\sim$10--100 copies per cell), whereas the logistic's resource-aware formulation better captures the sustained oscillations observed in single-cell microfluidic experiments~\cite{loinger2007stochastic,hernandez2023corrected}.

\subsubsection{Scaling for Maximal Gene Expression}

In certain scenarios, the absolute full expression $f^-(0, \theta, \lambda) = 1$ may be
desirable in $x = 0$ to match exactly the behavior of the Hill function. To address this
requirement, we can scale the logistic repression function to ensure full expression when no
repressor is present, perfectly aligning with $h^-(0, \theta, n) = 1$. The \textbf{scaled
logistic function} is defined as
\[
f^-_{\text{scaled}}(x, \theta, \lambda) = \big(1 + e^{-\lambda\theta}\big) \cdot \frac{1}{1 + e^{\lambda(x - \theta)}}.
\]

This scaling arises naturally from the requirement that $f^-_{\text{scaled}}(0, \theta, \lambda)
= 1$. Since the unscaled logistic satisfies
\[
f^-(0, \theta, \lambda) = \frac{1}{1 + e^{-\lambda\theta}},
\]
its reciprocal, $1 + e^{-\lambda\theta}$, provides exactly the normalization factor needed to
achieve unity at $x = 0$. As $x$ increases beyond zero, the denominator $1 + e^{\lambda(x -
\theta)}$ grows exponentially, driving the function toward zero while preserving the
characteristic sigmoidal shape and smooth decay properties. The scaling factor $e^{-\lambda\theta}$
diminishes rapidly with larger $\lambda$ and $\theta$ values. For instance:
\begin{itemize}
    \item With $\lambda = 4$ and $\theta = 3$: $e^{-\lambda\theta} = e^{-12} \approx 6.14 \times
    10^{-6}$, making $(1 + e^{-\lambda\theta}) \approx 1.0000061$.
    \item With $\lambda = 3$ and $\theta = 2$: $e^{-\lambda\theta} = e^{-6} \approx 0.00248$,
    giving $(1 + e^{-\lambda\theta}) \approx 1.00248$.
    \item With $\lambda = 1$ and $\theta = 2$: $e^{-\lambda\theta} = e^{-2} \approx 0.1353$,
    yielding $(1 + e^{-\lambda\theta}) \approx 1.1353$.
\end{itemize}

Thus, the scaling effect is noticeable only when $\lambda$ and $\theta$ are small. Under typical
parameter regimes used in gene regulatory network models ($\lambda \geq 2$, $\theta \geq 2$),
the scaling factor remains very close to 1, ensuring minimal deviation from the original logistic
shape. This means that for most practical applications with moderately steep transitions, the
scaled and unscaled logistic functions are nearly indistinguishable in shape, differing only by
a small vertical scaling that normalizes the maximum to exactly 1.

\begin{figure}[h]
\centering
\includegraphics[width=0.465\textwidth]{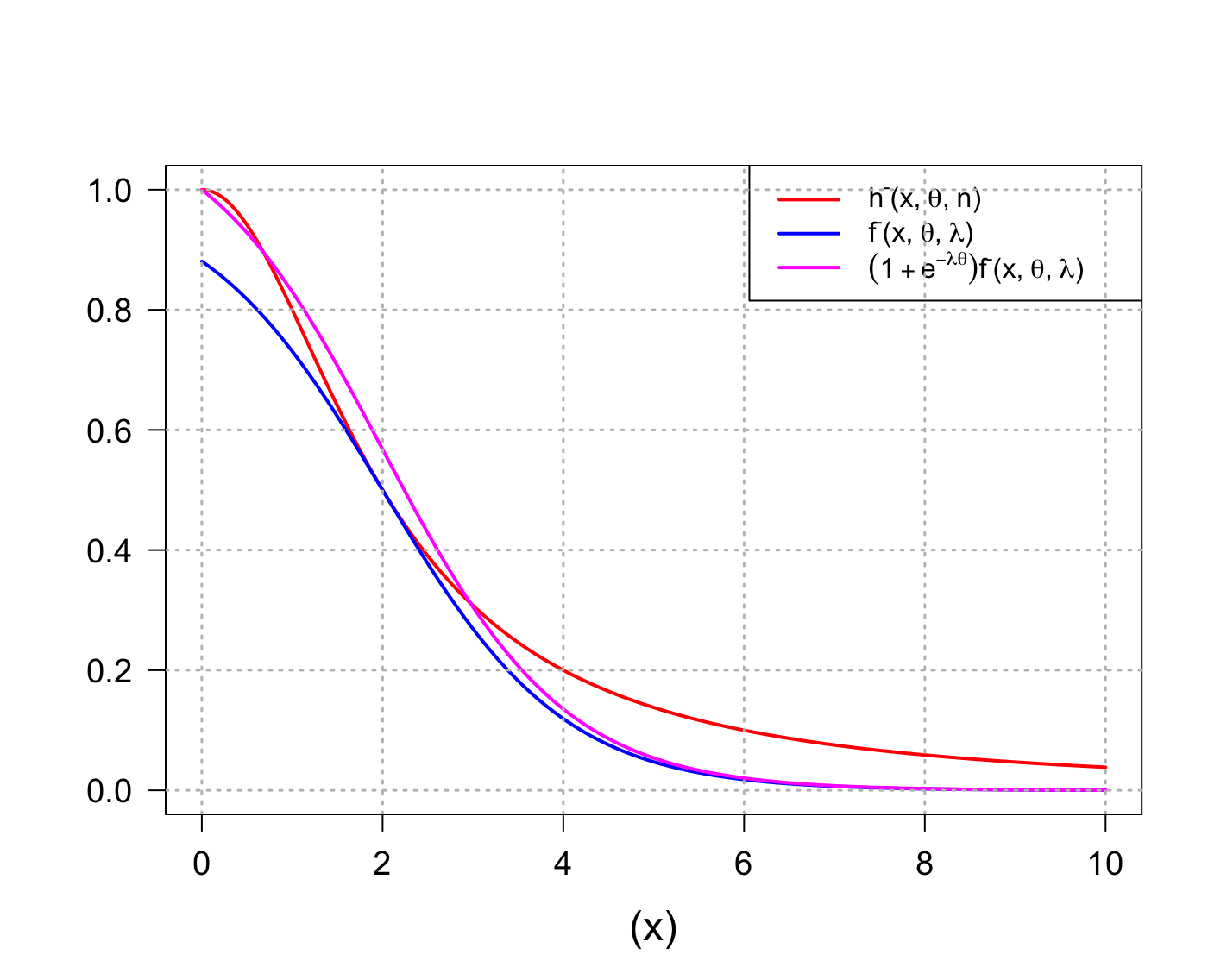}
\includegraphics[width=0.41\textwidth]{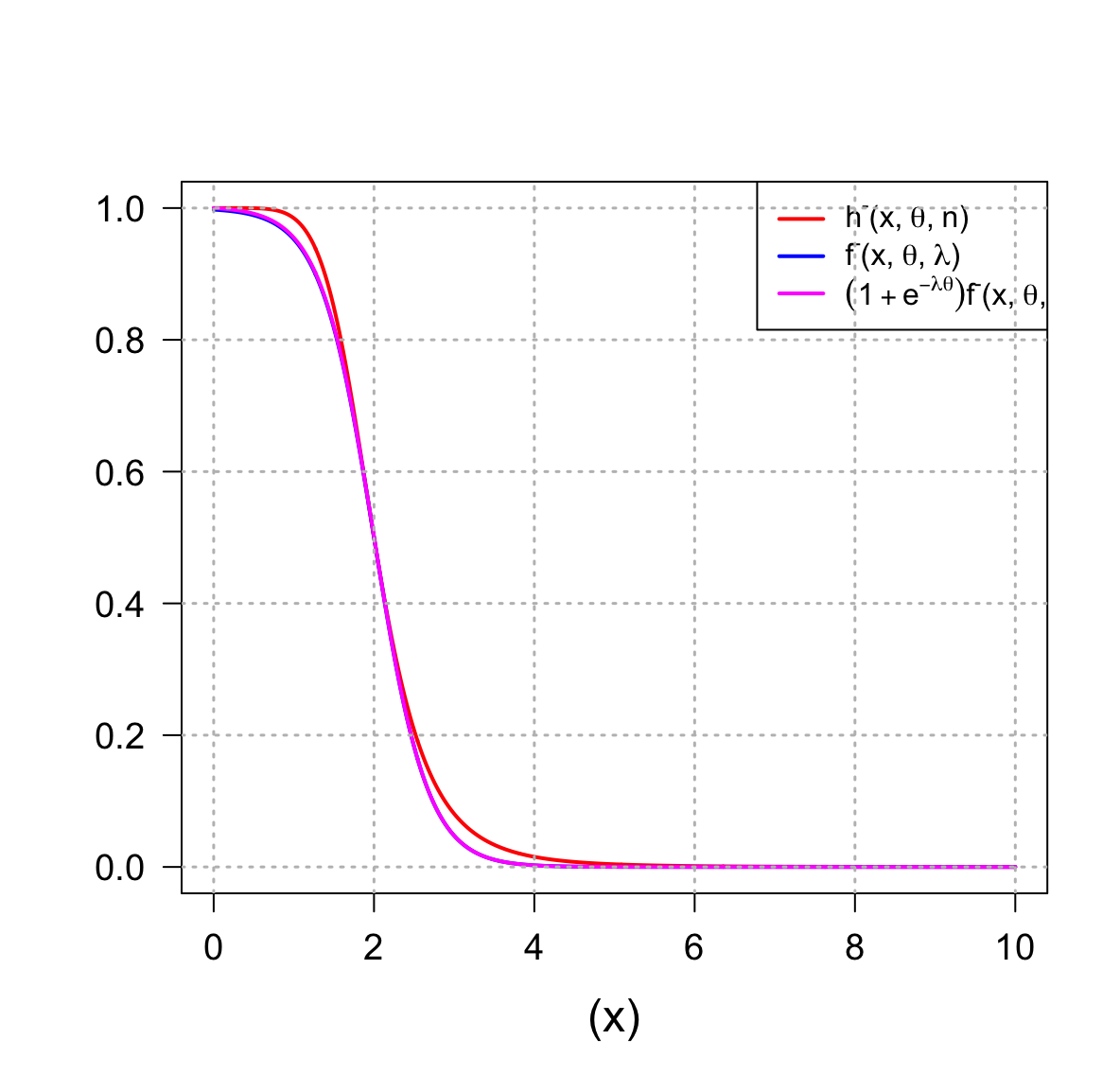}
\caption{Comparison of repression functions. \textbf{Left}: $n = 2$, $\theta = 2$, $\lambda =
1$. \textbf{Right}: $n = 6$, $\theta = 2$, $\lambda = 3$. Red line: Hill function $h^-(x,
\theta, n)$; blue line: logistic function $f^-(x, \theta, \lambda)$; magenta line: scaled
logistic $(1 + e^{-\lambda\theta})f^-(x, \theta, \lambda)$. The left panel shows moderate
baseline reduction ($\sim$88\% of maximum), suitable for systems with inherent inefficiencies
or resource limitations. The right panel shows near-perfect alignment among all three functions
when $\lambda$ and $\theta$ are large, with the scaling factor negligible.}
\label{fig:repression_functions}
\end{figure}

Figure~\ref{fig:repression_functions} compares the Hill, original logistic, and scaled logistic
repression functions across repressor concentrations from 0 to 10. The left panel employs
parameters $n = 2$, $\theta = 2$, and $\lambda = \frac{n}{\theta} = 1$; the right panel uses $n
= 6$, $\theta = 2$, and $\lambda = \frac{n}{\theta} = 3$, revealing how parameter variations
influence curve shapes and transition steepness. In the left panel ($\lambda = 1$, $\theta = 2$),
we compute:
\[
e^{-\lambda\theta} = e^{-2} \approx 0.1353 \quad \Rightarrow \quad f^-(0) \approx \frac{1}{1 +
0.1353} \approx 0.8808,
\]
representing approximately 88\% of maximum expression. This moderate reduction is well-suited
for modeling systems with inherent inefficiencies or resource limitations, such as bacterial
cells under nutrient stress, synthetic circuits approaching metabolic burden thresholds, or
promoters with significant constraints on basal activity. In the right panel ($\lambda = 3$,
$\theta = 2$), the scaling factor $(1 + e^{-6}) \approx 1.0025$ is negligible, and all three
functions align closely, demonstrating that the logistic and Hill functions exhibit nearly
identical behavior under typical steep-transition parameter regimes.

\begin{example}
Reconsidering the genetic oscillator model, in which two genes form a negative feedback loop
with gene 1 repressing gene 2 and gene 2 activating gene 1, we can integrate the scaled logistic
function to ensure exact unit maxima. The modified system becomes:
\begin{equation}
\begin{aligned}
\dot{x}_1 &= \kappa'_1 \frac{1}{1 + e^{\lambda(x_2 - \theta_2)}} - \gamma_1 x_1, \\
\dot{x}_2 &= \kappa_2 \frac{1}{1 + e^{-\lambda(x_1 - \theta_1)}} - \gamma_2 x_2,
\end{aligned}
\label{eq:scaled_oscillator_repression}
\end{equation}
where $\kappa'_1 = \kappa_1(1 + e^{-\lambda\theta_2})$. This adjustment scales the production
rate to ensure gene 1 achieves full expression when $x_2 = 0$ (no repressor present), while
gene 2 reaches full expression when $x_1$ is high (activator present). The scaled system closely
resembles the original; the primary differences lie in the constant multiplicative factor $(1 +
e^{-\lambda\theta_2})$ applied to the expression rate for repression. This factor enhances
consistency with systems that expect maximum output under optimal conditions, aligning the
logistic model with the Hill function's behavior at boundary points. Importantly, the oscillator
dynamics, including period, amplitude, and stability, remain fundamentally unchanged by this
scaling, as the factor affects only the absolute scale of expression of the repressed gene, not
the relative feedback strength that governs oscillatory behavior. The local stability analysis
presented above remains valid, with the trace and determinant of the Jacobian matrix retaining
their signs, ensuring the same qualitative dynamical behavior.
\end{example}

Scaling $f^-(x, \theta, \lambda)$ by $(1 + e^{-\lambda\theta})$ ensures absolute full expression
when repressor is absent, which is consistent with and often preferable for systems featuring
strong promoters where maximal output is expected under unrepressed conditions. It is important
to note that for activation Hill functions $h^+(x, \theta, n) = \frac{x^n}{x^n + \theta^n}$,
the full expression is \textit{never} achieved at any finite $x$; the function merely approaches
1 asymptotically as $x \to \infty$. Thus, the asymmetry between repression and activation,
where repression can achieve perfect repression at $x = 0$, but activation cannot achieve
perfect activation at any finite $x$, is inherent in both the Hill and logistic frameworks.

The logistic function provides an \textbf{exponentially decaying repression curve} that better mimics gradual, sigmoidal responses in cellular systems where repression builds progressively. The smoothness is adjustable via $\lambda$, allowing fine-tuned control over transition sharpness without relying exclusively on cooperativity (encoded by $n$ in the Hill equation), which may not always align with non-cooperative biological mechanisms, such as allosteric regulation, sequential binding without cooperativity, or indirect regulatory cascades.

The choice between scaled and unscaled logistic functions depends on the modeling context. The \textbf{unscaled logistic} $f^-(x,\theta,\lambda)$ starts slightly below 1 at low repressor levels, with this deviation accounting for real-world constraints (polymerase saturation, finite transcriptional machinery, resource limitations, transcriptional noise, and chromatin accessibility effects). It is preferable when modeling resource-limited, stochastic environments and when $\lambda\theta$ is moderate (e.g., $\lambda\theta < 4$). The \textbf{scaled logistic} $f^-_{\text{scaled}}(x,\theta,\lambda) = (1+e^{-\lambda\theta})f^-(x,\theta,\lambda)$ normalizes to exactly 1 at $x = 0$ while preserving the exponential decay shape, providing exact agreement with the Hill function's boundary behavior. Under typical parameter ranges ($\lambda\theta \geq 4$), the scaling factor satisfies $1 + e^{-\lambda\theta} \leq 1.02$, so the two forms are practically indistinguishable in dynamics; the scaled form is preferred when interfacing with Hill-based models or when comparative analysis requires exact normalization.

\textbf{Applications in Synthetic Biology.}
In synthetic circuits, tuning the unrepressed expression level is crucial for several reasons grounded in experimental observations:
\begin{itemize}
    \item \textbf{Metabolic flux optimization}: Overexpression of heterologous proteins diverts cellular resources from endogenous metabolism, reducing growth rate and product yield. The logistic's adjustable maximum allows modeling of realistic expression ceilings imposed by metabolic burden, enabling better predictions of sustainable production rates in industrial bioprocessing.

    \item \textbf{Toxicity avoidance}: Long-term expression in continuous bioreactors, therapeutic protein manufacturing, or biosensor deployment can become toxic if expression exceeds cellular capacity for protein folding, membrane integration, or metabolite processing. The unscaled logistic naturally incorporates a ``safety valve'' that prevents overestimation of sustainable expression levels.

    \item \textbf{Resource competition in multi-gene circuits}: In complex synthetic circuits (e.g., metabolic pathways with 5--10 enzymatic steps), incomplete repression leads to sub-maximal outputs due to competition for transcriptional and translational machinery. Engineered \textit{E.~coli} systems show that incomplete repression results in expression levels reaching only 70--90\% of theoretical maximum under realistic conditions~\cite{makrides1996strategies}, consistent with the unscaled logistic's natural expression ceiling.
\end{itemize}

In summary, the logistic repression function---scaled or unscaled---provides a mathematically elegant, biologically realistic, and computationally robust alternative to the Hill function: \textbf{(i) modeling realism}, since the unscaled form captures resource-limited and stochastic environments where $f^-(0) < 1$; \textbf{(ii) theoretical consistency}, since the scaled variant matches idealized Hill behavior while preserving the logistic's tractability; \textbf{(iii) parameter tunability}, since $\lambda$ and $\theta$ are independent (in contrast to the Hill function's coupled $n/(4\theta)$ slope); \textbf{(iv) computational advantages}, since the exponential form avoids the power-law singularities of $x^n$ near $x = 0$ for non-integer $n$; and \textbf{(v) control-theoretic compatibility}, since smooth derivatives of all orders enable feedback linearization, Lyapunov-based stability analysis, and observer design.

\section{General Logistic GRN Framework, Boolean Logic, and Theoretical Foundations}
\label{sec:advanced_analysis}

Our core proposal is to replace classical Hill-type sigmoidal functions with logistic functions, which furnish substantially improved analytical tractability across all dimensions of modeling: mathematical analysis becomes more transparent and rigorous, parameter estimation from experimental measurements becomes more robust, and the synthesis of feedback control systems and state observers becomes practical where previously it was intractable. The framework applies across a diverse spectrum of biological domains: gene regulatory networks, cellular proliferation dynamics, hematopoietic tissue formation, epidemiological disease transmission, immunological responses, pharmacodynamic drug action, neural circuit dynamics, tumor growth progression, and microbial bioprocess engineering, among many others. Beyond theoretical advances, this framework has immediate practical implications for synthetic biology and biotechnology: logistic-based models enable robust design of genetic circuits with predictable behavior under noisy cellular conditions, facilitate real-time control of engineered biosynthetic pathways through model predictive control strategies, and provide computationally efficient frameworks for genome-scale metabolic engineering, where thousands of regulatory interactions must be simulated simultaneously.

A distinguishing feature of our modeling philosophy is the \emph{explicit differentiation} between increasing and decreasing logistic functions: we deploy each form precisely where it is biologically appropriate, thereby maintaining clear regulatory interpretability. This stands in contrast to alternative formulations that unify all regulatory interactions through a single increasing functional form with signed weights~\cite{samuilik2022mathematical}, potentially obscuring biological meaning. Our product-of-logistics formulation preserves distinct functional forms for activation and repression, maintains regulator-specific thresholds that enable modular experimental characterization, and avoids parameter coupling in which threshold and weight adjustments compensate for each other during fitting.

\subsection{General Logistic-Based Gene Regulatory Network Model}

Let $\mathbf{x}(t) = (x_1(t), \ldots, x_n(t))^\top \in \mathbb{R}^n$ denote the state vector of gene expression levels in a network comprising $n$ genes, where $x_i(t)$ represents the concentration of the protein product of gene $i$ at time $t$. Each gene's dynamics obey a balance equation between synthesis and degradation:
\begin{equation}
\dot{x}_i(t) = \kappa_i\, f_i\!\big(\mathbf{x}(t)\big) - \gamma_i\, x_i(t), \qquad i = 1, \ldots, n,
\label{eq:logistic_system_general}
\end{equation}
where $\kappa_i > 0$ is the maximal production rate, $\gamma_i > 0$ is the degradation rate, and $f_i : \mathbb{R}^n \to (0,1)$ is the regulatory function that synthesizes the influences from all relevant activators and repressors acting on gene $i$.

\subsubsection{The product-of-logistics regulatory function}
We model parallel regulation, wherein multiple transcription factors simultaneously exert independent regulatory effects, some activating and others repressing the target gene's promoter. In this architecture, the regulatory function takes the product structure:
\begin{equation}
f_i(\mathbf{x}) = \prod_{m=1}^{M_i} g_{i,m}\!\big(x_{j(i,m)}\big),
\label{eq:parallel_regulation_product}
\end{equation}
where $M_i$ denotes the total number of regulatory inputs to gene $i$, and each factor is given by the logistic function:
\begin{equation}
g_{i,m}(s) = \frac{1}{1 + \exp\!\big[-\lambda_{i,m}\,\sigma_{i,m}\,(s - \theta_{i,m})\big]}.
\label{eq:logistic_factor}
\end{equation}
The parameters are:
\begin{itemize}
    \item $\lambda_{i,m} > 0$ controls the steepness of the regulatory response (analogous to the Hill coefficient $n$, matched via $\lambda_{i,m} = n/\theta_{i,m}$),
    \item $\theta_{i,m} \in \mathbb{R}_{>0}$ represents the regulatory threshold concentration (the concentration at which the regulatory effect is half-maximal),
    \item $\sigma_{i,m} \in \{+1,-1\}$ encodes the regulatory sign: $\sigma_{i,m} = +1$ for activation and $\sigma_{i,m} = -1$ for repression,
    \item the index map $j(i,m) \in \{1, \ldots, n\}$ connects the $m$-th regulatory input of gene $i$ to its corresponding state variable.
\end{itemize}
When $\sigma_{i,m} = +1$ (activation), \eqref{eq:logistic_factor} is the increasing logistic function $g_{i,m}(s) = f^+(s, \theta_{i,m}, \lambda_{i,m})$; when $\sigma_{i,m} = -1$ (repression), it is the decreasing logistic function $g_{i,m}(s) = f^-(s, \theta_{i,m}, \lambda_{i,m})$. Substituting \eqref{eq:parallel_regulation_product}--\eqref{eq:logistic_factor} into \eqref{eq:logistic_system_general}, the complete dynamical system is:
\begin{equation}
\dot{x}_i = \kappa_i \prod_{m=1}^{M_i} \frac{1}{1 + \exp\!\big[-\lambda_{i,m}\,\sigma_{i,m}\,(x_{j(i,m)} - \theta_{i,m})\big]} - \gamma_i\, x_i, \qquad i = 1, \ldots, n.
\label{eq:logistic_system_thm}
\end{equation}

The product structure~\eqref{eq:parallel_regulation_product} enjoys several desirable mathematical and biological properties:
\begin{enumerate}
    \item \textbf{Normalization:} Since each factor satisfies $g_{i,m} \in (0,1)$, the product satisfies $f_i(\mathbf{x}) \in (0,1)$ for all $\mathbf{x} \in \mathbb{R}^n$, ensuring that production rates are naturally bounded between $0$ and $\kappa_i$.
    \item \textbf{AND logic:} The product form models the requirement for the simultaneous satisfaction of all regulatory conditions (cooperative co-regulation), as each factor must be close to unity for $f_i$ to be large. This corresponds to a continuous relaxation of Boolean AND-gate logic.
    \item \textbf{Modular parameterization:} Each regulator $j(i,m)$ has its own independent threshold $\theta_{i,m}$ and steepness $\lambda_{i,m}$, directly determinable from dose-response curves for that particular interaction without compensatory adjustments.
    \item \textbf{Reduction to single-input case:} When $M_i = 1$, the model reduces to a single logistic factor, recovering the non-cooperative case.
    \item \textbf{Basal expression:} Since $g_{i,m}(s) > 0$ for all $s$, the product $f_i(\mathbf{x}) > 0$ for all $\mathbf{x}$, naturally capturing leaky promoter activity without artificial offsets.
\end{enumerate}

The parameters $\lambda_{i,m}$ may be set uniformly ($\lambda_{i,m} = \lambda$ for all $i,m$) for simplicity, or varied to capture heterogeneous regulatory sensitivities across different interactions. This systematic parameterization framework links logistic parameters $(\theta_{i,m}, \lambda_{i,m})$ to biological observables: dissociation constants $K_d$ measured via electrophoretic mobility shift assays, cooperativity coefficients extracted from Hill plots, half-maximal effective concentrations EC$_{50}$ from fluorescence reporter assays, and maximal expression levels quantified via qPCR or RNA-seq.

\subsection{Extension to Boolean Regulatory Logic}

The general framework~\eqref{eq:logistic_system_general} naturally accommodates the full range of \emph{logical rules} governing how multiple regulators collectively orchestrate gene expression, following the Boolean network formalism pioneered by Albert and Othmer~\cite{albert2003topology} for the \textit{Drosophila} segment polarity network. We denote by $\Phi : \mathcal{B} \to [0,1]$ the map that translates a Boolean regulatory rule into its continuous logistic approximation. The map $\Phi$ is defined recursively on the structure of the Boolean rule, with~\eqref{eq:parallel_regulation_product}--\eqref{eq:logistic_factor} providing the building blocks:

\begin{itemize}
    \item \textbf{Constant \textsc{False}} (gene permanently silenced):
    $\Phi(\textsc{False}) \equiv 0$, so~\eqref{eq:logistic_system_general} reduces to pure exponential decay $\dot{x}_i = -\gamma_i x_i$.

    \item \textbf{Constant \textsc{True}} (constitutive expression):
    $\Phi(\textsc{True}) \equiv 1$, giving monotone convergence to $x_i^* = \kappa_i/\gamma_i$.

    \item \textbf{Positive literal} (activation): variable $x_j$ is mapped to
    $\Phi(x_j) = g_{i,m}(x_j)$ with $\sigma_{i,m} = +1$.

    \item \textbf{NOT logic} (repression): negative literal $\lnot x_j$ is mapped to
    \begin{equation}
    \Phi(\lnot x_j) = g_{i,m}(x_j)\big|_{\sigma=-1} = f^-(x_j, \theta_j, \lambda) = 1 - f^+(x_j, \theta_j, \lambda),
    \label{eq:not_logic}
    \end{equation}
    consistent with the complement structure of Boolean negation (Section~\ref{sec:logistic_foundations}).

    \item \textbf{AND logic} (cooperative co-regulation): a conjunction $C(\mathbf{x}) = x_{i_1} \wedge \cdots \wedge x_{i_k}$ is mapped to the product of the corresponding logistic terms, which is exactly~\eqref{eq:parallel_regulation_product}:
    \[
    \Phi\!\left(x_{i_1} \wedge \cdots \wedge x_{i_k}\right) = \prod_{l=1}^{k} g_{i,m_l}(x_{i_l}),
    \]
    modeling the requirement for simultaneous satisfaction of all regulatory conditions.

    \item \textbf{OR logic} (independent activation): a disjunction $C_1(\mathbf{x}) \vee C_2(\mathbf{x})$ is mapped via the De~Morgan formula. Since $A \vee B = \lnot(\lnot A \wedge \lnot B)$, and applying \eqref{eq:not_logic} and the AND product rule:
    \begin{equation}
    \Phi(C_1 \vee C_2) = 1 - \big(1 - \Phi(C_1)\big)\big(1 - \Phi(C_2)\big).
    \label{eq:or_elementary}
    \end{equation}
    When each clause $C_k$ is a positive literal $x_{j_k}$, the complement $1-\Phi(C_k)$ coincides with the decreasing logistic $f^-(x_{j_k},\theta_{j_k},\lambda_{j_k})$ via \eqref{eq:not_logic}; for compound clauses, $1-\Phi(C_k)$ denotes the general complement of $\Phi(C_k)$, not necessarily a single decreasing logistic function.
    This generalizes recursively to $m$ independent clauses via the De~Morgan product formula:
    \begin{equation}
    \Phi\!\left(\bigvee_{k=1}^{m} C_k(\mathbf{x})\right) = 1 - \prod_{k=1}^{m}\big(1 - \Phi(C_k(\mathbf{x}))\big),
    \label{eq:demorgan_0}
    \end{equation}
    which ensures $\Phi \in [0,1]$ regardless of the number of independent regulatory pathways, in contrast to naive additive translations where $\sum_k \Phi(C_k)$ can exceed unity when multiple activators are simultaneously active, inflating production rates to $m\kappa_i$ and violating the biological bound $x_i^* \leq \kappa_i/\gamma_i$.
\end{itemize}

\begin{remark}[De Morgan representation of the OR gate]
\label{rem:demorgan}
The formula~\eqref{eq:demorgan_0} makes the De~Morgan structure manifest: Boolean OR is the complement of the conjunction of the complements. For two positive literals $x_1, x_2$, equation~\eqref{eq:or_elementary} reduces to $1-(1-f^+(x_1,\cdot))(1-f^+(x_2,\cdot)) = 1 - f^-(x_1)f^-(x_2)$, and the same decreasing logistic $f^-$ that encodes NOT logic~\eqref{eq:not_logic} thus also builds OR logic multiplicatively for positive literals. Wittmann \emph{et al.}~\cite{wittmann2009transforming} obtain the algebraically equivalent two-input instance $\Phi(x_1 \vee x_2) = x_1 + x_2 - x_1 x_2$ via multivariate polynomial interpolation of the Boolean OR gate, but do not identify it as a product of complemented logistic functions, nor state the general $m$-clause recursive product form~\eqref{eq:demorgan_0}. A concrete worked-out application of~\eqref{eq:demorgan_0} to the 11-gene Traynard mammalian cell-cycle Boolean network~\cite{Traynard2016}, including automatic translation to a continuous logistic ODE system and numerical integration with all variables remaining strictly bounded and non-negative, is given in the companion paper~\cite{belgacem2026logistic2}. The product representation~\eqref{eq:parallel_regulation_product} thus provides a unified encoding: AND logic corresponds directly to the product $\prod_m g_{i,m}$, while OR logic is built from this same structure via complementation~\eqref{eq:demorgan_0}.
\end{remark}

\subsection{Comparison with Weighted-Sum Logistic Formulations}
\label{sec:weighted_sum_comparison}

A widely used alternative aggregates all regulatory inputs into a single weighted sum passed through one increasing logistic function per gene~\cite{samuilik2022mathematical}:
\begin{equation}
\dot{x}_i = \kappa_i \cdot \frac{1}{1 + e^{-\mu_i\!\left(\sum_{j=1}^n w_{ij}x_j - \theta_i\right)}} - \gamma_i x_i, \qquad w_{ij} \in \mathbb{R}.
\label{eq:samuilik}
\end{equation}
While this unified form has attractive mathematical properties and preserves the $[0,1]$ bound, it carries structural costs that are not merely aesthetic.

\subsubsection{Critical-point pathology for repressors}
In \eqref{eq:samuilik}, repression is encoded by a negative weight $w_{ij} < 0$ inside an \emph{increasing} sigmoid. The midpoint of the sigmoid (where the output equals $1/2$ and the slope is steepest) satisfies $\sum_j w_{ij}x_j = \theta_i$. For a pure repressor with $w_{ij} < 0$, this places the critical point at $x_j = \theta_i/w_{ij} < 0$, strictly outside the biologically admissible domain $x_j \geq 0$. The repression function, therefore, remains nearly flat and close to zero throughout the physical domain, rendering it biologically inert. This is a structural consequence of encoding direction through weight signs inside a single increasing sigmoid; it is not correctable by recalibration.

By contrast, in the product-of-logistics formulation~\eqref{eq:logistic_system_thm}, repression is encoded through a \emph{decreasing} logistic factor $f^-(x_j, \theta_{i,m}, \lambda_{i,m})$, whose midpoint is at $x_j = \theta_{i,m} > 0$, directly interpretable as the $\text{IC}_{50}$ of the repressor.

\subsubsection{Threshold-scaling pathology}
The shared threshold in \eqref{eq:samuilik} is commonly prescribed as $\theta_i = \sum_j w_{ij}/2$~\cite{samuilik2022mathematical}. This threshold grows linearly with the number of activators (positive weights), becomes negative for pure-inhibitor gates, and loses all correspondence with any measurable molecular quantity. For an AND gate with one activator ($w_{i1} = +1$) and one repressor ($w_{i2} = -1$), the prescription yields $\theta_i = 0$, a threshold that corresponds to no measurable binding affinity or half-maximal concentration.

In the product-of-logistics framework, each regulator has its own independent threshold, $\theta_{i,m}$, directly determinable from single-regulator dose-response data, independent of all other regulators and of network size.

\subsubsection{Formal equivalence between weighted and fixed-weight product-of-logistics}
Within the product-of-logistics framework, incorporating explicit positive interaction weights $w_{ij} > 0$ is formally equivalent to the fixed-weight formulation after a parameter rescaling. For an activator with weight $w_{ij} > 0$:
\[
f^+(w_{ij}x_j, \theta_{ij}, \lambda_{ij}) = \frac{1}{1+e^{-\lambda_{ij}(w_{ij}x_j - \theta_{ij})}} = \frac{1}{1+e^{-\lambda'_{ij}(x_j - \theta'_{ij})}},
\]
where $\lambda'_{ij} = \lambda_{ij}w_{ij} > 0$ and $\theta'_{ij} = \theta_{ij}/w_{ij} > 0$. The effective threshold $\theta' = \theta/w > 0$ is always positive and biologically interpretable as $\text{EC}_{50}$. An identical rescaling holds for decreasing logistic repression terms. Consequently, the weighted and fixed-weight formulations produce identical trajectories after parameter estimation and differ only in how steepness and threshold information is distributed across parameters.

\subsection{Fundamental Mathematical Properties}

We now establish the rigorous theoretical foundations of system~\eqref{eq:logistic_system_thm}, exploiting the product structure~\eqref{eq:parallel_regulation_product}--\eqref{eq:logistic_factor} throughout the proofs. Theorem~\ref{thm:smoothness} below proves global existence, smoothness, and boundedness for the full multi-gene logistic GRN model; this complements prior global-dynamics results obtained for related smooth and reduced models---transcription--translation reductions~\cite{belgacem2018reduction,belgacem2013analysis,belgacem2013stability,belgacem2014stability}, concave gene-expression dynamics~\cite{belgacem2014mathematical}, and full-reversible enzymatic chains~\cite{belgacem2013global,belgacem2012full}---and lifts the pointwise logistic properties from Section~\ref{sec:logistic_foundations} to the network level with explicit, biologically interpretable constants.

\begin{theorem}[Global Existence, Smoothness, and Boundedness]
\label{thm:smoothness}
Consider the logistic gene regulatory network~\eqref{eq:logistic_system_thm} with regulatory functions~\eqref{eq:parallel_regulation_product}--\eqref{eq:logistic_factor}. The following properties hold for all parameter values $\kappa_i, \gamma_i, \lambda_{i,m} > 0$, $\theta_{i,m} \in \mathbb{R}$, and $\sigma_{i,m} \in \{+1,-1\}$:

\begin{enumerate}
\item \textbf{Infinite Smoothness.} The vector field $\mathbf{F} : \mathbb{R}^n \to \mathbb{R}^n$ defined by $F_i(\mathbf{x}) = \kappa_i f_i(\mathbf{x}) - \gamma_i x_i$ is of class $C^\infty$ on all of $\mathbb{R}^n$.

\item \textbf{Global Lipschitz Continuity and Unique Global Solutions.} The vector field $\mathbf{F}$ is globally Lipschitz continuous on $\mathbb{R}^n$ with explicit Lipschitz constant
\[
L_F \leq M := \max_{i=1,\ldots,n}\!\left(\kappa_i \sum_{j=1}^n L_i^j + \gamma_i\right), \quad L_i^j = \sum_{\substack{m=1 \\ j(i,m)=j}}^{M_i} \frac{\lambda_{i,m}}{4}.
\]
Consequently, for every initial condition $\mathbf{x}(0) \in \mathbb{R}^n$, there exists a unique solution $\mathbf{x}(t)$ defined for all $t \geq 0$.

\item \textbf{Positive Invariance and Uniform Boundedness.} The rectangular region
\[
\mathcal{B} := \prod_{i=1}^n \!\left[0,\, \frac{\kappa_i}{\gamma_i}\right] \subset \mathbb{R}^n
\]
is positively invariant under the flow of~\eqref{eq:logistic_system_thm}: if $\mathbf{x}(0) \in \mathcal{B}$, then $\mathbf{x}(t) \in \mathcal{B}$ for all $t \geq 0$, so each component satisfies $0 \leq x_i(t) \leq \kappa_i/\gamma_i$. Moreover, the positive orthant $\mathbb{R}^n_+$ is also positively invariant: if $\mathbf{x}(0) \in \mathbb{R}^n_+$, then $x_i(t) \geq 0$ for all $t \geq 0$ and $i = 1,\ldots,n$.
\end{enumerate}
\end{theorem}

Before presenting the proof, we recall the fundamental existence and uniqueness result that underlies our analysis.

\begin{theorem}[Picard--Lindel\"{o}f~\cite{picard1890memoire,lindelof1894application}]
\label{thm:picard_lindelof}
Consider $\dot{\mathbf{x}}(t) = \mathbf{F}(\mathbf{x}(t))$, $\mathbf{x}(t_0) = \mathbf{x}_0$, with $\mathbf{F} : D \subseteq \mathbb{R}^n \to \mathbb{R}^n$ continuous on an open set $D$.
\begin{enumerate}
    \item \textbf{Local existence and uniqueness:} If $\mathbf{F}$ is locally Lipschitz in a neighborhood of $\mathbf{x}_0$, there exists $\delta > 0$ and a unique solution on $[t_0, t_0+\delta]$.
    \item \textbf{Global existence and uniqueness:} If $\mathbf{F}$ is globally Lipschitz on $\mathbb{R}^n$, then for every $\mathbf{x}_0 \in \mathbb{R}^n$, the unique solution exists for all $t \geq t_0$.
\end{enumerate}
\end{theorem}

\begin{proof}[Proof of Theorem~\ref{thm:smoothness}]
The proof proceeds in three parts corresponding to the three claims, exploiting the product structure~\eqref{eq:parallel_regulation_product}--\eqref{eq:logistic_factor} throughout.

\medskip
\noindent\textbf{Part (i): Infinite Smoothness.}

Each scalar logistic factor in~\eqref{eq:logistic_factor} is:
\[
g_{i,m}(s) = \big(1 + e^{-\lambda_{i,m}\,\sigma_{i,m}\,(s-\theta_{i,m})}\big)^{-1}.
\]
This is a composition of three globally analytic functions: (a) the linear map $s \mapsto -\lambda_{i,m}\sigma_{i,m}(s-\theta_{i,m})$, which is $C^\infty$ on $\mathbb{R}$; (b) the exponential $u \mapsto e^u$, which is $C^\infty$ on $\mathbb{R}$ with values in $(0,\infty)$; and (c) the inversion $v \mapsto (1+v)^{-1}$, which is $C^\infty$ on $\{v \in \mathbb{R} : v \neq -1\}$, in particular on $(0,\infty)$. Since $1 + e^{-\lambda_{i,m}\sigma_{i,m}(s-\theta_{i,m})} > 1 > 0$ for all $s \in \mathbb{R}$, the denominator never vanishes, so $g_{i,m} \in C^\infty(\mathbb{R})$.

The regulatory function~\eqref{eq:parallel_regulation_product},
\[
f_i(\mathbf{x}) = \prod_{m=1}^{M_i} g_{i,m}\!\big(x_{j(i,m)}\big),
\]
is a finite product of $C^\infty$ functions of the coordinates $x_{j(i,m)}$, hence $f_i \in C^\infty(\mathbb{R}^n)$. Finally, $F_i(\mathbf{x}) = \kappa_i f_i(\mathbf{x}) - \gamma_i x_i$ is a linear combination of $C^\infty$ functions, so $F_i \in C^\infty(\mathbb{R}^n)$ and the entire vector field $\mathbf{F}$ is $C^\infty$ on $\mathbb{R}^n$.

\medskip
\noindent\textbf{Part (ii): Global Lipschitz Continuity.}

We bound the Jacobian $D\mathbf{F}(\mathbf{x})$ uniformly over all $\mathbf{x} \in \mathbb{R}^n$ by exploiting the uniform derivative bound established in Section~\ref{sec:logistic_foundations}.

\textit{Step 1: Derivative of a single logistic factor.}
For any factor $g_{i,m}$ from~\eqref{eq:logistic_factor}, differentiating using the chain rule:
\[
g_{i,m}'(s) = \lambda_{i,m}\,\sigma_{i,m}\, g_{i,m}(s)\big(1 - g_{i,m}(s)\big).
\]
Since $g_{i,m}(s) \in (0,1)$ for all $s \in \mathbb{R}$, the product $g_{i,m}(1 - g_{i,m}) \leq 1/4$ (with equality when $g_{i,m} = 1/2$, i.e., at $s = \theta_{i,m}$), and therefore:
\begin{equation}
|g_{i,m}'(s)| = \lambda_{i,m}\,g_{i,m}(s)\big(1-g_{i,m}(s)\big) \leq \frac{\lambda_{i,m}}{4} \quad \text{for all } s \in \mathbb{R}.
\label{eq:bound_gim}
\end{equation}

\textit{Step 2: Partial derivatives of the product regulatory function.}
Applying the product rule to~\eqref{eq:parallel_regulation_product} with respect to $x_j$:
\[
\frac{\partial f_i}{\partial x_j}(\mathbf{x}) = \sum_{\substack{m=1\\j(i,m)=j}}^{M_i} g_{i,m}'\!\big(x_{j(i,m)}\big) \prod_{\substack{q=1\\q\neq m}}^{M_i} g_{i,q}\!\big(x_{j(i,q)}\big).
\]
Since each factor satisfies $g_{i,q} \in (0,1)$, the product over $q \neq m$ is bounded by $1$. Applying~\eqref{eq:bound_gim}:
\begin{equation}
\left|\frac{\partial f_i}{\partial x_j}(\mathbf{x})\right| \leq \sum_{\substack{m=1\\j(i,m)=j}}^{M_i} \frac{\lambda_{i,m}}{4} =: L_i^j \quad \text{for all } \mathbf{x} \in \mathbb{R}^n.
\label{eq:bound_partial_fi}
\end{equation}

\textit{Step 3: Uniform bound on the Jacobian of $\mathbf{F}$.}
The entries of the Jacobian $D\mathbf{F}(\mathbf{x})$ are:
\[
\frac{\partial F_i}{\partial x_j}(\mathbf{x}) = \kappa_i \frac{\partial f_i}{\partial x_j}(\mathbf{x}) - \gamma_i\,\delta_{ij},
\]
where $\delta_{ij}$ is the Kronecker delta. Using~\eqref{eq:bound_partial_fi}:
\[
\left|\frac{\partial F_i}{\partial x_j}(\mathbf{x})\right| \leq \begin{cases} \kappa_i\, L_i^j & j \neq i, \\ \kappa_i\, L_i^i + \gamma_i & j = i. \end{cases}
\]
Taking the row-sum maximum over all $i$ and $j$, we obtain the uniform bound on the operator norm of the Jacobian (the $\ell^\infty$ operator norm equals the maximum row sum of absolute entries):
\[
\|D\mathbf{F}(\mathbf{x})\|_\infty \leq M := \max_{i=1,\ldots,n}\!\left(\kappa_i \sum_{j=1}^n L_i^j + \gamma_i\right) < \infty,
\]
for all $\mathbf{x} \in \mathbb{R}^n$, where the finiteness of $M$ follows from the finiteness of $n$, $M_i$, and the parameters. By the mean value theorem applied componentwise, for any $\mathbf{x}, \mathbf{y} \in \mathbb{R}^n$:
\[
\|\mathbf{F}(\mathbf{x}) - \mathbf{F}(\mathbf{y})\|_\infty \leq M \|\mathbf{x} - \mathbf{y}\|_\infty,
\]
so $\mathbf{F}$ is globally Lipschitz with constant $L_F \leq M$ in the $\ell^\infty$ norm (and, since all norms on $\mathbb{R}^n$ are equivalent, in every norm with a constant of comparable order). Global existence and uniqueness of solutions follow from Theorem~\ref{thm:picard_lindelof}, Part~2.

\begin{remark}[Lipschitz continuity of the Jacobian]
\label{rem:lipschitz_jacobian}
The Jacobian $D\mathbf{F}$ is itself globally Lipschitz continuous (i.e., the Hessian $D^2\mathbf{F}$ is uniformly bounded). This follows from bounding the second-order partial derivatives of $f_i$ from~\eqref{eq:parallel_regulation_product}. For $j \neq k$:
\[
\frac{\partial^2 f_i}{\partial x_k\,\partial x_j} = \sum_{\substack{m:\\j(i,m)=j}} \sum_{\substack{p:\\j(i,p)=k}} g_{i,m}'\, g_{i,p}' \prod_{\substack{q \neq m,p}} g_{i,q},
\]
which, by~\eqref{eq:bound_gim} and $g_{i,q} < 1$, is bounded in absolute value by $(\Lambda_i/4)^2$ where $\Lambda_i = \max_m \lambda_{i,m}$. For $j = k$, the second-order terms also involve $g_{i,m}''$. Differentiating the self-referential derivative $g_{i,m}'(s) = \lambda_{i,m}\sigma_{i,m}\,g_{i,m}(1-g_{i,m})$ once more (see Section~\ref{sec:logistic_foundations}):
\[
|g_{i,m}''(s)| = \lambda_{i,m}^2\,|g_{i,m}(1-g_{i,m})(1-2g_{i,m})| \leq \lambda_{i,m}^2\,\rho,
\]
where $\rho = \max_{g \in (0,1)} |g(1-g)(1-2g)| = \sqrt{3}/18 \approx 0.096$. Hence, all second-order partials are uniformly bounded, so $\|D^2\mathbf{F}(\mathbf{x})\| \leq K < \infty$ for the explicit constant:
\[
K = \kappa_{\max} \cdot \max_i\!\left(|\mathcal{M}_i|^2 \frac{\Lambda_i^2}{16} + |\mathcal{M}_i|\,\Lambda_i^2\,\rho\right), \quad \kappa_{\max} = \max_i \kappa_i,
\]
where $|\mathcal{M}_i| = M_i$ is the number of regulatory inputs. This establishes that $D\mathbf{F}$ is globally Lipschitz with constant $L_{DF} = K$, which implies, among other properties, that sensitivity equations are globally well-posed and Gauss-Newton-type parameter estimation algorithms converge reliably.
\end{remark}

\medskip
\noindent\textbf{Part (iii): Positive Invariance and Uniform Boundedness.}

Since each logistic factor satisfies $g_{i,m}(s) \in (0,1)$ for all $s \in \mathbb{R}$, the product~\eqref{eq:parallel_regulation_product} satisfies:
\[
0 < f_i(\mathbf{x}) < 1 \quad \text{for all } \mathbf{x} \in \mathbb{R}^n,
\]
so $0 < \kappa_i f_i(\mathbf{x}) < \kappa_i$ for all $\mathbf{x}$.

\textit{Nonnegativity.} Suppose at some time $t_0$ a component satisfies $x_i(t_0) = 0$. Then:
\[
\dot{x}_i(t_0) = \kappa_i f_i\!\big(\mathbf{x}(t_0)\big) - \gamma_i \cdot 0 = \kappa_i f_i\!\big(\mathbf{x}(t_0)\big) > 0.
\]
The vector field points strictly into $\mathbb{R}_+^n$ at every boundary face $\{x_i = 0\}$, so no component can become negative.

\textit{Upper bound.} Suppose at some time $t_1$ a component satisfies $x_i(t_1) = \kappa_i/\gamma_i$. Then:
\[
\dot{x}_i(t_1) = \kappa_i f_i\!\big(\mathbf{x}(t_1)\big) - \gamma_i \cdot \frac{\kappa_i}{\gamma_i} = \kappa_i\!\left(f_i\!\big(\mathbf{x}(t_1)\big) - 1\right) < 0,
\]
since $f_i < 1$ strictly. The vector field points strictly inward at every ceiling face $\{x_i = \kappa_i/\gamma_i\}$.

By these two inequalities, the vector field $\mathbf{F}$ points strictly inward on every face of the box $\mathcal{B} = \prod_{i=1}^n [0, \kappa_i/\gamma_i]$. By standard invariance arguments (see, e.g., \cite{brezis2010functional}), $\mathcal{B}$ is positively invariant: any trajectory starting in $\mathcal{B}$ remains in $\mathcal{B}$ for all $t \geq 0$. The nonnegativity argument above also shows that $\mathbb{R}^n_+$ is positively invariant.

This completes the proof.
\end{proof}

\begin{remark}[Implications of the Lipschitz structure for applications]
\label{rem:lipschitz_implications}
The global Lipschitz continuity of $\mathbf{F}$ and its Jacobian $D\mathbf{F}$ established above have several concrete consequences:
\begin{enumerate}
    \item \textbf{Numerical stability:} Standard ODE solvers (Runge-Kutta, Adams-Bashforth) exhibit uniform step-size control without the refinements forced by Hill function singularities near $x = 0$.
    \item \textbf{Sensitivity analysis:} The variational equation $\dot{Z} = D\mathbf{F}(\mathbf{x}(t))Z$ is globally well-posed; perturbations in initial conditions grow at most as $e^{L_F t}$.
    \item \textbf{Observer design:} High-gain observers with convergence rate $e^{-(L-L_F)t}$ can be synthesized using $L_F = M$ as an explicit, computable gain bound.
    \item \textbf{Parameter estimation:} Gauss-Newton and Levenberg-Marquardt algorithms converge reliably, as the Fisher information matrix is well-conditioned (bounded gradient norms $|\partial g_{i,m}/\partial s| \leq \lambda_{i,m}/4$).
\end{enumerate}
These properties contrast sharply with Hill-function-based models, where non-integer exponents cause $|h^{+\prime}(x)| \to \infty$ as $x \to 0^+$, making $L_F = \infty$ and invalidating all of the above guarantees.
\end{remark}

\subsection{Numerical Instability of Hill Functions as a Unifying Diagnosis}
\label{sec:hill_pathologies}

Section~\ref{sec:practical_consequences} traced the various Hill-function pathologies into their analytical failure modes; here we collect their numerical-simulation consequences into a single forensic narrative. The various failure modes of Hill-function-based ODE systems are not independent accidents; they are all manifestations of a single underlying cause: the ODE system derived from Hill functions with real-valued $n$ is \emph{numerically unstable} in a precise sense.

A numerical integration method applied to $\dot{\mathbf{x}} = \mathbf{F}(\mathbf{x})$, $\mathbf{x}(0) = \mathbf{x}_0$ is locally stable only when $\mathbf{F}$ is Lipschitz-continuous in a neighborhood of the trajectory, so that the classical theory of ODEs guarantees existence, uniqueness, and continuous dependence on initial data~\cite{hairer1993solving}. As established in Section~\ref{sec:numerical_integration}, the Hill vector field with non-integer $n$ fails this requirement on a neighborhood of the positive-orthant boundary in two distinct, compounding ways. For $0 < n < 1$, the first partial derivative of $F_i$ with respect to $x_i$ diverges as $x_i \to 0^+$, so $\mathbf{F}$ is not locally Lipschitz at the boundary; the Lipschitz constant $L = \sup\|\nabla\mathbf{F}\|$ is infinite near the origin and the standard error bound $\|\mathbf{e}(t)\| \leq \|\mathbf{e}(0)\| e^{Lt}$ provides no useful guarantee. For $n > 1$, $\mathbf{F}$ is locally Lipschitz but only $C^{\lfloor n \rfloor}$: derivatives of order greater than $\lfloor n \rfloor$ diverge as $x_i \to 0^+$, so standard convergence theorems for Runge--Kutta methods of order $p > \lfloor n \rfloor$ cease to apply near the boundary, and adaptive solvers that estimate higher derivatives produce spuriously large local error estimates and are forced to reduce step size dramatically. A second, complementary failure occurs in both regimes whenever a trajectory overshoots zero: by floating-point rounding of order $10^{-15}$, the expression $x_i^n$ becomes complex, and the right-hand side ceases to be a well-defined real-valued vector field altogether.

The consequence is that, \emph{in the presence of real-valued $n$, the Hill-function ODE system does not admit a Lipschitz-continuous right-hand side with uniformly bounded derivatives of all orders on any neighborhood of the positive-orthant boundary}, the very region that low-expression states and transient dynamics necessarily explore. Standard convergence theorems for Runge--Kutta and multistep methods do not apply; the computed solution has no guaranteed accuracy bound, even in exact arithmetic; and in finite-precision arithmetic, the simulation will fail or produce incorrect results with near-certainty for any trajectory passing through a neighborhood of zero, which is virtually inevitable for randomly initialized high-dimensional systems.

This is not a deficiency of any particular solver implementation. It is a structural property of the Hill function itself when $n \notin \mathbb{N}$, and it persists regardless of the programming language, the tolerance settings, or the integration algorithm used.

\subsubsection*{The cascade of failures}

Once $n \notin \mathbb{N}$, the instability manifests through the following chain of observable consequences.

\subsubsection{Loss of smoothness at $x=0$}
For $n \in (k, k+1)$, the $(k{+}1)$-th derivative of $h^+(x;\theta,n)$ diverges as $x \to 0^+$. The ODE vector field is therefore only $C^k$, not $C^\infty$. Runge--Kutta methods of order $p > k$ lose their theoretical convergence rate near the origin; adaptive solvers that internally estimate higher-order derivatives generate spuriously large error estimates and are forced to reduce the step size towards zero or terminate prematurely.

\subsubsection{Complex-valued arithmetic upon any negative overshoot}
Standard arithmetic defines $x^n$ for non-integer $n$ and $x < 0$ as a complex number via $x^n = |x|^n e^{i\pi n}$. Any adaptive solver that allows a trajectory to cross zero---even by a rounding error of order $10^{-15}$---will immediately corrupt the entire right-hand side of the ODE system with imaginary parts, rendering it undefined as a real vector. This failure occurs identically in Python/SciPy, MATLAB/\texttt{ode45}, Julia/DifferentialEquations, Mathematica/\texttt{NDSolve}, and any other general-purpose solver operating under standard floating-point conventions.

\subsubsection{Trajectory corruption and premature domain truncation}
Once any $x_i$ overshoots to a negative value, the complex-valued right-hand side corrupts every subsequent integration step. Critically, the corruption begins \emph{silently}: imaginary parts can be many orders of magnitude smaller than the state variables themselves (with magnitude on the order of $|x_{\rm overshoot}|^n$, which for small overshoots and moderate $n$ can fall well below typical numerical noise thresholds) and may appear long before any solver warning fires. The solver continues building its numerical solution past the first warning until complex arithmetic fully overwhelms step-size control, producing an integration domain that terminates well before the intended final time $t_f$. The entire trajectory from the first contaminated step onward is a corrupted path, not a solution to the true biological ODE system, regardless of whether any given query time lies inside or outside the solver's domain. The companion paper~\cite{belgacem2026logistic2} reports a direct empirical demonstration of this failure mode on an 80-gene Boolean-derived ODE system with non-integer cooperativity $n = 3.509$, where the Hill solver entered silent complex-valued contamination at $t \approx 52.64$ and terminated near $t \approx 63$--$65$ out of an intended horizon $t_f = 200$, while the logistic formulation, with identical parameters and initial conditions, completed the full integration without a single warning and with all 80 state variables non-negative and bounded throughout.

\subsubsection{Silent failure at all observation times}
This is the most insidious failure mode. Solver warnings typically appear only in a message log, not in the returned solution object, and do not prevent the solver from returning a seemingly well-formed interpolating function. A modeler who does not read every solver message will query that object at observation times $t^* < t_{\rm end}$, receive numerically well-defined values, and proceed to analyze them as solutions to the biological ODE system---while they are evaluations of a corrupted path that has been tracking complex arithmetic since the first overshoot. When downstream code queries the object at $t^* > t_{\rm end}$, the solver silently returns polynomial extrapolations that include large unphysical values (negative concentrations, values orders of magnitude above biological bounds). The caller receives a number; there is no exception. In both cases, the pipeline produces results that are numerically well-defined but biologically meaningless, detectable only by reading every warning emitted during integration and at each query time.

\subsubsection*{Why standard workarounds are inadequate}

Practitioners sometimes attempt to circumvent these problems by clamping $x$ to $[0,\infty)$ at each step (introducing a non-smooth kink at $x = 0$ that creates a new discontinuity in the vector field), replacing $x^n$ by $|x|^n$ (which makes $h^+$ an even function near zero, destroying the monotone sigmoid), adding a small offset $\epsilon > 0$ to $x$ (removing the singularity but introducing an arbitrary parameter with no biological interpretation), rounding $n$ to the nearest integer (changing the fitted parameter value, potentially crossing a bifurcation boundary), or tightening solver tolerances (postponing but not preventing the failure). None of these strategies resolves the underlying issue; they replace one numerical difficulty with another.

By contrast, the logistic function~\eqref{eq:logistic_system_thm} is globally $C^\infty$ and globally Lipschitz with constant $L_F \leq M$ (Theorem~\ref{thm:smoothness}), so the standard stability theory applies everywhere on $\mathbb{R}^N$, including near and below zero. The logistic system is therefore numerically stable by construction.

\section{Conclusion}
\label{sec:conclusion}

This paper establishes logistic functions as mathematically principled and computationally robust alternatives to Hill functions in gene regulatory network modeling. The increasing Hill function $h^+(x,\theta,n) = x^n/(x^n + \theta^n)$ is replaced by $f^+(x,\theta,\lambda) = 1/(1 + e^{-\lambda(x-\theta)})$, and the decreasing Hill function $h^-(x,\theta,n) = \theta^n/(x^n + \theta^n)$ by $f^-(x,\theta,\lambda) = 1/(1 + e^{\lambda(x-\theta)})$, with the parameter-matching relationship $\lambda = n/\theta$ equating the slopes of the two functions at the half-maximal concentration $x = \theta$. All formal results are purely mathematical---proofs, inequalities, and structural analyses that hold independently of any particular dataset, solver, or parameter regime.

The central technical claim is that each analytical pathology of the Hill function with non-integer $n$ admits a precise resolution by a corresponding property of the logistic function. The derivative singularities at the origin---$h^{+\prime}(x) \to \infty$ as $x \to 0^+$ for $0 < n < 1$, with the $(\lfloor n \rfloor{+}1)$-th derivative diverging on $n \in (k, k{+}1)$---are replaced by the uniformly bounded derivative $|\partial f^\pm/\partial x| \leq \lambda/4$. The hypergeometric or beta-function antiderivatives required by the Hill family for general $n$ collapse to the elementary form $\int f^+(x)\,dx = \lambda^{-1}\ln(1+e^{\lambda(x-\theta)}) + C$. The multivalued fractional-power inversion $(h^+)^{-1}(y) = \theta(y/(1-y))^{1/n}$ is replaced by the single-valued, smooth logit inverse $\theta + \lambda^{-1}\ln(y/(1-y))$. The Hill small-$n$ expansion $h^+ \approx 1/2 + (n/4)\ln(x/\theta)$, which diverges logarithmically as $x \to 0^+$ and so cannot serve as a global linearization, is replaced by the analogous logistic small-$\lambda$ expansion $f^+ \approx 1/2 + (\lambda/4)(x-\theta)$, which is linear and globally bounded. Finally, the absorbing state $h^+(0) = 0$ is replaced by the nonzero basal output $f^+(0,\theta,\lambda) = 1/(1+e^{\lambda\theta}) > 0$, which captures the persistent leaky transcription documented across gene regulatory systems without invoking ad hoc additive offsets. These pointwise properties are then lifted to the network level by Theorem~\ref{thm:smoothness}, which guarantees globally unique, smooth, and uniformly bounded solutions for the product-of-logistics GRN model with explicit Lipschitz constant $L_F \leq M = \max_i(\kappa_i \sum_j L_i^j + \gamma_i)$, while the Boolean-logic extension via the De~Morgan product formula $\Phi(\bigvee_k C_k) = 1 - \prod_k(1-\Phi(C_k))$ produces logic gates that remain within $[0,1]$ by construction.

The functional identity $h^+(x,\theta,n) = \sigma(n\ln(x/\theta))$ shows that the Hill function is the logistic evaluated at the log-ratio $\ln(x/\theta)$ rather than the additive deviation $x-\theta$, but the logarithmic change of variables introduces a state-dependent factor $e^{-s_i}$ on the production side that prevents the corresponding ODE models from being equivalent. The two formulations therefore encode genuinely different biological hypotheses---multiplicative-increment versus additive-threshold sensitivity---and the structural advantages of the logistic framework ($C^\infty$ regularity, finite global Lipschitz constant, nonzero basal output, parameter decoupling) hold regardless of which sensitivity interpretation a given system requires.

Several directions remain open. Systematic computational benchmarking should quantify convergence rates, Fisher information conditioning, and wall-clock integration times for Hill versus logistic formulations on realistic biological datasets. Bifurcation analysis of canonical GRN motifs---toggle switches, repressilators, and multi-stable circuits---should yield explicit stability boundaries using the closed-form Taylor coefficients derived here, complementing recent work on chaos and bifurcations in piecewise-linear gene-network models~\cite{belgacem2025chaos,farcot2019chaos}. Closed-form identifiability analysis, optimal experimental design, and observer and feedback-controller synthesis for synthetic gene circuits~\cite{belgacem2020control,chambon2020qualitative} stand to benefit directly from the global Lipschitz and $C^\infty$ structure established here, as does the numerical simulation of delay-differential GRN models in which time-lagged regulation interacts with smoothness-dependent stability criteria. Regularization techniques for switching/Zeno hybrid systems~\cite{belgacem2019probabilistic} provide a natural bridge from piecewise-linear approximations to the smooth logistic formulation. Applying the Boolean-to-ODE translation of Section~\ref{sec:advanced_analysis} to published cell-cycle and differentiation networks, together with attractor analysis of the resulting ODE systems, constitutes a further rich direction in which the present mathematical foundation can be tested at biological scale.

\section*{Declarations}

\subsection*{Availability of data, materials, and code}
No experimental datasets were generated or analyzed during this study. All theoretical results, mathematical proofs, and derivations are fully self-contained within the manuscript, with explicit formulas enabling independent reproduction of the results. No numerical simulations or computational code were used in this article.

\subsection*{Competing interests}
The author declares that he has no competing interests.

\subsection*{Funding}
Not applicable. This research received no specific grant from any funding agency in the public, commercial, or not-for-profit sectors.

\subsection*{Ethics approval and consent to participate}
Not applicable.

\subsection*{Authors' contributions}
Not applicable (single author).

\subsection*{Acknowledgements}
Not applicable.

\bibliographystyle{sn-mathphys-num}
\bibliography{mybibfile_0}

\end{document}